\documentclass[english]{article}
\usepackage[T1]{fontenc}
\usepackage[latin9]{inputenc}
\usepackage{geometry}
\usepackage{float}
\usepackage{booktabs}
\usepackage{mathtools}
\usepackage{amsmath}
\usepackage{amsthm}
\usepackage{amssymb}
\usepackage{graphicx}

\makeatletter

\providecommand{\tabularnewline}{\\}
\floatstyle{ruled}
\newfloat{algorithm}{tbp}{loa}
\providecommand{\algorithmname}{Algorithm}
\floatname{algorithm}{\protect\algorithmname}

\numberwithin{figure}{section}
\numberwithin{equation}{section}
\newcommand{\lyxaddress}[1]{
	\par {\raggedright #1
	\vspace{1.4em}
	\noindent\par}
}
\theoremstyle{plain}
\newtheorem{thm}{\protect\theoremname}
\theoremstyle{plain}
\newtheorem{lem}[thm]{\protect\lemmaname}
\theoremstyle{plain}
\newtheorem{prop}[thm]{\protect\propositionname}
\theoremstyle{remark}
\newtheorem*{rem*}{\protect\remarkname}
\theoremstyle{plain}
\newtheorem{assumption}[thm]{\protect\assumptionname}

\usepackage{algorithm,algpseudocode}

\usepackage{tikz}
\usetikzlibrary{shapes.geometric, arrows.meta,arrows, positioning}
\tikzstyle{block} = [rectangle, draw, fill=blue!20, text width=5em, text centered, rounded corners, minimum height=2em]
\tikzstyle{line} = [draw, -latex']
\tikzstyle{plainline} = [draw]
\tikzset{
	pb/.style={draw=none, text centered, minimum height=3em},
	block/.style={rectangle, draw, fill=blue!20, text width=10em, font=\footnotesize,  text centered, rounded corners, minimum height=3em},
	lblock/.style={rectangle, draw, fill=blue!20, text width=15em, ,font=\footnotesize, text centered, rounded corners, minimum height=.8em},
	doublearrow/.style={draw, latex'-latex'}  
}
\tikzset{multiline node/.style={ 
		align=center,
		text width=5cm 
	},
	line/.style={
		draw, -Latex  
	},
}

\ifdefined\showcaptionsetup
 \PassOptionsToPackage{caption=false}{subfig}
\fi
\usepackage{subfig}
\makeatother

\usepackage{babel}
\providecommand{\assumptionname}{Assumption}
\providecommand{\lemmaname}{Lemma}
\providecommand{\propositionname}{Proposition}
\providecommand{\remarkname}{Remark}
\providecommand{\theoremname}{Theorem}

\begin{document}
\title{Data Assimilation Models for Computing Probability Distributions of
Complex Multiscale Systems}
\author{Di Qi\textsuperscript{a} and Jian-Guo Liu\textsuperscript{b} }
\maketitle

\lyxaddress{\textsuperscript{a}Department of Mathematics, Purdue University,
150 North University Street, West Lafayette, IN 47907, USA}

\lyxaddress{\textsuperscript{b}Department of Mathematics and Department of Physics,
Duke University, Durham, NC 27708, USA}
\begin{abstract}
We introduce a data assimilation strategy aimed at accurately
capturing key non-Gaussian structures in probability distributions
using a small ensemble size. A major challenge in statistical forecasting of nonlinearly coupled multiscale systems is mitigating the large errors that arise when computing  high-order statistical
moments. To address this issue, a high-order stochastic-statistical
modeling framework is proposed that integrates statistical data assimilation into
finite ensemble predictions. The method effectively reduces the
approximation errors in finite ensemble estimates of non-Gaussian distributions by employing
a filtering update step that incorporates observation data in leading
moments to refine the high-order statistical feedback. Explicit filter operators
are derived from intrinsic nonlinear coupling structures, allowing
straightforward numerical implementations. We demonstrate the performance
of the proposed method through extensive numerical experiments on
a prototype triad system. The triad system offers an instructive and
computationally manageable platform mimicking essential aspects of
nonlinear turbulent dynamics. The numerical results show that the
statistical  data assimilation algorithm consistently captures the
mean and covariance, as well as various non-Gaussian probability
distributions exhibited in different statistical regimes of the triad system. The modeling
framework can serve as a useful tool for efficient sampling and reliable
forecasting of complex probability distributions commonly encountered
in a wide variety of applications involving multiscale coupling and
nonlinear dynamics.
\end{abstract}

\section{Introduction\protect\label{sec:Introduction}}

Predicting the distinct statistical behaviors observed in nonlinear
dynamical systems involving multiple spatial and temporal scales remains
a fundamental challenge across various natural and engineering problems
\cite{frisch1995turbulence,majda2006nonlinear,pedlosky2013geophysical,palmer2019stochastic}.
One primary difficulty arises from accurately quantifying the multiscale
nonlinear interactions between the large-scale mean state and small-scale
stochastic fluctuations amplified by inherent instability. Such interactions
often lead to non-Gaussian probability distributions characterized
by high-order statistics and intermittent extreme events, driven by
the intricate multiscale coupling mechanism \cite{cousins2014quantification,tong2021extreme,qi2022anomalous}.
Developing efficient computational algorithms capable of capturing
these critical non-Gaussian probabilistic features remains a central
issue in practical applications \cite{leutbecher2008ensemble,surace2019avoid,gao2023transition}.
Ensemble-based methods combined with data assimilation strategies
\cite{leutbecher2008ensemble,law2015data,cleary2021calibrate,bach2023filtering}
have been successfully applied for recovering leading-order statistics
in linear dynamical systems from noisy and partial observations. However,
as nonlinear coupling effects become dominant, low-order approaches
such as Kalman filters often suffer inherent difficulties and fail
to capture the essential higher-order moment statistics \cite{majda2018strategies,qi2023high,qi2023data}.
As a result, accurate and efficient methods for quantification and
prediction of these high-order statistics and the associated non-Gaussian
probability distributions are still needed for reliable forecasting
of the complex phenomena.

We consider the nonlinear statistical forecast problem formulated
as the following general stochastic dynamical equation (SDE) \cite{majda2018strategies,qi2024coupled}
describing the uncertainty evolution of the random state $u\in\mathbb{R}^{d}$
starting from random initial state $u\left(0;\omega\right)\sim\mu_{0}$
and driven by stochastic forcing and nonlinear interactions 
\begin{equation}
\frac{\mathrm{d}u}{\mathrm{d}t}=\Lambda u+B\left(u,u\right)+F\left(t\right)+\sigma\left(t\right)\dot{W}\left(t;\omega\right).\label{eq:model_general}
\end{equation}
On the right hand side of the above equation \eqref{eq:model_general},
the first term, $\Lambda=L-D$, represents linear dispersion and dissipation
effects, where $L^{*}=-L$ is an energy-conserving skew-symmetric
operator; and $D<0$ is a negative definite operator. Inhomogeneous
forcing effects are introduced in a deterministic component, $F$,
and a stochastic component represented by a Gaussian random process,
$\sigma\left(t\right)\dot{W}\left(t;\omega\right)$. Most importantly,
nonlinear coupling effect has a non-negligible contribution in the
dynamical system introduced via a quadratic form, $B\left(u,u\right)$,
which satisfies the energy conservation law by $u\cdot B\left(u,u\right)=0$.
The model structures in \eqref{eq:model_general} are representative
in a wide variety of multiscale systems in many fields \cite{lesieur1987turbulence,kalnay2003atmospheric,pedlosky2013geophysical}.
In computing key statistical predictions in the model state $u$,
the low-order moments become intricately connected to the high-order
statistical information due to the nonlinear coupling $B\left(u,u\right)$.
Conventional ensemble Kalman filters \cite{evensen1994sequential,doucet2001sequential,bocquet2017degenerate}
neglecting the higher-order moments information usually become insufficient.
In this case even with low dimensionality $d$, finite ensemble approximation
frequently suffers from collapse of particles, with the group of particles
concentrating in the center region of the PDF and failing to capture
the outliers charactering the key non-Gaussian statistics and extreme
events \cite{mohamad2018sequential,tong2021extreme,qi2023data}. Thus
effective algorithms require to capture the entire probability density
functions (PDFs) including high-order information using a moderate
ensemble size to maintain the affordable computational cost. 

In this paper, we introduce a practical modeling and computational
strategy designed to accurately capture the probability distributions
and key statistical characteristics of the solution to \eqref{eq:model_general}.
Based on the theoretical framework presented in \cite{qi2024coupled}
and the statistical filtering approach in \cite{bach2023filtering},
we develop a detailed data assimilation algorithm aimed at achieving
accurate statistical prediction in finite ensemble approximation of
potentially highly non-Gaussian probability distributions. To cope
with the computational limitation in practical applications, the ensemble
of simulated samples needs to be constrained in a small size. We propose
to correct the large fluctuating errors that commonly appear using
small ensemble size by exploiting partial observation data of the
low-order statistical moments. An effective data assimilation algorithm
is then formulated to improve the model accuracy by capturing the
higher-order moment information and reduce the high computational
cost at the same time. In particular, we conduct detailed numerical
study on the proposed ensemble data assimilation scheme based on the
systematic statistical modeling framework and applied on a representative
triad system \cite{majda2019linear} with multiple distinctive statistical
regimes. 

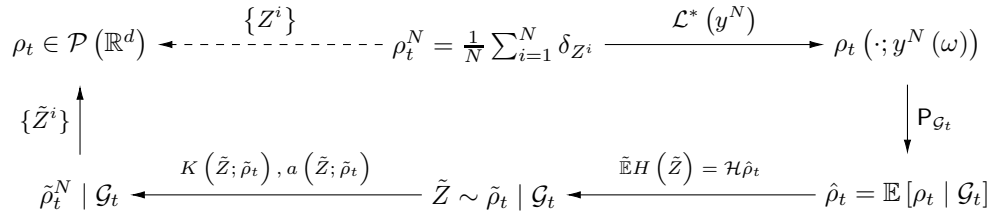
\begin{figure}[h]
\centering
\begin{tikzpicture}[node distance = 5cm, auto] 
		\node [pb] (target) {$\rho_{t}\in\mathcal{P}\left(\mathbb{R}^{d}\right)$};
		\node [pb,right of=target,node distance=5.5cm] (discrete) {$\rho_{t}^{N}=\frac{1}{N}\sum_{i=1}^{N}\delta_{Z^{i}}$};
		\node [pb, right of=discrete, node distance=5.5cm] (stoch) {$\rho_{t}\left(\cdot;y^{N}\left(\omega\right)\right)$};

		\node [pb,  below of=stoch, node distance=2.cm] (filter) {$\hat{\rho}_{t}=\mathbb{E}\left[\rho_{t}\mid\mathcal{G}_{t}\right]$};
		\node [pb,  below of=discrete, node distance=2.cm] (meanfield) {$\tilde{Z} \sim \tilde{\rho}_{t}\mid\mathcal{G}_{t}$};
		\node [pb,  below of=target, node distance=2.cm] (model) {$\tilde{\rho}_{t}^{N}\mid\mathcal{G}_{t}$};
 
		\path[draw, -{Latex[length=2mm,width=1mm]}, dashed] (discrete) -- (target) node[midway, above,font=\small] {$\left\{ Z^{i}\right\}$}; 
		\path[draw, -{Latex[length=2mm,width=1mm]}] (discrete) -- (stoch) node[midway, above,font=\small] {$\mathcal{L}^{*}\left(y^{N}\right)$};
		\path[draw, -{Latex[length=2mm,width=1mm]}] (stoch) -- (filter) node[midway, right,font=\small] {$\mathsf{P}_{\mathcal{G}_{t}}$};
		\path[draw, -{Latex[length=2mm,width=1mm]}] (filter) -- (meanfield) node[midway, above,font=\scriptsize] {$\tilde{\mathbb{E}}H\left(\tilde{Z}\right)=\mathcal{H}\hat{\rho}_{t}$};
		\path[draw, -{Latex[length=2mm,width=1mm]}] (meanfield) -- (model) node[midway, above,font=\scriptsize] {$K\left(\tilde{Z};\tilde{\rho}_{t}\right), a\left(\tilde{Z};\tilde{\rho}_{t}\right)$}; 
		\path[draw, -{Latex[length=2mm,width=1mm]}] (model) -- (target) node[midway, left,font=\small] {$\{ \tilde{Z}^{i}\}$}; 
		
\end{tikzpicture}

\caption{Flow chart illustrating ideas in constructing the data assimilation
model for statistical forecast.\protect\label{fig:Diagram}}
\end{figure}

The main ideas in constructing the data assimilation model using statistical
observation data is illustrated in the flow chart in Figure~\ref{fig:Diagram}.
We propose to compute the probability distribution $\rho$ of the
model state $u$ in \eqref{eq:model_general} using the more tractable
stochastic-statistical equations \eqref{eq:full_model}. The target
probability distribution will be approximated by an empirical probability
distribution $\rho^{N}$ through an interacting particle simulation
of the stochastic coefficients $Z$. However, in practice this ensemble-based
approach will often become insufficient to accurately capture the
essential PDF structures when only a small sample size $N$ is available.
To study the uncertainty from a finite sample size, we consider the
continuous distribution $\rho\left(\cdot;y^{N}\right)$ of each particle
$Z^{i}$ as a $\mathcal{P}\left(\mathbb{R}^{d}\right)$-valued random
field defined by \eqref{eq:filter_dyn}. The randomness is introduced
due to the finite sample estimation of the leading moments $y^{N}\left(\omega\right)$
in the statistical equation \eqref{eq:obs_approx}. This leads to
a natural filtering problem to find the optimal state estimation of
the random field $\rho$ based on the observation data $\mathcal{G}_{t}$
generated by the low-order statistical moments with respect to the
probability distribution. The optimal filter solution $\hat{\rho}$
then can be found through the projection on the space of $\mathcal{G}_{t}$-measurable
square-integrable random fields and is given by the Kalman-Bucy filter
\eqref{eq:KB-pdf} as an infinite-dimensional functional equations.
To propose an efficient computational strategy to solve $\hat{\rho}$,
a new stochastic process $\tilde{Z}\sim\tilde{\rho}$ is introduced
so that they can provide consistent high-order statistics $\ensuremath{\tilde{\mathbb{E}}H=\mathcal{H}\hat{\rho}}$
according to the observation function $H$. In particular, the governing
SDE for the new process $\tilde{Z}$ \eqref{eq:filtering_full} is
derived with explicit forms of the filtering coefficients $a,K$ in
\eqref{eq:filter_scheme}. Finally, this new probability distribution
$\tilde{\rho}$ can be computed efficiently by a finite sample approximation
$\tilde{\rho}^{N}$ which provides accurate high-order statistical
consistency and generates samples giving a better representation of
the target distribution of the model. 

In the structure of this paper, we first review the general multiscale
modeling framework and develop the data assimilation model based on
the coupling structure in the stochastic-statistical model in Section~\ref{sec:An-integrated-multiscale}.
Then, the detailed the ensemble data assimilation equations involving
the explicit filtering operators and practical computational algorithms
are constructed in Section~\ref{sec:Ensemble-data-assimilation}.
The performance of the new data assimilation model and its skill in
recovering both leading-order mean and covariance and the crucial
higher-order statistical feedbacks are extensively tested under a
representative prototype triad system demonstrating different statistical
regimes in Section~\ref{sec:Numerical-performance}. A summarizing
discussion and potential future research directions are given in Section~\ref{sec:Summary}.
Additional proofs of the results used in the main text and more details
on the intuition and applications of the triad system are listed in
Appendices~\ref{sec:Proofs-of-theorems} and~\ref{sec:Details-on-triad}.

\section{An integrated multiscale modeling framework with data assimilation\protect\label{sec:An-integrated-multiscale}}

To start with, we describe the multiscale modeling strategy for the
statistical solution to the general system \eqref{eq:model_general}.
In particular, the resulting coupled stochastic-statistical equations
and the associated ensemble approximation can be naturally combined
with data assimilation for improved sampling of the target probability
distributions.

\subsection{The coupled stochastic-statistical formulation for multiscale dynamics}

To characterize the uncertainty in the stochastic model state, the
solution $u$ is represented as a random field (denoted by $\omega$)
and decomposed into the multiscale composition of a statistical mean
state $\bar{u}=\mathbb{E}\left(u\right)$ and stochastic fluctuations
$u^{\prime}$ in a high-dimensional representation under a proper
orthonormal basis $\left\{ \hat{v}_{k}\right\} _{k=1}^{d}$, that
is,
\begin{equation}
u\left(t;\omega\right)=\bar{u}\left(t\right)+u^{\prime}\left(t;\omega\right)=\bar{u}\left(t\right)+\sum_{k=1}^{d}Z_{k}\left(t;\omega\right)\hat{v}_{k}.\label{eq:decomp}
\end{equation}
Above, $\bar{u}$ represents the statistical mean field of the dominant
largest scale structure; and $Z\left(t;\omega\right)=\left\{ Z_{k}\left(t;\omega\right)\right\} _{k=1}^{d}$
are the stochastic processes characterizing the uncertainty in the
fluctuation processes $u^{\prime}$ on each eigenmodes $\hat{v}_{k}$.
Such decomposition is commonly used, for example, in the describing
the zonal jets in geophysical turbulence and the coherent radial flow
in fusion plasmas \cite{salmon1998lectures,nicholson1983introduction}.

Under the decomposition \eqref{eq:decomp}, we can separate the full
statistics in the original stochastic state $u$ by the leading two
statistical moments $\bar{u},R$ and a mean-zero stochastic process
$Z$ governed by a coupled statistical and stochastic dynamical equations.
In particular, the \emph{statistical dynamical equations} describing
the evolution of the mean $\bar{u}\left(t\right)\in\mathbb{R}^{d}$
and covariance $R\left(t\right)\in\mathbb{R}^{d\times d}$ can be
derived as\addtocounter{equation}{0}\begin{subequations}\label{eq:full_model} 
\begin{equation}
\begin{aligned}\frac{\mathrm{d}\bar{u}}{\mathrm{d}t} & =\Lambda\bar{u}+B\left(\bar{u},\bar{u}\right)+F+\sum_{k,l=1}^{d}B\left(\hat{v}_{k},\hat{v}_{l}\right)\mathbb{E}\left(Z_{k}Z_{l}\right),\\
\frac{\mathrm{d}R}{\mathrm{d}t} & =L\left(\bar{u}\right)R+RL^{T}\left(\bar{u}\right)+Q_{\sigma}+Q_{F}\left(\mathbb{E}\left(Z\otimes Z\otimes Z\right)\right).
\end{aligned}
\label{eq:dyn_stat}
\end{equation}
Accordingly, the stochastic process $Z\left(t;\omega\right)\in\mathbb{R}^{d}$
satisfies the following \emph{stochastic differential equation }coupled
with the leading-order statistical moments $\left(\bar{u},R\right)$
solved from \eqref{eq:dyn_stat}
\begin{equation}
\mathrm{d}Z=L\left(\bar{u}\right)Z\mathrm{d}t+Q_{v}\left(Z\otimes Z-R\right)\mathrm{d}t+\sigma\mathrm{d}W.\label{eq:dyn_stoc}
\end{equation}
\end{subequations}In the statistical equations \eqref{eq:dyn_stat},
we define the stochastic forcing $Q_{\sigma}=\sigma\sigma^{T}$ from
the white noise process, the coupling coefficients $L\left(\bar{u}\right)$
due to interactions between the mean and covariance, and $Q_{F}$
due to the higher-order moments feedback from the triad modes $Z\otimes Z\otimes Z=\left\{ Z_{m}Z_{n}Z_{k}\right\} _{m,n,k}$.
In stochastic equation \eqref{eq:dyn_stoc}, $Q_{v}$ represents the
stochastic quadratic coupling between modes from $Z\otimes Z=ZZ^{T}$.
The explicit expressions for these coefficients can be found according
to the nonlinear coupling function $B\left(u,u\right)$ for all the
modes $k,l=1,\cdots,d$ as
\begin{equation}
\begin{aligned}L_{kl}\left(\bar{u}\right)= & \hat{v}_{k}\cdot\left[\Lambda\hat{v}_{l}+B\left(\bar{u},\hat{v}_{l}\right)+B\left(\hat{v}_{l},\bar{u}\right)\right],\\
Q_{F,kl}= & \sum_{m,n=1}^{d}\left[\gamma_{kmn}\mathbb{E}\left(Z_{m}Z_{n}Z_{l}\right)+\gamma_{lmn}\mathbb{E}\left(Z_{m}Z_{n}Z_{k}\right)\right],\\
Q_{v,k}= & \sum_{m,n=1}^{d}\gamma_{kmn}\left(Z_{m}Z_{n}-R_{mn}\right),
\end{aligned}
\label{eq:coupling_coeff}
\end{equation}
with the coupling coefficients $\gamma_{kmn}=\hat{v}_{k}\cdot B\left(\hat{v}_{m},\hat{v}_{n}\right)$. 

The above \emph{coupled stochastic-statistical equations} \eqref{eq:full_model}
provide a self-consistent closed system for recovering the statistical
solution of $u$. The leading-order moments $\bar{u}$ and $R$ are
solved by the statistical equations \eqref{eq:dyn_stat} involving
the higher-order moments of the stochastic coefficients $Z$, and
all high-order statistical information is recovered through the law
$\rho$ of the stochastic process $Z$ from \eqref{eq:dyn_stoc} dependent
on the solutions $\left\{ \bar{u},R\right\} $. More detailed discussions
on the derivation and advantages of this new formulation can be found
in \cite{qi2024coupled}. It demonstrates that this new set of equations
provides consistent statistical solutions with the original system
\eqref{eq:model_general}, while enjoys additional advantages of more
adaptive to various model reduction and data assimilation strategies
\cite{majda2018strategies,majda2019linear}.

\subsection{Predicting probability density functions using statistical observation
data}

A practical approach for numerically implementing the coupled stochastic-statistical
equations \eqref{eq:full_model} is to adopt a particle approximation
to the probability distribution of the stochastic variable $Z$. Thus
the expectations required in the statistical equations \eqref{eq:dyn_stat}
can be estimated through an empirical average of the samples $\mathbf{Z}=\left\{ Z^{i}\right\} _{i=1}^{N}$
\begin{equation}
\rho^{N}\left(z,t\right)=\frac{1}{N}\sum_{i=1}^{N}\delta\left(z-Z^{i}\left(t\right)\right),\quad\mathbb{E}^{N}\left(f\left(\mathbf{Z}\right)\right)=\frac{1}{N}\sum_{i=1}^{N}f\left(Z^{i}\right).\label{eq:emprical_pdf}
\end{equation}
Therefore, the statistical solution can be computed by solving the
following equations as an interacting particle system by evolving
the ensemble $\mathbf{Z}$
\begin{equation}
\begin{aligned}\frac{\mathrm{d}Z^{i}}{\mathrm{d}t}= & \:L\left(\bar{u}^{N}\right)Z^{i}+Q_{v}\left(Z^{i}\otimes Z^{i}-R^{N}\right)+\sigma\dot{W}^{i},\quad i=1,\cdots,N,\\
\frac{\mathrm{d}\bar{u}_{}^{N}}{\mathrm{d}t}= & \:\Lambda\bar{u}^{N}+B\left(\bar{u}^{N},\bar{u}^{N}\right)+\sum_{k,l}B\left(\hat{v}_{k},\hat{v}_{l}\right)\mathbb{E}^{N}\left(\mathbf{Z}_{k}\otimes\mathbf{Z}_{l}\right)+F,\\
\frac{\mathrm{d}R^{N}}{\mathrm{d}t}= & \:L\left(\bar{u}^{N}\right)R_{}^{N}+R_{}^{N}L^{T}\left(\bar{u}^{N}\right)+Q_{F}\left(\mathbb{E}^{N}\left(\mathbf{Z}\otimes\mathbf{Z}\otimes\mathbf{Z}\right)\right)+Q_{\sigma}\\
 & \qquad+\epsilon^{-1}\left(\mathbb{E}^{N}\left[\mathbf{Z}\otimes\mathbf{Z}\right]-R_{t}^{N}\right).
\end{aligned}
\label{eq:num_model}
\end{equation}
Several modifications are introduced in the numerical model \eqref{eq:num_model}
compared to the original equations \eqref{eq:full_model}. Instead
of computing the exact law $\rho\left(z,t\right)$ of the stochastic
process $Z$ by solving the PDE \eqref{eq:dyn_pdf}, a finite particle
approximation in the form of \eqref{eq:emprical_pdf} is used to estimate
the crucial higher moments feedback in the mean and covariance equations.
The samples are generated by a McKean-Vlasov SDE implicitly dependent
on all the sample trajectories through the statistical solutions $\bar{u}^{N},R^{N}$.
In addition, a relaxation term with an additional parameter $\epsilon>0$
is added to the covariance equation for $R^{N}$ to enforce consistency
in the finite particle approximation of the covariance. It is found
that this term is essential for maintaining stable numerical especially
with strong mean-fluctuation coupling from the term $L\left(\bar{u}_{t}\right)$
(see Figure~\ref{fig:relaxation} in Section~\ref{subsec:Typical-statistical-regimes}).

It can be shown \cite{graham1996asymptotic} that the empirical measure
$\rho^{N}$ converges to the law of each $Z^{i}$, $\rho^{N}\rightarrow\rho$,
as $N\rightarrow\infty$. The solution of the probability distribution
$\rho$ of $Z^{i}$ is given by the corresponding Fokker-Planck equation
\begin{equation}
\frac{\partial\rho_{}}{\partial t}=\mathcal{L}^{*}\left(\bar{u},R\right)\rho\coloneqq-\nabla\cdot\left[L\left(\bar{u}\right)z\rho+Q_{v}\left(zz^{T}-R\right)\rho\right]+\frac{1}{2}\nabla\cdot\left[\nabla\cdot\left(Q_{\sigma}\rho\right)\right],\label{eq:dyn_pdf}
\end{equation}
where $\mathcal{L}_{}^{*}$ is the adjoint of the generator $\mathcal{L}$
that is dependent on the mean $\bar{u}_{}$ and covariance $R$. Still,
a major difficulty remains if only a very small number of samples
$N$ is affordable to estimate the empirical distribution $\rho^{N}$.
Furthermore, the internal instability (that is, the positive eigenvalues
in $L\left(\bar{u}\right)$ due to the mean-fluctuation coupling $B\left(\bar{u},\hat{v}_{l}\right)\cdot\hat{v}_{k}$)
may lead to fast growth of the sample errors and quick divergence
of the numerical solutions (see \cite{majda2019linear,qi2024unambiguous}
and Figure~\ref{fig:Statistical-forecasts} in Section~\ref{subsec:Direct-prediction}).
This sets an inherent obstacle for efficient computation of the statistical
solutions in the multiscale coupling system. To address this inherent
difficulty, we assume that additional observation data $y=\left\{ \bar{u},R\right\} $
containing only the first two moments is available to improve the
prediction of the probability distribution $\rho^{N}$ from the stochastic
samples. Especially when the nonlinear coupling plays a dominant role
in the dynamics, many non-Gaussian features will emerge in the probability
distribution $\rho$. Our goal is then to filter the non-Gaussian
PDF $\rho$ of $Z$ containing crucial higher-order statistics by
taking into account observations from the leading two moments that
are easy to access.

In order to introduce the data assimilation strategies to improve
the prediction of the computational model, we rewrite the coupled
stochastic-statistical equations \eqref{eq:num_model} as a conditional
linear system about the probability distribution $\rho$ of the ensemble
member $Z^{i}$ and the empirical moments $\bar{u}^{N},R^{N}$ also
become stochastic processes
\begin{equation}
\begin{aligned}\partial_{t}\rho= & \mathcal{L}^{*}\left(\bar{u}^{N},R^{N}\right)\rho^ {},\\
\mathrm{d}\bar{u}^{N}= & \left[\mathbb{E}H^{m}\left(Z\right)+h_{m}\left(\bar{u}^{N}\right)\right]\mathrm{d}t+\Gamma_{m}^{N}\mathrm{d}B_{m},\\
\mathrm{d}R^{N}= & \left[\mathbb{E}H^{v}\left(Z\right)+h_{v}\left(\bar{u}^{N},R^{N}\right)\right]\mathrm{d}t+\Gamma_{v}^{N}\mathrm{d}B_{v}.
\end{aligned}
\label{eq:filter_dyn}
\end{equation}
Above, $\mathcal{L}\left(\bar{u},R\right)$ is given by the infinitesimal
generator in \eqref{eq:dyn_pdf} and $\rho$ is the continuous density
solution. In the equations for $\bar{u}^{N},R^{N}$, we summarize
all the deterministic terms (that is, all terms in \eqref{eq:num_model}
beside the ones with $\mathbb{E}^{N}$) in the functions $h_{m}$
and $h_{v}$ respectively. Higher-order moment feedbacks can be then
written as expectations with respect to the continuous probability
distribution $\rho$. The explicit expressions for the observation
functions $H^{m}\in\mathbb{R}^{d}$ and $H^{v}\in\mathbb{R}^{d\times d}$
can be found from the expressions in \eqref{eq:coupling_coeff} as
quadratic and cubic functions about $z$
\begin{equation}
H_{k}^{m}\left(z\right)=\sum_{p,q=1}^{d}\gamma_{kpq}z_{p}z_{q},\quad H_{kl}^{v}\left(z\right)=\sum_{p,q=1}^{d}\left(\gamma_{kpq}z_{p}z_{q}z_{l}+\gamma_{lpq}z_{p}z_{q}z_{k}\right).\label{eq:operators_obs}
\end{equation}
Importantly, notice that the additional noise terms with coefficients
$\left(\Gamma_{m}^{N},\Gamma_{v}^{N}\right)$ are introduced in the
observed states $\left(\bar{u}^{N},R^{N}\right)$ to calibrate the
errors from finite ensemble approximation. In fact, we can assume
that the empirical average in the mean and covariance equations can
be both decomposed into the expectation with $\rho$ and the additional
noise as a correction term to the finite sample estimation
\begin{equation}
\mathbb{E}^{N}H\left(Z\right)\mathrm{d}t=\mathbb{E}H\left(Z\right)\mathrm{d}t+\Gamma^{N}\mathrm{d}B.\label{eq:obs_approx}
\end{equation}
Thus, the additional noise term $\Gamma^{N}\mathrm{d}B$ is used to
represent the fluctuating error from the $N$ samples approximation
to the true expectation $\mathbb{E}H\left(Z\right)$. It is confirmed
from the numerical tests in Section~\ref{subsec:Calibration-of-observation}
that \eqref{eq:obs_approx} offers desirable characterization of the
approximation errors in practical applications. In this way, the coupled
system \eqref{eq:filter_dyn} sets up a standard linear filtering
problem given by the signal process $\rho$ as a $\mathcal{P}\left(\mathbb{R}^{d}\right)$-valued
stochastic process and the observation process $\mathcal{G}_{t}=\sigma\left\{ y\left(s\right),s\leq t\right\} $
with $y\left(t\right)=\left\{ \bar{u}^{N}\left(t\right),R^{N}\left(t\right)\right\} $
satisfying the linear equation with respect to the signal process
$\rho$
\begin{equation}
\mathrm{d}y_{}^ {}=\left[\mathcal{H}\rho+h\left(y^ {}\right)\right]\mathrm{d}t+\Gamma\mathrm{d}B,\quad\mathcal{H}\rho=\mathbb{E}H\left(Z\right)=\int H\left(z\right)\rho\left(z\right)\mathrm{d}z,\label{eq:obs_general}
\end{equation}
where $\mathcal{H}$ becomes a linear operator acting on the probability
density $\rho$ with $H=\left[H^{m},H^{v}\right]$, $h=\left[h_{m},h_{v}\right]$,
and $\Gamma\mathrm{d}B=\left[\Gamma_{m}\mathrm{d}B_{m},\Gamma_{v}\mathrm{d}B_{v}\right]$. 

Applying the Kalman-Bucy filter in the infinite-dimensional Hilbert
space \cite{liptser2013statistics,curtain1975infinite} for the stochastic
process $\rho$ in \eqref{eq:filter_dyn} conditional on the observation
processes in \eqref{eq:obs_general}, we find the \emph{optimal high-order
filter solution} $\hat{\rho}=\mathbb{E}\left[\rho_{}\mid\mathcal{G}_{t}\right]$
satisfying the following closed system of functional equations
\begin{equation}
\begin{aligned}\mathrm{d}\hat{\rho}_{}= & \:\mathcal{L}^{*}\left(y_{}\right)\hat{\rho}_{}\mathrm{d}t+\mathcal{\hat{C}}\mathcal{H}^{*}\Gamma^{-2}\left\{ \mathrm{d}y_{t}-\left[\mathcal{H}\hat{\rho}+h\left(y_{}\right)\right]\mathrm{d}t\right\} ,\\
\mathrm{d}\hat{\mathcal{C}}_{}= & \left[\mathcal{L}_{}^{*}\left(y_{}\right)\mathcal{\hat{C}}_{}+\mathcal{\hat{C}}\mathcal{L}\left(y\right)\right]\mathrm{d}t-\mathcal{\hat{C}}_{}\mathcal{H}^{*}\Gamma^{-2}\mathcal{H}\mathcal{\hat{C}}\mathrm{d}t,
\end{aligned}
\label{eq:KB-pdf}
\end{equation}
where $\hat{\mathcal{C}}\left(\omega\right):L^{2}\left(\mathbb{R}^{d}\right)\rightarrow L^{2}\left(\mathbb{R}^{d}\right)$
is the self-adjoint covariance operator with $\hat{\mathcal{C}}^{*}=\hat{\mathcal{C}}$.
The idea of filtering the probability distributions starts from the
Fokker-Planck filter in \cite{bach2023filtering}, and a systematic
filtering model is developed in \cite{qi2024coupled} based on the
specific nonlinear coupling structure in the stochastic-statistical
model \eqref{eq:full_model}. Still, it remains intractable to directly
solve the functional system \eqref{eq:KB-pdf}. The final step is
to construct effective ensemble solvers for the above optimal filter
solution $\hat{\rho}$.

\subsection{The approximate ensemble filter with consistent high-order statistics}

As a final step, we introduce a practical strategy to efficiently
compute the optimal filtering solution $\hat{\rho}$ based on the
observed statistics. The idea is to construct a surrogate process
$\tilde{Z}$ so that the corresponding probability distribution of
$\tilde{Z}\sim\tilde{\rho}$ can serve as an effective representation
of the optimal filter solution $\hat{\rho}$. Then, efficient particle
approaches can be adopted to capture the probability distribution
of $\tilde{Z}$ instead of solving the infinite-dimensional equations
\eqref{eq:KB-pdf}.

Associated with the forecast equation \eqref{eq:num_model} for the
stochastic process $Z$, given $N$ particles for the sampling solution
of the stochastic state, $\tilde{\mathbf{Z}}=\left\{ \tilde{Z}^{i}\right\} _{i=1}^{N}$,
we can construct the following filtering updating equation with an
additional update according to the observation in the second line
\begin{equation}
\begin{aligned}\mathrm{d}\tilde{Z}^{i}= & L\left(\bar{u}_{}^{N}\right)\tilde{Z}_{}^{i}\mathrm{d}t+Q_{v}\left(\tilde{Z}_{}^{i}\otimes\tilde{Z}^{i}-R^{N}\right)\mathrm{d}t+\sigma\mathrm{d}\tilde{W}_{}^{i}\\
+ & a^{m}\left(\tilde{Z}^{i};\tilde{\rho}\right)\mathrm{d}t+K^{m}\left(\tilde{Z}^{i};\tilde{\rho}\right)\mathrm{d}I^{m}+a^{v}\left(\tilde{Z}^{i};\tilde{\rho}\right)\mathrm{d}t+K^{v}\left(\tilde{Z}^{i};\tilde{\rho}\right)\mathrm{d}I^{v},
\end{aligned}
\label{eq:filtering_full}
\end{equation}
where the innovations $I^{m},I^{v}$ for the statistical observations
are defined based on the observation data $\mathrm{d}y=\left(\mathrm{d}\bar{u},\mathrm{d}R\right)$
as
\begin{equation}
\begin{aligned}\mathrm{d}I^{m}\left(t\right)= & \:\mathrm{d}\bar{u}-\left[H^{m}\left(\tilde{Z}^{i}\right)+h_{m}\left(\bar{u}_{}^{N}\right)\right]\mathrm{d}t,\\
\mathrm{d}I^{v}\left(t\right)= & \:\mathrm{d}R-\left[H^{v}\left(\tilde{Z}^{i}\right)+h_{v}\left(\bar{u}_{}^{N},R_{}^{N}\right)\right]\mathrm{d}t.
\end{aligned}
\label{eq:model_innov}
\end{equation}
In the new filter equation \eqref{eq:filtering_full} for the stochastic
process $\tilde{Z}$, the first line follows the same dynamical equation
\eqref{eq:dyn_stoc} as the forecast step, while the second line introduces
additional control correction based on the observation data. The new
functionals known as the Kalman gain $K$ and the drift $a$ are defined
based on the probability distribution of $\tilde{Z}^{i}\sim\tilde{\rho}$,
and needs to be solved by the equations \eqref{eq:filter_ops} shown
below. Notice that the filter equation \eqref{eq:filtering_full}
is also dependent on the first two moments $\left(\bar{u}^{N},R^{N}\right)$
which can be solved by integrating the two statistical equations in
\eqref{eq:num_model}. In addition, we need to introduce continuous
observations to estimate $\mathrm{d}y_{t}$ in the filtering scheme.
Assume that the observation data $y_{n}=y\left(t_{n}\right)$ comes
at times $t_{n}=n\Delta t$ with a short observation interval $\Delta t$.
We can approximate the increment at $t_{n}$ in the observation data
from the linear interpolation for $t\in\left(t_{n},t_{n+1}\right]$
\begin{equation}
\mathrm{d}y\left(t\right)\simeq\Delta y_{n}=y_{n+1}-y_{n}.\label{eq:obs_interp}
\end{equation}
Furthermore, even though consecutive observations are required in
a short interval for the estimate of $\Delta y_{n}$, we may not require
to have continuous observation data at each time updating step. The
updating step of data assimilation in the second line of the equation
\eqref{eq:filtering_full} can be applied only at the steps when observation
data $\Delta y_{n}$ is received.

Finally, the development of effective filtering scheme then relies
on the construction of the Kalman gain operators $\left(K^{m},K^{v}\right)$
and the drift terms $\left(a^{m},a^{v}\right)$ in the second line
of the filter update in \eqref{eq:filtering_full} based on the observation
functions $H^{m},H^{v}$. In general, these terms should be chosen
so that the probability distribution $\tilde{\rho}$ of the constructed
stochastic $\tilde{Z}^{i}$ can correctly reflect the optimal filter
solution $\hat{\rho}$ in \eqref{eq:KB-pdf}. From the standard procedure
of developing the mean field evolution equations \cite{calvello2022ensemble,pathiraja2021mckean},
they can be solved according to the probability distribution $\tilde{\rho}$
of $\tilde{Z}$, such that
\begin{equation}
-\nabla\cdot\left(K^{T}\tilde{\rho}\right)=\tilde{\rho}\Gamma^{-2}\left[H\left(\tilde{Z}\right)-\tilde{\mathbb{E}}\left(H\right)\right],\quad a=\nabla\cdot\left(K\Gamma^{2}K^{T}\right)-K\Gamma^{2}\nabla\cdot K^{T}.\label{eq:filter_ops}
\end{equation}
Theoretical analysis \cite{qi2024coupled} shows that the high-order
filtering equation \eqref{eq:filtering_full} generates consistent
statistics with the optimal filter solution \eqref{eq:KB-pdf} in
the analysis update if the above conditions \eqref{eq:filter_ops}
are satisfied. Importantly, high-order statistics according to the
nonlinear observation operators $H$ in \eqref{eq:operators_obs},
$\tilde{\mathbb{E}}\left[H\left(\tilde{Z}\right)\right]=\mathcal{H}\hat{\rho}$,
are preserved in the filter updating equation. This indicates that
the new filtering model for $\tilde{Z}$ is able to capture the crucial
high-order statistics in the optimal filter solution $\hat{\rho}$
rather than only the first two moments in the linear ensemble Kalman
filters. However, efficient ways to compute the Kalman gain and drift
operators through \eqref{eq:filter_ops} are still needed without
losing the essential high-order moments information. In the next section,
we will propose easy-to-implement schemes to compute these key filter
operators $K$ and $a$ without the need to solve the distribution
function $\tilde{\rho}$. 

\section{Ensemble data assimilation schemes maintaining high-order statistics\protect\label{sec:Ensemble-data-assimilation}}

In this section, we construct a practical numerical scheme for the
ensemble filtering equation \eqref{eq:filtering_full}. The goal is
to generate a better empirical representation \eqref{eq:emprical_pdf}
of the probability distribution using only a small ensemble size.
The accurate computation of the model statistics requires that the
non-Gaussian statistics involved in the observation function $H^{m}$
and $H^{v}$ are properly represented in the filtering update. This
leads to several key treatments in the construction of the filtering
operators. In particular, we exploit the detailed structures of the
observation functions to derive explicit expressions for the functions
$a_{}^{m},K_{}^{m}$ and $a_{}^{v},K^{v}$ in the filter equation. 

\subsection{Construction of detailed filter operators with nonlinear observation
functions}

Assuming that the first $s$ components of the mean and covariance
are observed, we can derive the filter operators $\left(a,K\right)$
by exploiting the specific quadratic and cubic structures of the observation
functions $H^{m}\in\mathbb{R}^{s}$ and $H^{v}\in\mathbb{R}^{s\times s}$
for all the observed modes $1\leq k,l\leq s$
\begin{equation}
\begin{aligned}H_{k}^{m}\left(z\right)= & \sum_{p,q}\gamma_{kpq}z_{p}z_{q},\\
H_{kl}^{v}\left(z\right)= & \sum_{p,q}\gamma_{kpq}z_{p}z_{q}z_{l}+\gamma_{lpq}z_{p}z_{q}z_{k},
\end{aligned}
\label{eq:obs_func}
\end{equation}
where we define the coefficient $\gamma_{kpq}=\hat{v}_{k}\cdot B\left(\hat{v}_{p},\hat{v}_{q}\right)$
according to \eqref{eq:operators_obs}. The following property of
the observation functions $H^{m},H^{v}$ can be found from direct
computation using the above quadratic and cubic structures.
\begin{lem}
\label{lem:observation-functions}The observation functions $H^{m}$
and $H^{v}$ defined in \eqref{eq:obs_func} satisfy
\begin{equation}
\begin{aligned}z\cdot\nabla H^{m}\left(z\right)=\sum_{j}z_{j}\partial_{j}H^{m}= & \:2H^{m}\left(z\right),\\
z\cdot\nabla H^{v}\left(z\right)=\sum_{j}z_{j}\partial_{j}H^{v}= & \:3H^{v}\left(z\right).
\end{aligned}
\label{eq:obs_symm}
\end{equation}
\end{lem}

With the above symmetry in the observation functions \eqref{eq:obs_symm},
we are able to find explicit expressions for the Kalman gain and drift
operators that enable efficient computation of these terms. In the
following, we summarize the useful results and formulas and put detailed
derivations in Appendix \ref{sec:Proofs-of-theorems}.

\subsubsection{The explicit forms of the Kalman gain and drift operators}

With the above explicit relations for the observation functions, we
can first find special solutions to the equation \eqref{eq:filter_ops}
to recover the Kalman gain operator $K\left(z;\tilde{\rho}\right)$,
that is,
\[
-\nabla\cdot\left(K^{T}\tilde{\rho}\right)=\tilde{\rho}\Gamma^{-2}\left[H\left(\tilde{Z}\right)-\bar{H}\right],
\]
where $H=H^{m}$ or $H^{v}$ is the stretched vector and $\bar{H}=\tilde{\mathbb{E}}\left[H\left(\tilde{Z}\right)\right]$.
Still, we would like to avoid directly solving the above equation
since the probability distribution $\tilde{\rho}$ is usually intractable
and can be only estimated from an ensemble approach. By multiplying
$H$ on both sides, the identity for $K$ implies a necessary condition
\begin{equation}
\tilde{\mathbb{E}}\left[K^{T}\nabla H\right]=\Gamma^{-2}C_{}^{H},\label{eq:kalman_gain}
\end{equation}
where $C^{H}=\tilde{\mathbb{E}}\left[\left(H-\bar{H}_{}\right)\left(H-\bar{H}\right)^{T}\right]$
is the second-order moment of $H$ with respect to $\tilde{\rho}$.
Then we may solve instead the above equation to find proper candidate
for the gain function $K_{t}$. In this way, we can compute the detailed
expressions for the Kalman gain operators based on the explicit forms
and their particular symmetry in the observation functions \eqref{eq:obs_symm}.
\begin{prop}
\label{prop:Kalman_gain}The following Kalman gain operators $K^{m}\in\mathbb{R}^{d\times s}$
and $K^{v}\in\mathbb{R}^{d\times s^{2}}$ satisfy the equation \eqref{eq:kalman_gain}
\begin{equation}
\begin{aligned}K^{m}\left(Z\right)= & \frac{1}{2}Z\left[H^{m}\left(Z\right)-\bar{H}^{m}\right]^{T}\Gamma_{m}^{-2},\\
K_{}^{v}\left(Z\right)= & \frac{1}{3}Z\left[H^{v}\left(Z\right)-\bar{H}^{v}\right]^{T}\Gamma_{v}^{-2}.
\end{aligned}
\label{eq:kalman_algm}
\end{equation}
with the state vector $Z\in\mathbb{R}^{d}$, observation function
$H^{m}\in\mathbb{R}^{s}$, $H^{v}\in\mathbb{R}^{s^{2}}$, observation
noises $\Gamma_{m}\in\mathbb{R}_{\mathrm{sym}}^{s\times s}$, $\Gamma_{v}\in\mathbb{R}_{\mathrm{sym}}^{s^{2}\times s^{2}}$,
and $\bar{H}=\tilde{\mathbb{E}}\left[H\left(Z\right)\right]$.
\end{prop}

Next, the function $a\left(z;\tilde{\rho}\right)$ can be directly
solved from the explicit solutions of $K_{}$ in \eqref{eq:kalman_algm}
according to 
\begin{equation}
a=\nabla\cdot\left(K\Gamma^{2}K^{T}\right)-K\Gamma^{2}\nabla\cdot\left(K^{T}\right).\label{eq:drift}
\end{equation}
We can also find the explicit forms of the drift terms through direct
computation using \eqref{eq:drift} and again the specific structures
in observation functions \eqref{eq:obs_symm}.
\begin{prop}
\label{prop:drift}With the forms of $K^{m}$and $K^{v}$ in \eqref{eq:kalman_algm},
the corresponding drift terms satisfying \eqref{eq:drift} can be
found as
\begin{equation}
\begin{aligned}a^{m}\left(Z\right)= & \frac{1}{4}Z\left[H^{m}\left(Z\right)-\bar{H}^{m}\right]^{T}\Gamma_{m}^{-2}\left[3H^{m}\left(Z\right)-\bar{H}^{m}\right],\\
a^{v}\left(Z\right)= & \frac{1}{9}Z\left[H^{v}\left(Z\right)-\bar{H}^{v}\right]^{T}\Gamma_{v}^{-2}\left[4H^{v}\left(Z\right)-\bar{H}^{v}\right].
\end{aligned}
\label{eq:drift_algm}
\end{equation}
\end{prop}

Notice that the solution to \eqref{eq:kalman_gain} only satisfies
a necessary condition for the original equation for the Kalman gain.
Still, it already accounts for the crucial high-order statistics with
respect to $\tilde{\rho}$ involving in the nonlinear observation
functions \eqref{eq:obs_func}. Therefore, the achieved explicit forms
of Kalman gain and drift terms \eqref{eq:kalman_algm} and \eqref{eq:drift_algm}
can serve as suitable candidate for the construction of high-order
filtering schemes. It shows to be a better choice than that in the
standard EnKF scheme (shown next in Section~\ref{subsec:Comparison-EnKF}
and the numerical comparisons in Section~\ref{sec:Numerical-performance})
which only considers Gaussian projection of $\tilde{\rho}$ thus neglects
the crucial high-order statistics information in $H$. 

\subsubsection{Numerical implementation of the filter operators}

Based on the explicit expressions of the filtering operators in \eqref{eq:kalman_algm}
and \eqref{eq:drift_algm}, we are able to construct direct algorithms
for effective implementation of the filtering scheme. At each time
updating step $t_{n}$, the mean and covariance can be computed by
integrating the statistical equations \eqref{eq:num_model} 
\begin{equation}
\begin{aligned}\Delta\bar{u}_{n}^{N}=\bar{u}_{n+1}^{N}-\bar{u}_{n}^{N}= & \int_{t_{n}}^{t_{n+1}}\left[\mathbb{E}^{N}H_{m}\left(\tilde{\mathbf{Z}}\left(s\right)\right)+h_{m}\left(\bar{u}^{N}\left(s\right)\right)\right]\mathrm{d}s,\\
\Delta R_{n}^{N}=R_{n+1}^{N}-R_{n}^{N}= & \int_{t_{n}}^{t_{n+1}}\left[\mathbb{E}^{N}H_{v}\left(\tilde{\mathbf{Z}}\left(s\right)\right)+h_{v}\left(\bar{u}^{N}\left(s\right),R^{N}\left(s\right)\right)\right]\mathrm{d}s.
\end{aligned}
\label{eq:update_stat}
\end{equation}
Above, the empirical expectation $\mathbb{E}^{N}\left(\cdot\right)$
is computed as in \eqref{eq:emprical_pdf} using the ensemble average
of all the simulated samples from the filter equation \eqref{eq:filtering_full}.
Then the filtering updating step combines the observation data $\left(\Delta\bar{u},\Delta R\right)$
in \eqref{eq:obs_interp} and the model forecast \eqref{eq:update_stat}
to get an optimal estimate for the ensemble distribution for $\tilde{Z}^{i}$
as 
\begin{equation}
\begin{aligned}\tilde{Z}_{n+1}^{i}= & \hat{Z}_{n+1}^{i}+\left(a^{m}\mathrm{\Delta}t+K^{m}\Delta I^{m}\right)+\left(a^{v}\mathrm{\Delta}t+K^{v}\Delta I^{v}\right),\\
\hat{Z}_{n+1}^{i}= & \tilde{Z}_{n}^{i}+L\left(\bar{u}_{n}^{N}\right)\tilde{Z}_{n}^{i}\mathrm{\Delta}t+Q_{v}\left(\tilde{Z}_{n}^{i}\otimes\tilde{Z}_{n}^{i}-R^{N}\right)\mathrm{\Delta}t+\sigma\Delta\tilde{W}_{n}^{i}.
\end{aligned}
\label{eq:filtering_update}
\end{equation}
Above, we split the filtering procedure in the standard two-step process,
where $\hat{Z}^{i}$ gets the forecast step update for the stochastic
state then the prior forecast is corrected through the filtering operators
when the observation data is available. The following proposition
provides the explicit expressions for directly computing the filter
updates using the samples and observation data.
\begin{prop}
\label{prop:filtering_update}Given the observations $\left(\Delta\bar{u},\Delta R\right)$
and the model predicted increments $\left(\Delta\bar{u}^{N},R^{N}\right)$,
the filtering updating step in the filtering equation \eqref{eq:filtering_update}
can be computed directly based on the ensemble member $\tilde{Z}^{i}$
as
\begin{equation}
\begin{aligned}a^{m}\mathrm{\Delta}t+K^{m}\Delta I^{m}= & \frac{1}{2}\left[\tilde{Z}_{}^{i}H_{m}^{\prime T}\left(\tilde{Z}_{}^{i}\right)\Gamma_{m}^{-2}\right]\left(\Delta\bar{u}_{}-\Delta\bar{u}_{}^{N}\right)\\
 & +\frac{\Delta t}{2}\left[\tilde{Z}^{i}H_{m}^{\prime T}\left(\tilde{Z}^{i}\right)\Gamma_{m}^{-2}\right]\bar{H}_{m}^ {}+\frac{\Delta t}{4}\left[\tilde{Z}^{i}H_{m}^{\prime T}\left(\tilde{Z}^{i}\right)\Gamma_{m}^{-2}H_{m}^{\prime}\left(\tilde{Z}^{i}\right)\right],\\
a^{v}\mathrm{\Delta}t+K^{v}\Delta I^{v}= & \frac{1}{3}\left[\tilde{Z}^{i}H_{v}^{\prime T}\left(\tilde{Z}^{i}\right)\Gamma_{v}^{-2}\right]\left(\Delta R-\Delta R^{N}\right)\\
 & +\frac{\Delta t}{3}\left[\tilde{Z}^{i}H_{v}^{\prime T}\left(\tilde{Z}^{i}\right)\Gamma_{v}^{-2}\right]\bar{H}_{v}+\frac{\Delta t}{9}\left[\tilde{Z}^{i}H_{v}^{\prime T}\left(\tilde{Z}^{i}\right)\Gamma_{v}^{-2}H_{v}^{\prime}\left(\tilde{Z}^{i}\right)\right].
\end{aligned}
\label{eq:filter_scheme}
\end{equation}
with $\bar{H}_{}^ {}=\mathbb{E}^{N}\left[H\left(\tilde{\mathbf{Z}}\right)\right]$
and $H^{\prime}=H\left(\tilde{Z}\right)-\bar{H}$.
\end{prop}

Using the explicit formulas in \eqref{eq:filtering_update}, we can
directly update each filtering sample $\tilde{Z}^{i}$ during the
time updating interval containing the higher-order moments information
in the observation functions $H_{m}$ and $H_{v}$. To summarize,
the dynamical equations for the particles are coupled through the
empirical average among all the samples according to Algorithm~\ref{alg:filter_full}.

\begin{algorithm}
\begin{algorithmic}[1] 
\Ensure{Introduce discrete time step $\Delta t$  with $M\Delta t=T$. The sequence of statistical observations are given by the increments of the mean and covariance $\Delta y_{n}=\left\{ \Delta\bar{u}_{n},\Delta R_{n}\right\}$ measured at the time instants $t_{n}=n\Delta t$.}
\Require{At initial time $t=0$, draw an ensemble of samples $\left\{\tilde{Z}_{0}^{i}\right\}_{i=1}^{N}$ from the initial distribution $\tilde{\rho}_0$.}

\For{$n = 0$ while $n < M$, during the time updating interval $t\in\left[t_{n},t_{n+1}\right]$.}
	\State{Integrate the samples to the next time step $\left\{\hat{Z}_{n+1}^{i}\right\}$ using forecast model in the second equation of \eqref{eq:filtering_update}.}
	\State{Integrate statistical mean and covariance to $\left\{\bar{u}^{N}_{n+1}, R^{N}_{n+1}\right\}$ by  \eqref{eq:update_stat} using the average of all samples.}
	\State{Compute the filtering update terms with the explicit formulas in \eqref{eq:filter_scheme} with the observation data $\Delta y_{n}$.}
    \State{Update the samples $\left\{\tilde{Z}_{n+1}^{i}\right\}$ from the prior states $\left\{\hat{Z}_{n+1}^{i}\right\}$  using the first equation of \eqref{eq:filtering_update}.}
    
\EndFor

\end{algorithmic}

\caption{Ensemble probability filter with statistical observations\protect\label{alg:filter_full}}
\end{algorithm}

\begin{rem*}
1. In practical implementations, it is observed that some outliers
of the samples $\tilde{Z}^{i}$ may introduce occasional instability
by creating some extremely large values in the high-order terms in
\eqref{eq:filter_scheme}. To improve stability in the highly unstable
regime, it is found useful to use the expectation values $\mathbb{E}^{N}\left[\tilde{Z}H_{}^{\prime T}\left(\tilde{Z}\right)\Gamma^{-2}\right]$
and $\mathbb{E}^{N}\left[\tilde{Z}H_{}^{\prime T}\left(\tilde{Z}^ {}\right)\Gamma^{-2}H_{}^{\prime}\left(\tilde{Z}\right)\right]$
instead of each ensemble evaluation to improve filter stability without
sacrificing too much accuracy.

2. High computational cost may still exist for solving a high-dimensional
SDE \eqref{eq:filtering_full}. This difficulty could be addressed
by the reduced-order algorithms such as the random batch methods \cite{qi2023high,qi2023random}.
We aim to combine the efficient forecast models with the data assimilation
strategy in high-dimensional problems in the following-up research.
\end{rem*}

\subsection{Comparison with ensemble Kalman filter schemes\protect\label{subsec:Comparison-EnKF}}

In comparison, we also illustrate the common strategy used in ensemble
Kalman filters \cite{calvello2022ensemble}. Assuming that the Kalman
gain $K$ is a deterministic matrix with no randomness, this leads
to the following choice of the deterministic Kalman gain matrix independent
of $\tilde{Z}$, and a zero drift term due to the constant Kalman
gain according to the equation \eqref{eq:filter_ops}
\begin{equation}
K=\tilde{\mathbb{E}}\left[\tilde{Z}\left(H\left(\tilde{Z}\right)-\bar{H}\right)^{T}\right]\Gamma^{-2}=\tilde{C}^{ZH}\Gamma^{-2},\quad\mathrm{and}\quad a=0.\label{eq:kalman_determ}
\end{equation}
where the covariance matrix $\tilde{C}^{ZH}$ is given by the cross-covariance
between the process $\tilde{Z}$ and the observation function $H\left(\tilde{Z}\right)$.
The ensemble Kalman filter scheme then follows with the following
filter equations 
\begin{equation}
\begin{aligned}\tilde{Z}_{n+1}^{i}= & \tilde{Z}_{n}^{i}+L\left(\bar{u}_{n}^{N}\right)\tilde{Z}_{n}^{i}\mathrm{\Delta}t+Q_{v}\left(\tilde{Z}_{n}^{i}\otimes\tilde{Z}_{n}^{i}-R^{N}\right)\mathrm{\Delta}t+\sigma\Delta\tilde{W}_{n}^{i}\\
+ & \tilde{C}^{ZH_{m}}\Gamma_{m}^{-2}\left\{ \Delta\bar{u}-\left[H_{m}\left(\tilde{Z}_{n}^{i}\right)+h_{m}\right]\Delta t\right\} +\tilde{C}^{ZH_{v}}\Gamma_{v}^{-2}\left\{ \Delta R-\left[H_{v}\left(\tilde{Z}_{n}^{i}\right)+h_{v}\right]\Delta t\right\} ,
\end{aligned}
\label{eq:EnKF}
\end{equation}
where the covariance $\hat{C}^{ZH}=\tilde{\mathbb{E}}^{N}\left[\hat{\mathbf{Z}}H^{\prime}\left(\hat{\mathbf{Z}}\right)^{T}\right]$
is computed according to the empirical average among all the sample
forecast $\tilde{\mathbf{Z}}_{n+1}$ from the second equation of \eqref{eq:filtering_update}.
The above updating scheme \eqref{eq:EnKF} is usually referred to
as the the ensemble Kalman filter (EnKF). It has been shown that the
EnKF approach can effectively drive the probability density functions
to the equilibrium such as using the ensemble Fokker-Planck filter
\cite{bach2023filtering}.

However, in our modeling framework involving the coupled stochastic-statistical
system, the high-order moments are playing a central role as the high-order
feedbacks in the statistical equations \eqref{eq:dyn_stat} for accurate
statistical prediction. Notice that in the EnKF approach, constant
Kalman gains $K_{m}$ and $K_{v}$ matrices are used independent of
each sample realization. It adopts the Gaussian projection on the
stochastic process $Z$ thus only statistics up to the second-order
moments are considered in the filter update. As a result, this approximation
deliberately neglected the crucial high-order statistics contained
in the observation functions $H_{m}$ and $H_{v}$. Compared with
the more precisely calibrated filter operators \eqref{eq:filter_scheme},
directly applying the EnKF in the coupled stochastic-statistical model
\eqref{eq:num_model} may miss the crucial high-order statistical
information in the sampled observation $H\left(\tilde{Z}^{i}\right)$
thus lead to larger errors and instability in the filter updates.
The degeneracy of particles failing to capture the essential high-order
information in the EnKF prediction is demonstrated in Figure~\ref{fig:Sampling}
and \ref{fig:Prediction-3mom} from direct numerical tests.

\subsection{Convergence of the statistical ensemble filter approximation}

In this final part, we discuss the convergence of the discrete numerical
scheme for the statistical estimates using finite ensemble approximation.
Let $\tilde{\rho}^{N}\left(z,t\right)=\frac{1}{N}\sum_{i=1}^{N}\delta\left(z-\tilde{Z}^{i}\left(t\right)\right)$
be the random field from the finite ensemble approximation of the
$N$ stochastic samples from \eqref{eq:filtering_update}, and $\tilde{\rho}\left(z,t\right)$
is the corresponding continuous distribution of the $\mathcal{P}\left(\mathbb{R}^{d}\right)$-valued
random field from the law of each stochastic process $\tilde{Z}_{t}$
in \eqref{eq:filtering_full} conditional on the observations $\mathcal{G}_{t}$. 

First, assume that the initial samples $\left\{ \tilde{Z}_{0}^{i}\right\} _{i=1}^{N}$
are drawn i.i.d. from the initial distribution $\tilde{\rho}_{0}$,
and the a unique solution exists for the McKean-Vlasov SDE \eqref{eq:filtering_full}
for each sample $\tilde{Z}^{i}\sim\tilde{\rho}$. From the established
conclusions from the limit of the $N$ interacting particle system
\cite{oelschlager1984martingale,graham1996asymptotic}, the empirical
probability distribution $\tilde{\rho}^{N}$ estimated with finite
samples will converge to the continuous measure $\tilde{\rho}$ as
the ensemble size $N\rightarrow\infty$. In addition, we will need
the following assumptions on the structures of the coupled stochastic-statistical
equations \eqref{eq:full_model}
\begin{assumption}
\label{assu:model_disc}Assume that the model dynamics functions $B:\mathbb{R}^{d}\times\mathbb{R}^{d}\rightarrow\mathbb{R}^{d}$
and $L:\mathbb{R}^{d}\rightarrow\mathbb{R}^{d\times d}$ in the mean
and covariance equations \eqref{eq:dyn_stat} are Lipschitz continuous,
that is, there is a constant $\beta>0$ so that
\[
\left|B\left(u,u\right)-B\left(v,v\right)\right|\leq\beta\left|u-v\right|,\quad\left\Vert L\left(u\right)-L\left(v\right)\right\Vert \leq\beta\left|u-v\right|.
\]
In addition, the nonlinear coefficients $\gamma_{kmn}=\hat{v}_{k}\cdot B\left(\hat{v}_{m},\hat{v}_{n}\right)$
in \eqref{eq:dyn_stoc} are uniformly bounded, that is, there exists
a constant $C>0$, so that for all $k,m,n$
\[
\left|\gamma_{kmn}\right|\leq C.
\]
\end{assumption}

Given the observations $\mathcal{G}_{t}$, we have for any test function
$\varphi\in C_{b}^{2}\left(\mathbb{R}^{d}\right)$ the empirical measure
$\tilde{\rho}^{N}$ converges to the continuous distribution for each
sample $\tilde{\rho}$
\begin{equation}
\left\langle \tilde{\rho}_{}^{N},\varphi\right\rangle =\frac{1}{N}\sum_{i=1}^{N}\varphi\left(\tilde{Z}^{i}\right)\rightarrow\mathbb{E}\left[\varphi\left(\tilde{Z}\right)\mid\mathcal{G}_{t}\right]=\left\langle \tilde{\rho},\varphi\right\rangle ,\label{eq:lln}
\end{equation}
a.s. as $N\rightarrow\infty$. Furthermore, there is the error estimate
for the empirical estimate $\left\langle \tilde{\rho}_{t}^{N},\varphi\right\rangle \coloneqq\mathbb{E}^{N}\varphi=\frac{1}{N}\sum\varphi\left(\tilde{Z}^{i}\left(t\right)\right)$
and $\left\langle \tilde{\rho}_{t},\varphi\right\rangle \coloneqq\tilde{\mathbb{E}}\left[\varphi\left(\tilde{Z}\left(t\right)\right)\right]$
for $T>0$
\begin{equation}
\mathbb{E}\left[\sup_{0\leq t\leq T}\left|\left\langle \tilde{\rho}_{t}^{N},\varphi\right\rangle -\left\langle \tilde{\rho}_{t},\varphi\right\rangle \right|^{2}\right]\leq\frac{C_{T}}{N}\left\Vert \varphi\right\Vert _{\infty}^{2}.\label{eq:conv_pdf}
\end{equation}
Proofs on \eqref{eq:lln} and \eqref{eq:conv_pdf} follow from the
law of large numbers and can be found in such as Theorem 9.18 of \cite{bain2009fundamentals}.

Next, we consider the finite ensemble and discrete time estimation
of the statistical mean and covariance states in the filtering solution.
The statistical equations for the continuous solutions $\left(\bar{u},R\right)$
can be written based on their explicit coupling dynamics \eqref{eq:dyn_stat}
and the high-order terms according to the observation functions $H^{m},H^{v}$
in \eqref{eq:operators_obs} with respect to the continuous probability
$\tilde{Z}\sim\tilde{\rho}$
\begin{equation}
\begin{aligned}\frac{\mathrm{d}\bar{u}}{\mathrm{d}t}= & \:\Lambda\bar{u}\left(t\right)+B\left(\bar{u}\left(t\right),\bar{u}\left(t\right)\right)+F+\tilde{\mathbb{E}}H^{m}\left(\tilde{Z}\left(t\right)\right),\\
\frac{\mathrm{d}R_{}}{\mathrm{d}t}= & \:L\left(\bar{u}\left(t\right)\right)R\left(t\right)+R\left(t\right)L\left(\bar{u}\left(t\right)\right)^{T}+Q_{\sigma}+\tilde{\mathbb{E}}H^{v}\left(\tilde{Z}\left(t\right)\right).
\end{aligned}
\label{eq:model_conti}
\end{equation}
On the other hand, the numerical mean and covariance estimates $\left(\bar{u}^{N,\delta},R^{N,\delta}\right)$
from the discrete time numerical updates and with the ensemble approximation
$\left\{ \tilde{Z}^{i}\right\} $ are computed from the equations
\eqref{eq:update_stat} with respect to the discrete empirical distribution
$\tilde{\rho}^{N}$
\begin{equation}
\begin{aligned}\frac{\mathrm{d}\bar{u}^{N,\delta}}{\mathrm{d}t}= & \:\Lambda\bar{u}^{N,\delta}\left(\tau\left(t\right)\right)+B\left(\bar{u}^{N,\delta}\left(\tau\left(t\right)\right),\bar{u}^{N,\delta}\left(\tau\left(t\right)\right)\right)+F+\frac{1}{N}\sum_{i=1}^{N}H^{m}\left(\tilde{Z}^{i}\left(\tau\left(t\right)\right)\right),\\
\frac{\mathrm{d}R^{N,\delta}}{\mathrm{d}t}= & \:L\left(\bar{u}^{N,\delta}\left(\tau\left(t\right)\right)\right)R^{N,\delta}\left(\tau\left(t\right)\right)+R^{N,\delta}\left(\tau\left(t\right)\right)L\left(\bar{u}^{N,\delta}\left(\tau\left(t\right)\right)\right)^{T}+Q_{\sigma}+\frac{1}{N}\sum_{i=1}^{N}H^{v}\left(\tilde{Z}^{i}\left(\tau\left(t\right)\right)\right),
\end{aligned}
\label{eq:model_stat_ens}
\end{equation}
where the forward Euler scheme is adopted with the discrete time update
using a constant $\tau\left(t\right)=n\Delta t$ during the time interval
$t\in\left[t_{n},t_{n+1}\right]$. Above in \eqref{eq:model_conti}
and \eqref{eq:model_stat_ens}, we neglect the last relaxation term
with $\epsilon$ since it will automatically vanish with the resulting
consistency. Notice that $\bar{u}^{N,\delta},R^{N,\delta}$ and $\bar{u},R$
are stochastic processes due to the random samples $\left\{ \tilde{Z}^{i}\right\} $
and the conditional expectation dependent on the observations $\mathcal{G}_{t}$.
We have the following result for the convergence of finite ensemble
$N$ and finite time step $\Delta t$ approximation to the continuous
model prediction.
\begin{thm}
\label{thm:conv_stat}If Assumption~\ref{assu:model_disc} is satisfied
and under the same initial condition, the statistical solution $\left(\bar{u}_{n}^{N,\delta},R_{n}^{N,\delta}\right)=\left(\bar{u}^{N,\delta}\left(t_{n}\right),R^{N,\delta}\left(t_{n}\right)\right)$
of the finite ensemble model \eqref{eq:model_stat_ens} with discrete
time step $\Delta t$ converges to the true statistical state $\left(\bar{u}_{n},R_{n}\right)=\left(\bar{u}\left(t_{n}\right),R\left(t_{n}\right)\right)$
of the continuous model \eqref{eq:model_conti} with the error estimates
\begin{equation}
\begin{aligned}\mathbb{E}\left[\sup_{n\Delta t\leq T}\left|\bar{u}_{n}^{N,\delta}-\bar{u}_{n}\right|^{2}\right] & \leq\left(C_{1,T}\Delta t^ {}+\frac{C_{2,T}}{N}\right)\left\Vert H^{m}\right\Vert _{\infty},\\
\mathbb{E}\left[\sup_{n\Delta t\leq T}\left\Vert R_{n}^{N,\delta}-R_{n}\right\Vert ^{2}\right] & \leq\left(C_{1,T}^{\prime}\Delta t^ {}+\frac{C_{2,T}^{\prime}}{N}\right)\left(\left\Vert H^{m}\right\Vert _{\infty}+\left\Vert H^{v}\right\Vert _{\infty}\right),
\end{aligned}
\label{eq:stat_bnds}
\end{equation}
where $C_{1,T},C_{2,T},C_{1,T}^{\prime},C_{2,T}^{\prime}$ are constants
depending on the final time $T$.
\end{thm}

\begin{proof}
First, considering the mean equations in \eqref{eq:model_conti} and
\eqref{eq:model_stat_ens} from the same initial state, we have
\[
\bar{u}^{N,\delta}\left(t\right)-\bar{u}\left(t\right)=\int_{0}^{t}\left[M\left(\bar{u}^{N,\delta}\left(\tau\left(s\right)\right)\right)-M\left(\bar{u}\left(s\right)\right)\right]\mathrm{d}s+\int_{0}^{t}\left[\left\langle \tilde{\rho}_{\tau\left(s\right)}^{N},H^{m}\right\rangle -\left\langle \tilde{\rho}_{s},H^{m}\right\rangle \right]\mathrm{d}s,
\]
where we define $M\left(u\right)=\Lambda u+B\left(u,u\right)+F$ and
assume that the forcing $F$ and $Q_{\sigma}$ are constants for simplicity.
Using the Lipschitz condition for $M$ from Assumption \ref{assu:model_disc}
and applying Cauchy-Schwarz inequality, there is
\begin{align}
\mathbb{E}\left[\sup_{t\leq T}\left|\bar{u}^{N}\left(t\right)-\bar{u}\left(t\right)\right|^{2}\right] & \leq2T\beta^{2}\mathbb{E}\int_{0}^{T}\left|\bar{u}^{N,\delta}\left(\tau\left(s\right)\right)-\bar{u}\left(s\right)\right|^{2}\mathrm{d}s+2T\mathbb{E}\int_{0}^{T}\left|\left\langle \tilde{\rho}_{\tau\left(s\right)}^{N},H^{m}\right\rangle -\left\langle \tilde{\rho}_{s},H^{m}\right\rangle \right|^{2}\mathrm{d}s\nonumber \\
 & \leq C_{1}T\int_{0}^{T}\mathbb{E}\left[\sup_{s^{\prime}\leq s}\left|\bar{u}^{N,\delta}\left(s^{\prime}\right)-\bar{u}\left(s^{\prime}\right)\right|^{2}\right]\mathrm{d}s+C_{2}T^{2}\Delta t\left\Vert H^{m}\right\Vert _{\infty}^{2}+\frac{C_{3}T^{2}}{N}\left\Vert H^{m}\right\Vert _{\infty}^{2}.\label{eq:est_m}
\end{align}
Above in the first term of the last inequality, we estimate the error
by comparing the discretized time solution $\bar{u}^{N,\delta}\left(\tau\left(t\right)\right)$
following \eqref{eq:model_stat_ens} with the corresponding continuous
time solution $\bar{u}^{N,\delta}\left(t\right)$ 
\[
\left|\bar{u}^{N,\delta}\left(\tau\left(t\right)\right)-\bar{u}^{N,\delta}\left(t\right)\right|^{2}\leq\left|t-\tau\left(t\right)\right|^{2}\left[\left|M\left(\bar{u}^{N,\delta}\left(\tau\left(s\right)\right)\right)\right|^{2}+\left|\left\langle \tilde{\rho}_{\tau\left(t\right)}^{N},H^{m}\right\rangle \right|^{2}\right]\leq\Delta t^{2}\left(\left\Vert M\right\Vert _{\infty}^{2}+\left\Vert H^{m}\right\Vert _{\infty}^{2}\right).
\]
Thus the error estimation follows
\[
\left|\bar{u}^{N,\delta}\left(\tau\left(s\right)\right)-\bar{u}\left(s\right)\right|^{2}\leq2\left|\bar{u}^{N,\delta}\left(\tau\left(s\right)\right)-\bar{u}^{N,\delta}\left(s\right)\right|^{2}+2\left|\bar{u}^{N,\delta}\left(s\right)-\bar{u}\left(s\right)\right|^{2}\leq C\Delta t^{2}\left\Vert H^{m}\right\Vert _{\infty}^{2}+2\sup_{s^{\prime}\leq s}\left|\bar{u}^{N,\delta}\left(s^{\prime}\right)-\bar{u}\left(s^{\prime}\right)\right|^{2}.
\]
And for the second term involving expectation of $H^{m}$, the convergence
of the empirical measure \eqref{eq:conv_pdf} gives
\begin{align*}
\mathbb{E}\left[\left|\left\langle \tilde{\rho}_{\tau\left(t\right)}^{N},H^{m}\right\rangle -\left\langle \tilde{\rho}_{t},H^{m}\right\rangle \right|^{2}\right] & \leq2\mathbb{E}\left[\left|\left\langle \tilde{\rho}_{\tau\left(t\right)}^{N},H^{m}\right\rangle -\left\langle \tilde{\rho}_{t}^{N},H^{m}\right\rangle \right|^{2}\right]+2\mathbb{E}\left[\left|\left\langle \tilde{\rho}_{t}^{N},H^{m}\right\rangle -\left\langle \tilde{\rho}_{t},H^{m}\right\rangle \right|^{2}\right]\\
 & \leq\left(C\Delta t+\frac{C_{T}}{N}\right)\left\Vert H^{m}\right\Vert _{\infty}^{2}.
\end{align*}
Finally applying Gr\"{o}nwall's inequality to \eqref{eq:est_m},
we get the mean state estimate in \eqref{eq:stat_bnds}.

Next, under a similar fashion, we have
\begin{align*}
R^{N,\delta}\left(t\right)-R\left(t\right)= & \int_{0}^{t}L\left(\bar{u}\left(s\right)\right)\left[R^{N,\delta}\left(\tau\left(s\right)\right)-R\left(s\right)\right]\mathrm{d}s+\int_{0}^{t}\left[L\left(\bar{u}^{N,\delta}\left(\tau\left(s\right)\right)\right)-L\left(\bar{u}\left(s\right)\right)\right]R\left(s\right)\mathrm{d}s\\
 & +\int_{0}^{t}\left[L\left(\bar{u}^{N,\delta}\left(\tau\left(s\right)\right)\right)-L\left(\bar{u}\left(s\right)\right)\right]\left[R^{N,\delta}\left(\tau\left(s\right)\right)-R\left(s\right)\right]\mathrm{d}s+c.c.\\
 & +\int_{0}^{t}\left[\left\langle \tilde{\rho}_{\tau\left(s\right)}^{N},H^{v}\right\rangle -\left\langle \tilde{\rho}_{s},H^{v}\right\rangle \right]\mathrm{d}s.
\end{align*}
Above, $c.c.$ represents the symmetric terms from the transposes
$RL\left(\bar{u}\right)^{T}$. Again, using the Lipschitz condition
of $L$ in Assumption \ref{assu:model_disc}, $\left\Vert L\left(u\right)\right\Vert \leq\beta\left|u\right|+\beta_{1}$,
we can compute errors from the covariance equation 
\begin{align*}
\mathbb{E}\left[\sup_{t\leq T}\left\Vert R^{N,\delta}\left(t\right)-R\left(t\right)\right\Vert ^{2}\right] & \leq C_{1}T\beta^{2}\mathbb{E}\left[\sup_{t\leq T}\left|\bar{u}\right|^{2}\int_{0}^{T}\left\Vert R^{N,\delta}\left(\tau\left(s\right)\right)-R\left(s\right)\right\Vert ^{2}\mathrm{d}s\right]\\
 & +C_{2}T\beta^{2}\mathbb{E}\left[\sup_{t\leq T}\left|\bar{u}^{N,\delta}\left(t\right)-\bar{u}\left(t\right)\right|^{2}\sup_{t\leq T}\left\Vert R\left(t\right)\right\Vert ^{2}\right]\\
 & +C_{3}\mathbb{E}\int_{0}^{T}\left|\left\langle \tilde{\rho}_{\tau\left(s\right)}^{N},H^{v}\right\rangle -\left\langle \tilde{\rho}_{s},H^{v}\right\rangle \right|^{2}\mathrm{d}s+C_{T}\Delta t^{2}.
\end{align*}
Using the uniform boundedness of $\bar{u},R$ and \eqref{eq:conv_pdf}
for $H^{v}$ together with the previous error estimate of the mean
state for $\mathbb{E}\left[\sup_{t\leq T}\left|\bar{u}^{N,\delta}\left(t\right)-\bar{u}\left(t\right)\right|^{2}\right]$,
we reach the final covariance error estimate in \eqref{eq:stat_bnds}.
\end{proof}
Theorem~ \ref{thm:conv_stat} guarantees that the discrete numerical
scheme of the approximating ensemble filter model can recover the
leading-order statistics in mean and covariance. It implies that the
performance of the ensemble filter estimation relies on the accurate
approximations of the expectation of the observation functions $H^{m},H^{v}$.
Usually, the higher-order moments in $H^{m}$ and $H^{v}$ in \eqref{eq:obs_func}
become extremely difficult to capture with a small sample size. This
leads to the rapidly growing model errors (as in Figure~\ref{fig:Statistical-forecasts}
shown in the numerical tests). On the other hand, the design of the
new high-order filtering scheme guarantees consistent statistics $\tilde{\mathbb{E}}\left[H\left(\tilde{Z}_{t}\right)\right]$
in the observation functions with the optimal filter solution $\hat{\rho}$
(see Theorem 7 in \cite{qi2024coupled}). This leads to the much improved
performance of the the statistical forecasts using the new filter
model.

\section{Numerical performance on the prototype triad system\protect\label{sec:Numerical-performance}}

Using the explicit equations \eqref{eq:filter_scheme}, we now test
the performance of the proposed data assimilation algorithm on a prototype
triad system with instructive implications to many practical applications.
The triad system \cite{majda2018strategies} is given by a three-dimensional
ODE system for the state $\mathbf{u}=\left(u_{1},u_{2},u_{3}\right)^{T}$
with both linear and nonlinear coupling combined with stochastic forcing
\begin{equation}
\begin{aligned}\frac{du_{1}}{dt}= & \lambda_{2}u_{3}-\lambda_{3}u_{2}-d_{1}u_{1}+B_{1}u_{2}u_{3}+\sigma_{1}\dot{W}_{1},\\
\frac{du_{2}}{dt}= & \lambda_{3}u_{1}-\lambda_{1}u_{3}-d_{2}u_{2}+B_{2}u_{3}u_{1}+\sigma_{2}\dot{W}_{2},\\
\frac{du_{3}}{dt}= & \lambda_{1}u_{2}-\lambda_{2}u_{1}-d_{3}u_{3}+B_{3}u_{1}u_{2}+\sigma_{3}\dot{W_{3}}.
\end{aligned}
\label{eq:triad}
\end{equation}
It can be seen that the above triad system \eqref{eq:triad} fits
into our general formulation \eqref{eq:model_general}, where the
model coefficients 
\[
\Lambda=\begin{bmatrix}-d_{1} & -\lambda_{3} & \lambda_{2}\\
\lambda_{3} & -d_{2} & -\lambda_{1}\\
-\lambda_{2} & \lambda_{1} & -d_{3}
\end{bmatrix},\quad B\left(\mathbf{u},\mathbf{u}\right)=\begin{bmatrix}B_{1}u_{2}u_{3}\\
B_{2}u_{3}u_{1}\\
B_{3}u_{1}u_{2}
\end{bmatrix},
\]
contain the linear skew-symmetric (off-diagonal) and dissipation (diagonal)
operator, together with the nonlinear quadratic coupling $B\left(\mathbf{u},\mathbf{u}\right)$
satisfying energy conservation with $B_{1}+B_{2}+B_{3}=0$. The triad
system can serve as an elementary building block of many more general
turbulent systems emphasizing the key energy conserving nonlinear
interactions. Though low-dimensional, this system can demonstrate
a wide variety of different statistical regimes (as shown next in
Figure~\ref{fig:stat_pdf}), making it a nice first test model for
a thorough study of the prediction skill of the proposal ensemble
filtering strategy in dealing with different statistical features.

\subsection{Typical statistical regimes in the triad system\protect\label{subsec:Typical-statistical-regimes}}

One attractive feature of the triad system \eqref{eq:triad} as a
prototype test model is that it is able to generate a wide variety
of dynamical regimes demonstrating distinctive statistical features
ranging from near-Gaussian to highly non-Gaussian probability distributions.
This sets up a desirable testbed for examining the skill of different
statistical prediction methods in dealing with vastly different statistical
dynamics.

\subsubsection{Statistical regimes with distinctive statistical features}

The interaction in triad systems \eqref{eq:triad} constitutes the
generic linear and nonlinear coupling mechanism between any three
modes in larger systems with quadratic nonlinearity. A direct three-dimensional
Galerkin truncation of many complex turbulent dynamics possesses the
energy-conserving nonlinearity as in the general formulation \eqref{eq:model_general}.
For example , a direct link can be built to interpret the triad system
as a prototype three-mode interaction with forward and backward energy
cascades in geophysical turbulence (see Appendix~\ref{subsec:A-direct-link}).
The random forcing together with the damping term simulates the inhomogeneous
effects of the interaction with other modes that are not resolved
in the projected three dimensional subspace. Thus, the stochastic
triad system can serve as a qualitative model for a wide variety of
turbulent phenomena regarding energy exchange and cascades and supply
important intuition for many phenomena \cite{frisch1995turbulence,gluhovsky1997interpretation,majda2002priori}.
They also provide elementary test models with subtle features for
prediction, uncertainty quantification, and state estimation. Additional
dynamical and statistical properties of the triad system are summarized
in Appendix~\ref{sec:Details-on-triad}.

In our testing cases, we consider the following three typical dynamical
regimes of the triad system \eqref{eq:triad} containing representative
statistical structures. Model parameters used for the three test regimes
are listed in Table~\ref{tab:model_arameters}. 
\begin{itemize}
\item \emph{Regime I: Near-Gaussian regime with equipartition of energy.}
This regime considers the convergence to a Gaussian equilibrium distribution
with the competition of linear and nonlinear effects. The equipartition
of energy, that is, $\frac{\sigma_{1}^{2}}{2d_{1}}=\frac{\sigma_{2}^{2}}{2d_{2}}=\frac{\sigma_{3}^{2}}{2d_{3}}=\sigma_{\mathrm{eq}}^{2}$,
is designed so that a Gaussian distribution, $p_{\mathrm{eq}}\sim\exp\left(-\frac{1}{2}\sigma_{\mathrm{eq}}^{-2}\mathbf{\left|u\right|}^{2}\right)$,
will be reached at the final equilibrium state. The linear and nonlinear
parameters are chosen to have comparable values to induce strong interactions
during the transient state;
\item \emph{Regime II: Nonlinear regime with forward energy cascade.} This
regime focuses on strong quadratic coupling with weak linear damping
and forcing effects. Skew-symmetric linear terms are set to be zero
and only small damping and noise effects are added. The first mode
$u_{1}$ is set to have large initial mean and covariance while the
other two modes $u_{2},u_{3}$ only have small initial values. This
induces strong energy cascades from $u_{1}$ to the other two less
energetic modes $u_{2},u_{3}$ driven by the dominant nonlinear coupling;
\item \emph{Regime III: Unstable regime with dual energy cascades.} This
regime is used to simulate the inherent internal instability observed
in turbulent systems. The instability is introduced by a negative
damping in the first mode $u_{1}$, while the other two modes $u_{2},u_{3}$
are stable with positive damping. On the other hand, the first mode
is weakly forced by stochastic forcing while the other two are strongly
excited by random noises. The nonlinear coupling first makes that
energy cascades forwardly from mode $u_{1}$ to the other less energetic
modes $u_{2},u_{3}$ then backwardly from the excited modes $u_{2},u_{3}$
back to $u_{1}$. 
\end{itemize}
The initial state $\mathbf{u}_{0}\sim\mathcal{N}\left(\bar{\mathbf{u}}_{0},\mathbf{r}_{0}\right)$
is set to satisfy an independent Gaussian distribution with mean $\bar{\mathbf{u}}_{0}$
and variance $\mathbf{r}_{0}$. The true statistical solutions of
the triad system \eqref{eq:triad} in the above dynamical regimes
are solved through direct Monte-Carlo simulations. We run an ensemble
of $\mathrm{MC}=1\times10^{5}$ particles, which shall be enough for
capturing the statistics in a three-dimensional phase space. A fourth-order
Runge-Kutta scheme with time step $\Delta t=1\times10^{-3}$ is used
to integrate the system in time. The stochastic forcing is simulated
through the standard Euler-Maruyama scheme. The initial ensemble is
chosen from a standard Gaussian random sampling with the mean $\bar{\mathbf{u}}_{0}$
and variance $\mathbf{r}_{0}$ listed in Table~\ref{tab:model_arameters}.
In particular, we choose $B_{1}>0$ and $B_{2},B_{3}<0$ to induce
nonlinear instability (see the stability analysis in \eqref{eq:triad_stability}).
The model is run up to a final time $T=10$ where near equilibrium
state is reached.

\begin{table}
\begin{centering}
{\footnotesize{}%
\begin{tabular}{ccccccc}
\toprule 
 & {\footnotesize$\left(B_{1},B_{2},B_{3}\right)$} & {\footnotesize$\left(\lambda_{1},\lambda_{2},\lambda_{3}\right)$} & {\footnotesize$\left(d_{1},d_{2},d_{3}\right)$} & {\footnotesize$\left(\sigma_{1},\sigma_{2},\sigma_{3}\right)$} & {\footnotesize$\bar{\mathbf{u}}_{0}$} & {\footnotesize$\mathbf{r}_{0}$}\tabularnewline
\midrule
\midrule 
{\footnotesize regime I} & {\footnotesize$\left(1,-0.6,-0.4\right)$} & {\footnotesize$\left(3,-2,-1\right)$} & {\footnotesize$\left(0.2,0.1,0.1\right)$} & {\footnotesize$\left(1.58,1.12,1.12\right)$} & {\footnotesize$\left(2,1.6,-2\right)$} & {\footnotesize$\left(0.5,0.5,1\right)$}\tabularnewline
\midrule 
{\footnotesize regime II} & {\footnotesize$\left(1,-0.6,-0.4\right)$} & {\footnotesize$\left(0,0,0\right)$} & {\footnotesize$\left(0.02,0.01,0.01\right)$} & {\footnotesize$\left(0.5,0.35,0.35\right)$} & {\footnotesize$\left(3,-0.1,0.1\right)$} & {\footnotesize$\left(0.5,0.01,0.01\right)$}\tabularnewline
\midrule 
{\footnotesize regime III} & {\footnotesize$\left(2,-1,-1\right)$} & {\footnotesize$\left(0.09,0.06,-0.03\right)$} & {\footnotesize$\left(-0.4,2,2\right)$} & {\footnotesize$\left(0.1,0.32,0.32\right)$} & {\footnotesize$\left(2,1,1.5\right)$} & {\footnotesize$\left(0.5,5,10\right)$}\tabularnewline
\bottomrule
\end{tabular}}{\footnotesize\par}
\par\end{centering}
\caption{Parameters for the triad system \eqref{eq:triad} in the three test
regimes.\protect\label{tab:model_arameters}}
\end{table}

\begin{figure}
\begin{centering}
\subfloat[Regime I]{\includegraphics[scale=0.36]{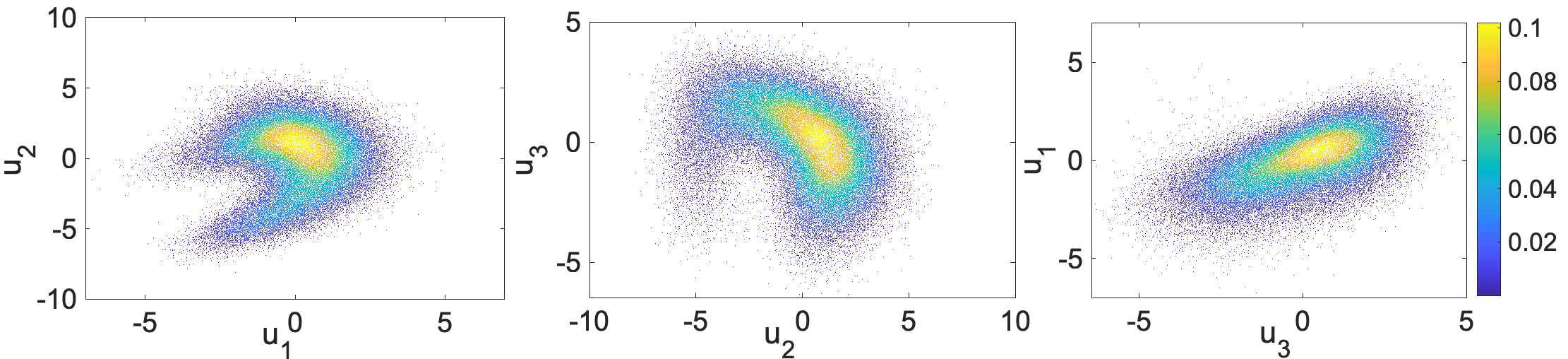}

}
\par\end{centering}
\begin{centering}
\vspace{-1.2em}\subfloat[Regime II]{\includegraphics[scale=0.36]{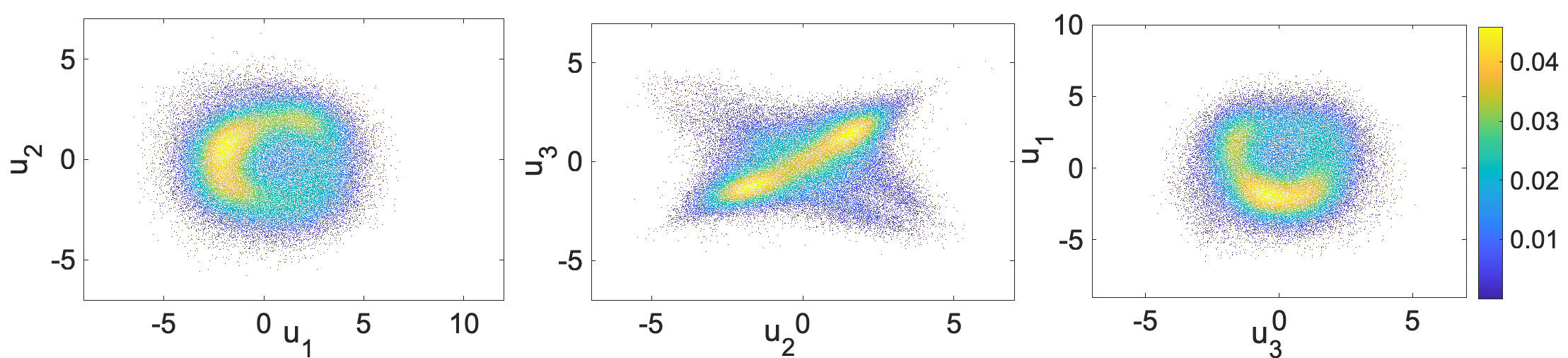}

}
\par\end{centering}
\begin{centering}
\vspace{-1.2em}\subfloat[Regime III]{\includegraphics[scale=0.36]{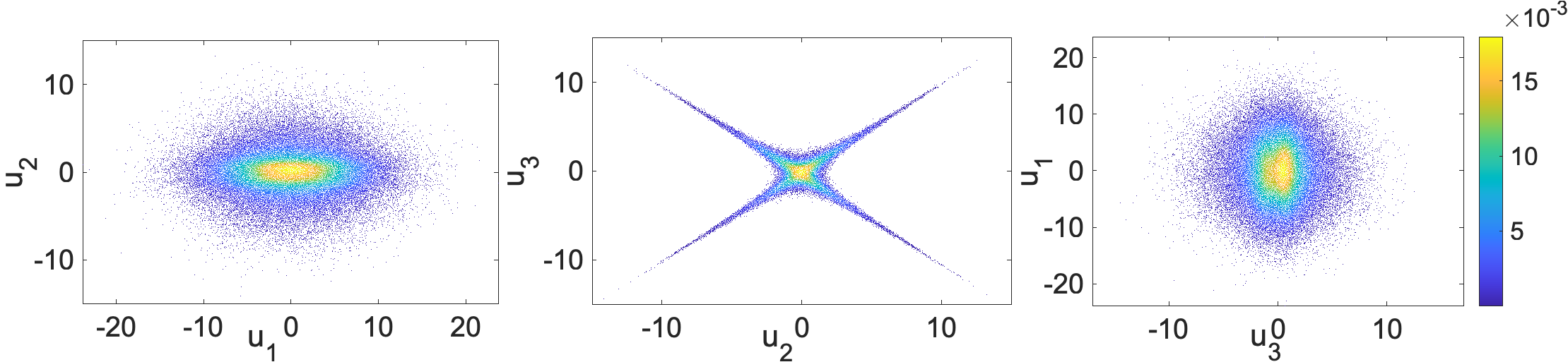}

}
\par\end{centering}
\caption{Joint PDFs at $t=5$ of triad modes $u_{1},u_{2},u_{3}$ in the three
test regimes shown in scatter plots from a direct MC simulation using
$\mathrm{MC}=1\times10^{5}$ samples. The density of particles is
represented by colors in the scatter plots. \protect\label{fig:stat_pdf}}
\end{figure}

The projected probability distributions of the triad state $p\left(\mathbf{u},t\right)$
captured by the direct Monte-Carlo simulations are demonstrated in
Figure~\ref{fig:stat_pdf}. Representative non-Gaussian probability
distributions are observed among all test regimes with distinctive
statistical structures. The first test regime is the simplest but
nevertheless representative showing the route of transient convergence
to equipartition of energy. Still, as we will show in the following
numerical tests in Section~\ref{subsec:Numerical-performance}, higher-order
moments are playing an pivoting role in this case and cannot be simply
ignored in determining the correct final near-Gaussian equilibrium
distribution. The second test regime emphasizes the nonlinear quadratic
coupling between the three modes, leading to more complicated non-Gaussian
probability distributions. In particular, we observe the wider spread
of the samples in the scatter plots representing extreme events that
are crucial but difficult to capture with a small ensemble. The third
test regime introduces stronger interactions and dual energy cascades
between the interacting modes. This leads to a strange attractor with
star-shaped joint-distribution, showing a strongly nonlinear regime
dominated by highly non-Gaussian statistics. Especially in this regime,
the negative damping $d_{1}=-0.4$ in the first mode introduces persistent
internal instability into the system. In addition, the skew-symmetric
linear interaction terms add extra emphasis on the cross-covariances.
This regime becomes especially interesting and challenging because
of the strong and inherent instability.

\subsubsection{Small ensemble prediction with the coupled stochastic-statistical
model\protect\label{subsec:Direct-prediction}}

We start with testing direct forecast of the coupled stochastic-statistical
model \eqref{eq:num_model} to capture key model statistics. Using
this coupled modeling framework, the statistical equations will be
used to compute the leading moments $\bar{u}^{N}$ and $R^{N}$, combined
with an ensemble simulation for the stochastic coefficients $\left\{ Z^{i}\right\} _{i=1}^{N}$
aiming to capture the high-order moments feedback. To cope with the
realistic scenario where only a small number of samples are affordable,
we check the model forecast skill using a moderate ensemble size $N=100$,
in contrast to the truth in the previous section generated by $\mathrm{MC}=1\times10^{5}$
samples. Due to the dominant nonlinear coupling terms, the higher-order
moments feedback are involved in the statistical equations \eqref{eq:dyn_stat},
requiring accurately capturing the non-Gaussian statistics from the
limited samples even only to predict the leading-order mean and covariance.
This sets an especially challenging problem demanding good characterization
of the non-Gaussian distributions (including the extreme outliers
observed in the PDFs in Figure~\ref{fig:stat_pdf}) using only the
small number of samples.

\begin{figure}
\includegraphics[scale=0.3]{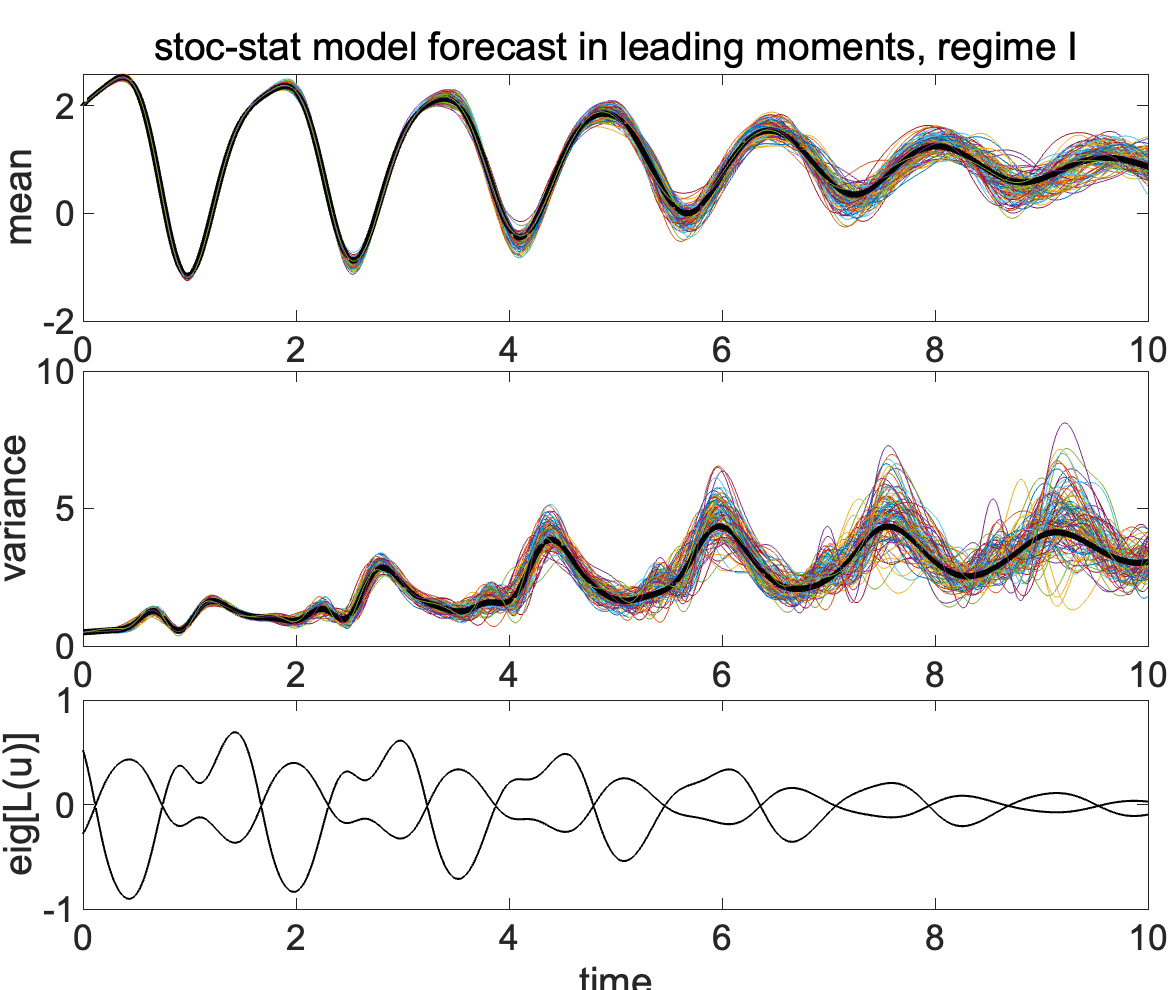}\includegraphics[scale=0.3]{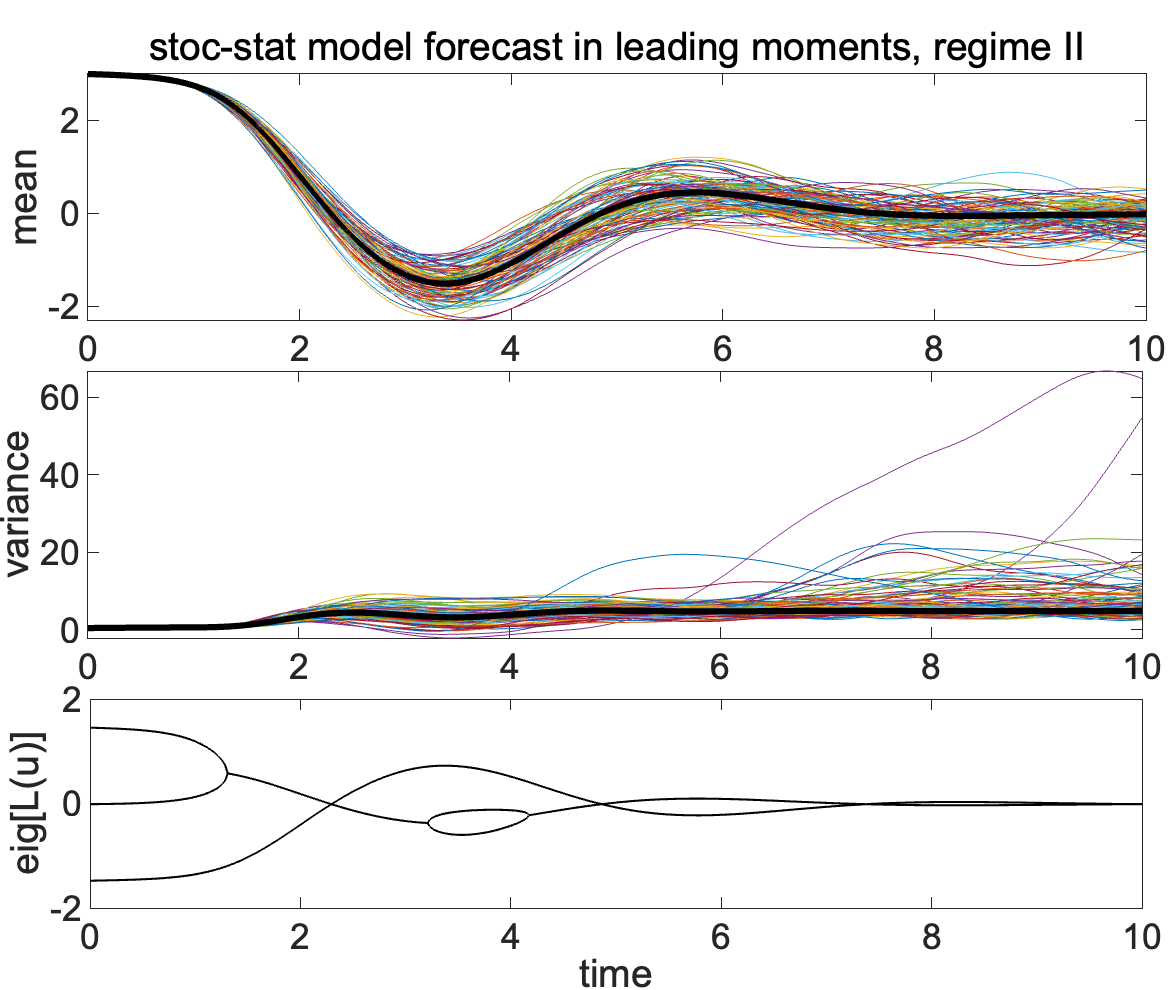}\includegraphics[scale=0.3]{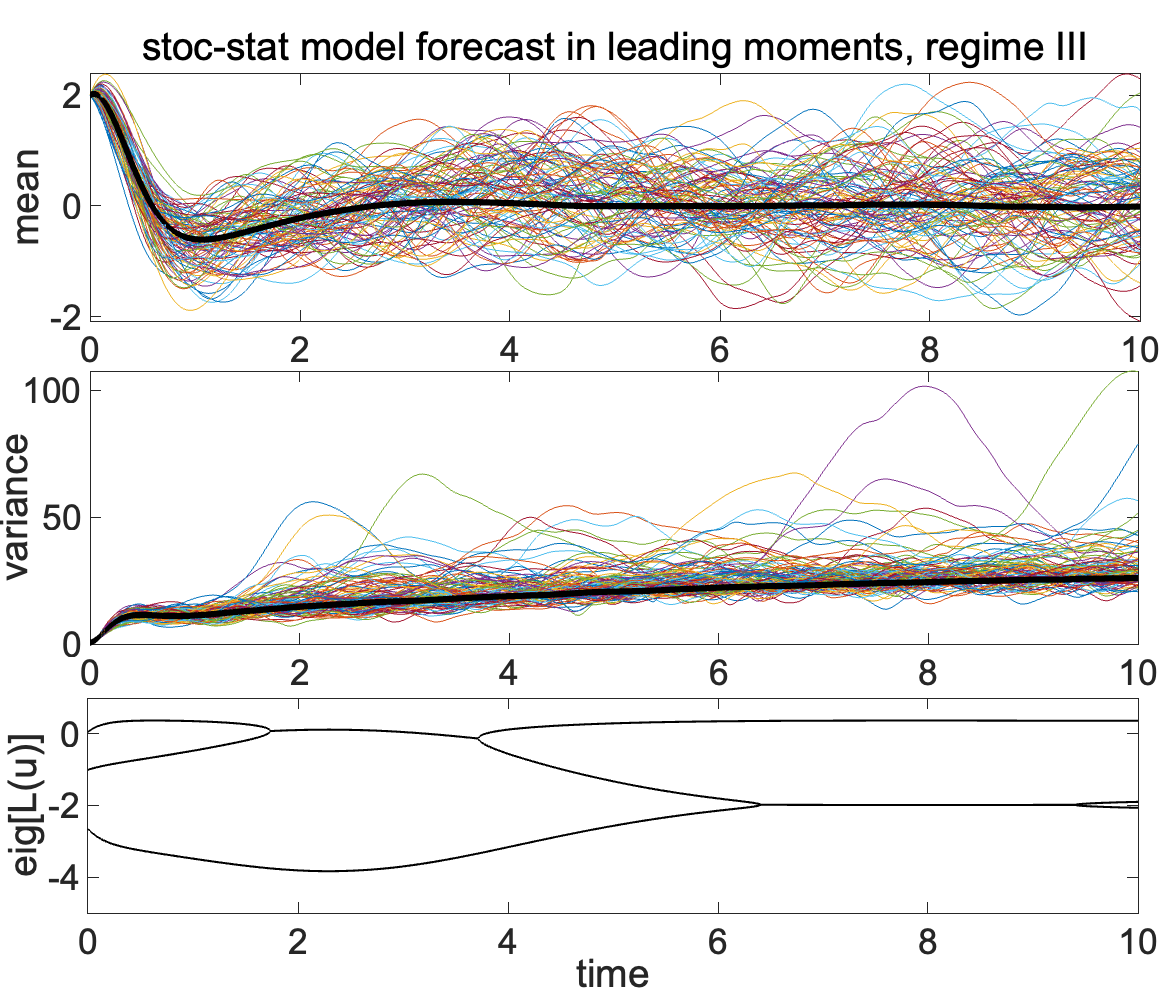}

\caption{Statistical forecasts using the stochastic-statistical model with
$N=100$ samples. Different realizations of the mean $\bar{u}_{1}$
and variance $r_{1}$ are plotted in comparison with the truth in
black lines. The third row plots the Lyapunov exponent of the system
indicating instability.\protect\label{fig:Statistical-forecasts}}

\end{figure}

First, the direct numerical predictions by running the coupled stochastic-statistical
model in the three test regimes are shown in Figure~\ref{fig:Statistical-forecasts}.
To demonstrate the unavoidable large amount uncertainty induced through
the small ensemble forecast, we plot multiple realizations of the
mean and variance trajectories $\left(\bar{u}^{N},R^{N}\right)$ from
different randomly sampled initial stochastic state $Z^{i}$ using
the small sample size $N=100$. Model errors in the mean and variance
forecasts are shown to rapidly grow in time among all three test regimes
starting from accurate initial states. To further illustrate the development
of such errors, we also plot the Lyapunov exponent, that is, the real
parts of eigenvalues of the linearized matrix $L\left(\bar{u}\right)$
in \eqref{eq:dyn_stat}, characterizing the inherent instability due
to the interaction with the mean state. Positive eigenvalues indicate
the unstable growth rate that amplifies small uncertainties in the
variance. It can be observed clearly that persistent instability maintains
in time amplifying the spread of different realizations of the solutions
in all three test regimes, especially in regimes II and III which
are experiencing stronger unstable growth during longer time periods.
In the exact equation \eqref{eq:dyn_stat}, such unstable growth rate
will be marginally balanced by the higher moments terms through the
nonlinear coupling between the states. However, due to the insufficient
representation of the highly non-Gaussian structures (as illustrated
later in Figure~\ref{fig:Sampling}) with the limited number of samples,
larger errors are introduced into the system. As a result, the trajectories
of the mean and variance fail to track the truth and quickly diverge
from the initial state. The accuracy of the predicted mean and variance
will improve if we increase the number of samples $N$ as indicated
in Theorem~\ref{thm:conv_stat}. However, this will usually require
an enormous sample size (due to the $T$ dependence in the coefficients
in \eqref{eq:stat_bnds}) even in this low-dimensional example, making
any direct numerical approach impractical. These examples offer a
typical illustration of the inherent difficulty in accurate prediction
of model statistics when only a small sample size is affordable.

\begin{figure}
\begin{centering}
\includegraphics[scale=0.3]{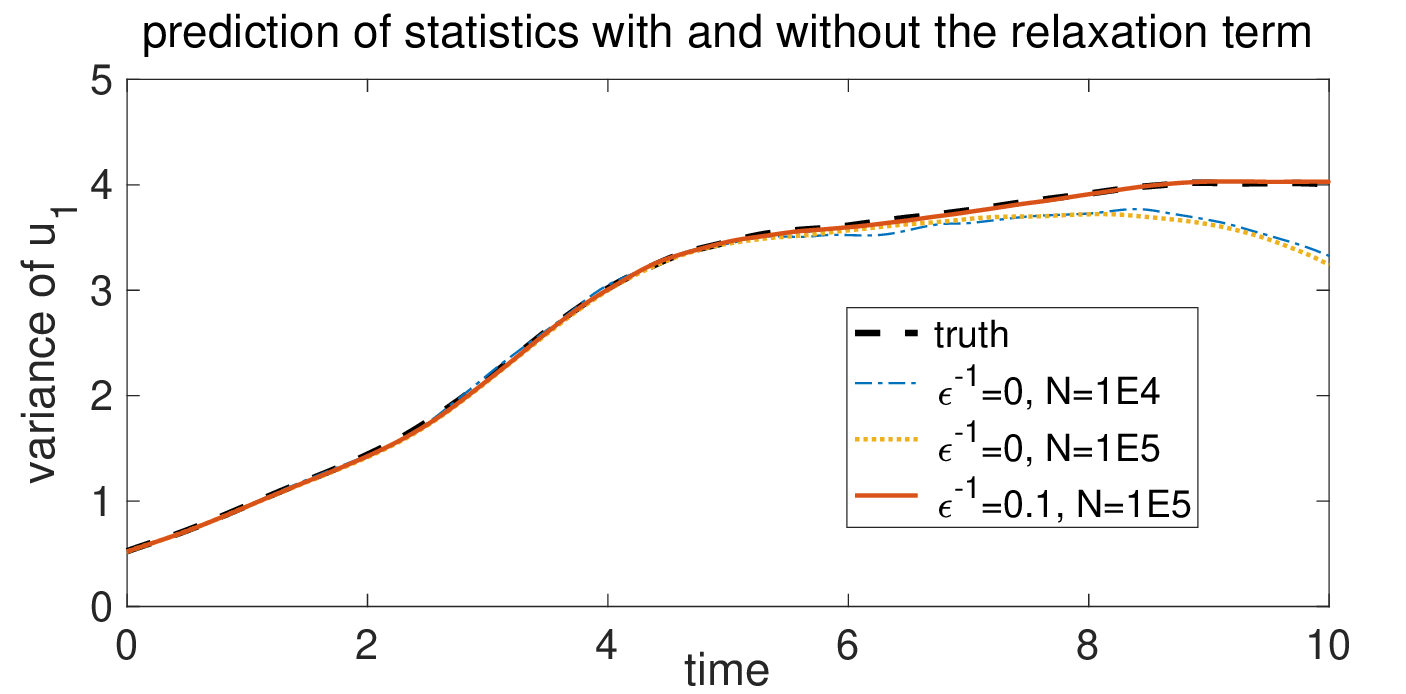}
\par\end{centering}
\caption{Model prediction of the variance in the most unstable mode $u_{1}$
with and without the additional relaxation term.\protect\label{fig:relaxation}}
\end{figure}

In addition, to further enforce the convergence of the numerical scheme,
we show that the additional relaxation term added in the covariance
equation of the numerical model \eqref{eq:num_model} is essential
especially in the regimes with stronger instability. In Figure~\ref{fig:relaxation},
we plot the model prediction of the variance in the most unstable
mode $u_{1}$ in regime III by directly applying the forecast model.
It shows that without the relaxation term $\epsilon^{-1}=0$ enforcing
the consistency between the sample approximation $\mathbb{E}^{N}\left[ZZ^{T}\right]$
and the covariance $R_{t}$, large numerical errors will start to
develop in time even using an extremely large sample size. This is
due to the persistent model instability among the states (see also
the last row of Figure~\ref{fig:Statistical-forecasts} for the internal
growth rate). On the other hand, it shows that the numerical errors
can be effectively corrected by just introducing a small relaxation
term with a small parameter $\epsilon^{-1}=0.1$, thus accurate statistical
convergence is guaranteed for the long-term time integration. In all
the following numerical tests, we adopt this small relaxation parameter
$\epsilon^{-1}=0.1$.

\subsection{Numerical performance using the data assimilation model\protect\label{subsec:Numerical-performance}}

Next, we demonstrate that the proposed data assimilation model can
effectively improve both stability and accuracy in the prediction
of the key statistics of the triad system. Furthermore, it shows that
maintaining the high-order correction terms constructed in the new
filtering equation \eqref{eq:filter_scheme} is essential to achieve
stable statistical prediction compared with the ensemble Kalman filters
\eqref{eq:EnKF} where only the low-order moments information is used.

\subsubsection{Calibration of observation noises\protect\label{subsec:Calibration-of-observation}}

\begin{figure}
\begin{centering}
\includegraphics[scale=0.3]{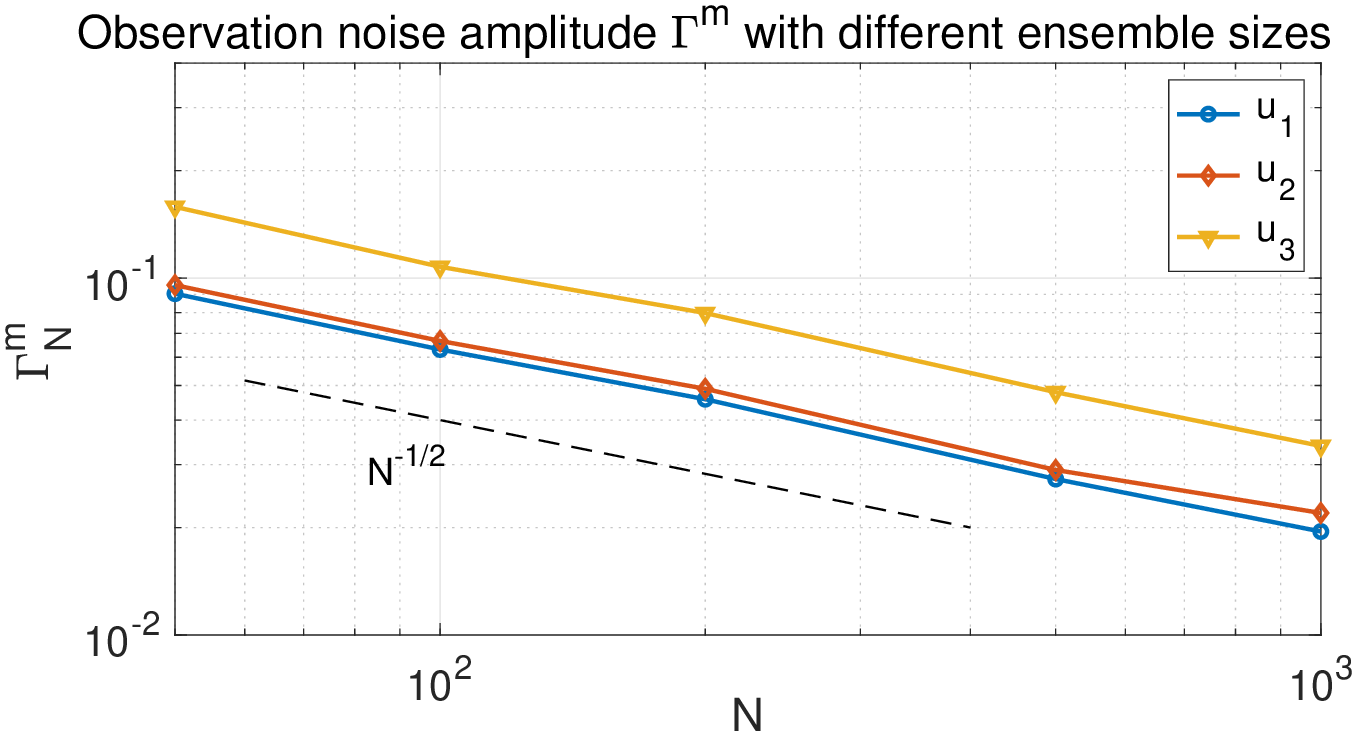}\includegraphics[scale=0.3]{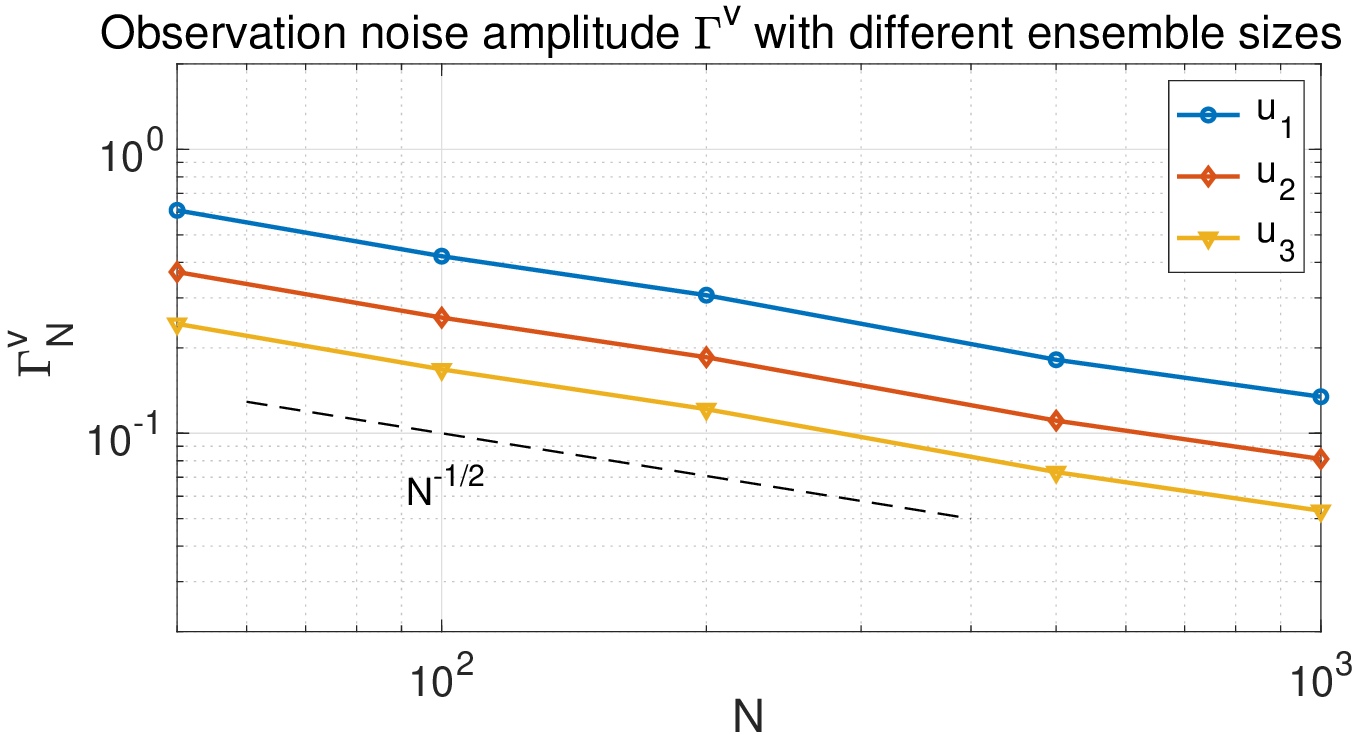}
\par\end{centering}
\caption{Estimate of the observation noise with different sample sizes $N$.
The noise parameters for the mean $\Gamma^{m}$ and covariance $\Gamma^{v}$
are computed based on the three modes of the triad system.\protect\label{fig:Estimate-obs-noise}}

\end{figure}

In setting up the filtering equations, we need to first estimate the
observation noises $\Gamma_{m}^{N}$ and $\Gamma_{v}^{N}$ in \eqref{eq:filter_dyn}
based on the finite ensemble size $N$. From the direct model simulations
for $\bar{u}^{N}$ and $R^{N}$ in Figure~\ref{fig:Statistical-forecasts},
it shows that it is reasonable to treat $y_{t}^{N}=\left(\bar{u}^{N},R^{N}\right)$
as a stochastic process and the randomness generated from the errors
in the finite ensemble approximation in \eqref{eq:obs_approx}, such
that the empirical estimate $\mathbb{E}^{N}H\mathrm{d}t=\mathbb{E}H\mathrm{d}t+\Gamma^{N}\mathrm{d}B$.
In general, we can only expect upper bounds for the errors in the
empirical averages as in \eqref{eq:stat_bnds} regarding to the sample
size $N$ and observation function $H$. Still in practical implementations,
it is sufficient to get an estimate of the noise levels of $\Gamma_{m}^{N}$
and $\Gamma_{v}^{N}$. In particular, we propose the following equations
for the observation states according to \eqref{eq:obs_general} where
we assume that error from the finite sample is dominant in the observation
equation, that is
\[
\begin{aligned}\mathrm{d}y_{t} & =\mathbb{E}H\mathrm{d}t,\\
\mathrm{d}y_{t}^{N} & =\mathbb{E}^{N}H=\mathbb{E}H\mathrm{d}t+\Gamma_{s}^{N}\mathrm{d}B.
\end{aligned}
\]
Above, $y_{t}$ is the true deterministic observed states $\bar{u},R$,
and $y_{t}^{N}$ is the stochastic observation process modeled with
an additional noise term accounting for the randomness with finite
ensemble estimation. Therefore, we find the following way to estimate
the observation errors by assuming that the noise amplitude remains
a constant in time for simplicity
\begin{equation}
\mathbb{E}\left(\left\Vert y_{t}^{N}-y_{t}\right\Vert ^{2}\right)\approx\int_{0}^{t}\left(\Gamma_{s}^{N}\right)^{2}\mathrm{d}s\approx t\left(\Gamma^{N}\right)^{2}.\label{eq:noise_estimate}
\end{equation}
In Figure~\ref{fig:Estimate-obs-noise}, we plot the estimated noise
amplitudes in the observed states of mean and variance using \eqref{eq:noise_estimate}
using different sample sizes $N$. The numerical results confirm the
$N^{-\frac{1}{2}}$ convergence rate \eqref{eq:stat_bnds} found in
Theorem~\ref{thm:conv_stat} depending on the sample size $N$. It
also provides a systematic way to estimate the observation noise level
in different filtering model simulations according to the scaling
law without repeating the different ensemble simulations many times. 

\subsubsection{Prediction of the statistical mean and covariance}

Now, we compare the performance of the high-order data assimilation
model in the triad system. To be adaptive to the general high-dimensional
systems, we focus on testing the forecast skill of the models using
a small sample size. As illustrated in Figure~\ref{fig:Statistical-forecasts},
this small sample size has already become insufficient to capture
the key statistical features in the simple triad system, and leads
to large fluctuating errors in the prediction of mean and covariance
by directly running the forecast model without using filtering. The
additional observation data $y_{n}=\left\{ \bar{u}_{n},R_{n}\right\} $
is then introduced aiming to correct the errors in the finite ensemble
forecast of the stochastic coefficients $\tilde{\mathbf{Z}}_{n}=\left\{ \tilde{Z}_{n}^{i}\right\} $
in \eqref{eq:filtering_full}. The goal is to generate representative
samples $\tilde{\mathbf{Z}}_{n}$ that can accurately characterize
the high-order moments and PDF structure of model states. Notice that
in the filter updating scheme only the information of derivatives
$\mathrm{d}\bar{u},\mathrm{d}R$ of leading-order moments are used
from the observation data for the updates of stochastic samples. The
model forecasts of the mean and covariance $\bar{u}^{N}$ and $R^{N}$
are still directly updated through the statistical equations \eqref{eq:num_model},
thus the observed mean and covariance are not involved in updating
the forecast mean and covariance. Therefore, their accuracy closely
replies on the the finite sample estimates of the higher-order moments
due to their nonlinear dynamics. In the following, we check the prediction
of mean and covariance using different models as an indicator for
the model skill to capture key high-order statistics in $\mathbf{\tilde{Z}}_{n}$.

\begin{figure}
\includegraphics[scale=0.3]{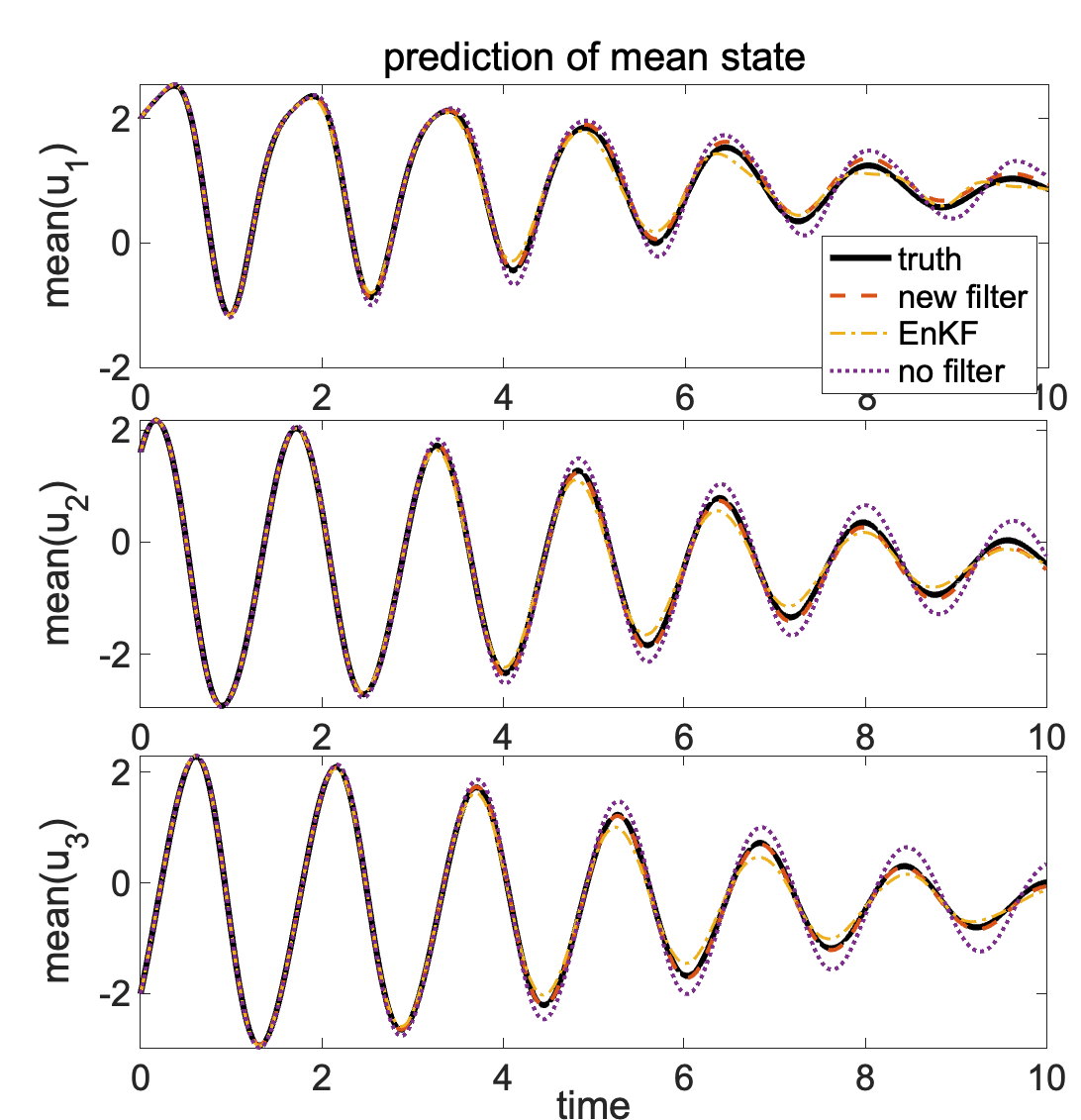}\includegraphics[scale=0.3]{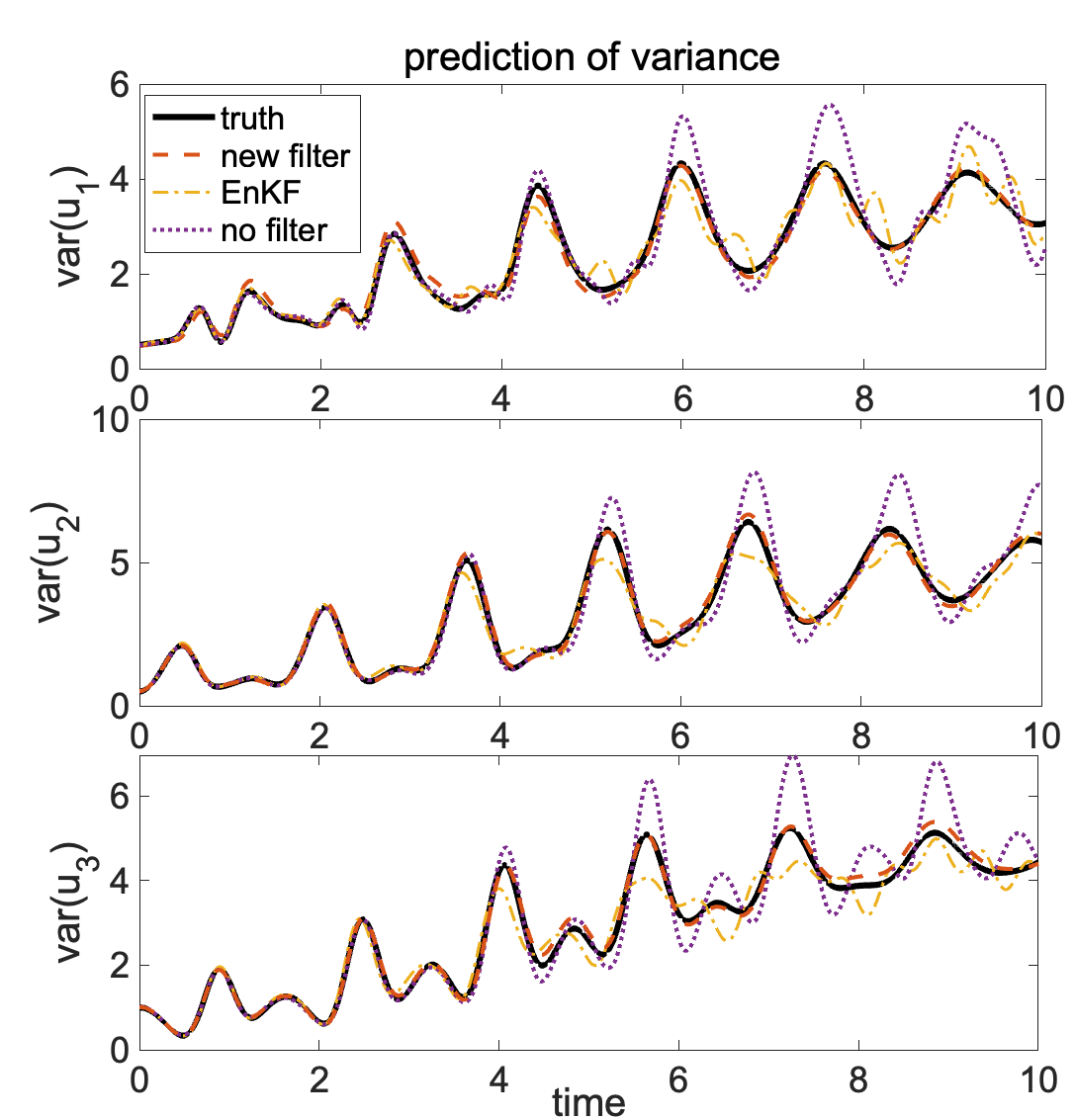}\includegraphics[scale=0.3]{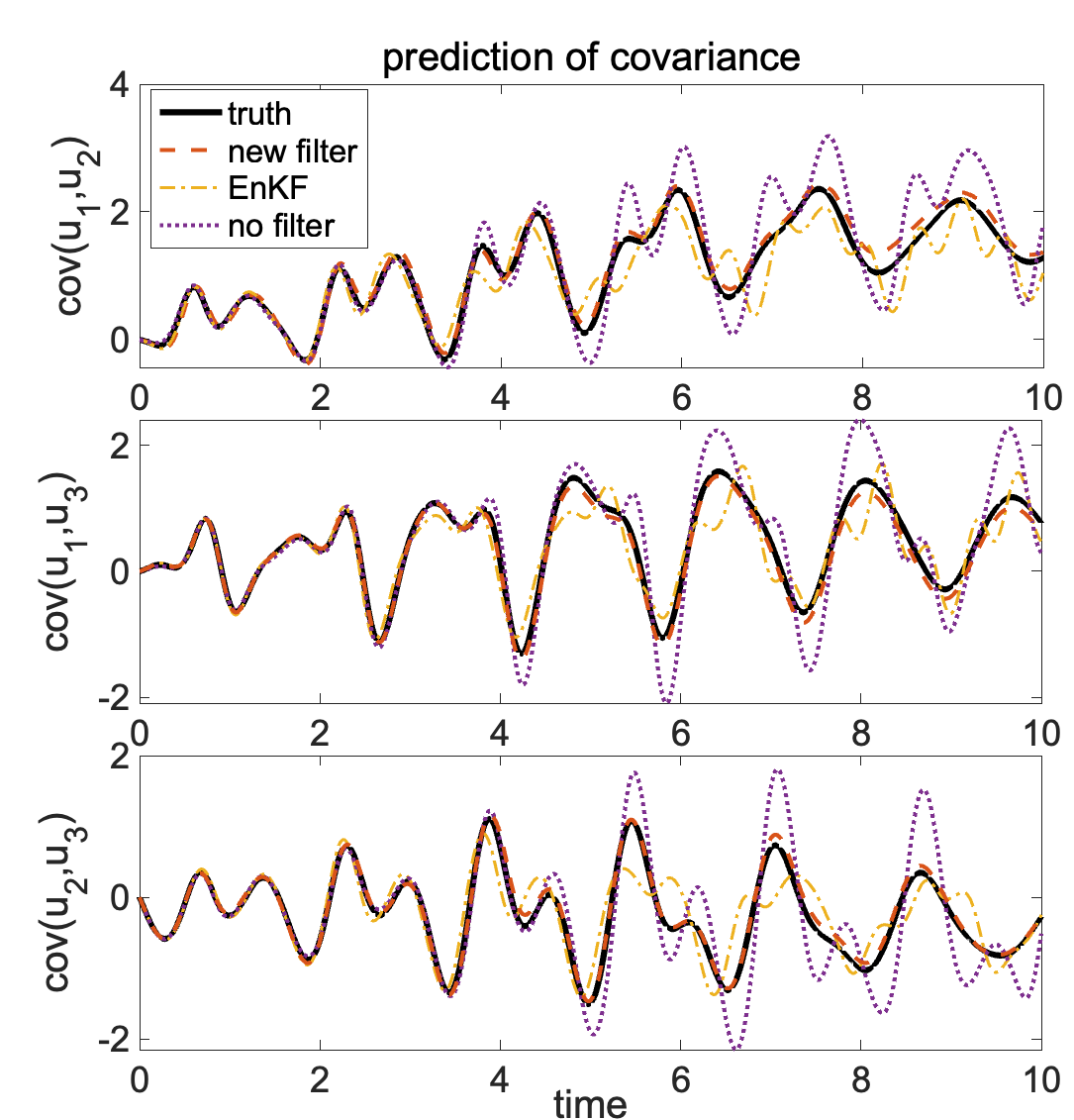}

\caption{Statistical prediction of the mean, variance, and covariance in regime
I of the triad system. Results from the high-order data assimilation
model are compared with the direct prediction without filter \eqref{eq:num_model}
and the standard EnKF \eqref{eq:EnKF} with $N=100$ samples. The
truth is generated with a direct MC approach with $\mathrm{MC}=1\times10^{5}$
samples.\protect\label{fig:Statistical-prediction1}}

\end{figure}

\begin{figure}
\includegraphics[scale=0.3]{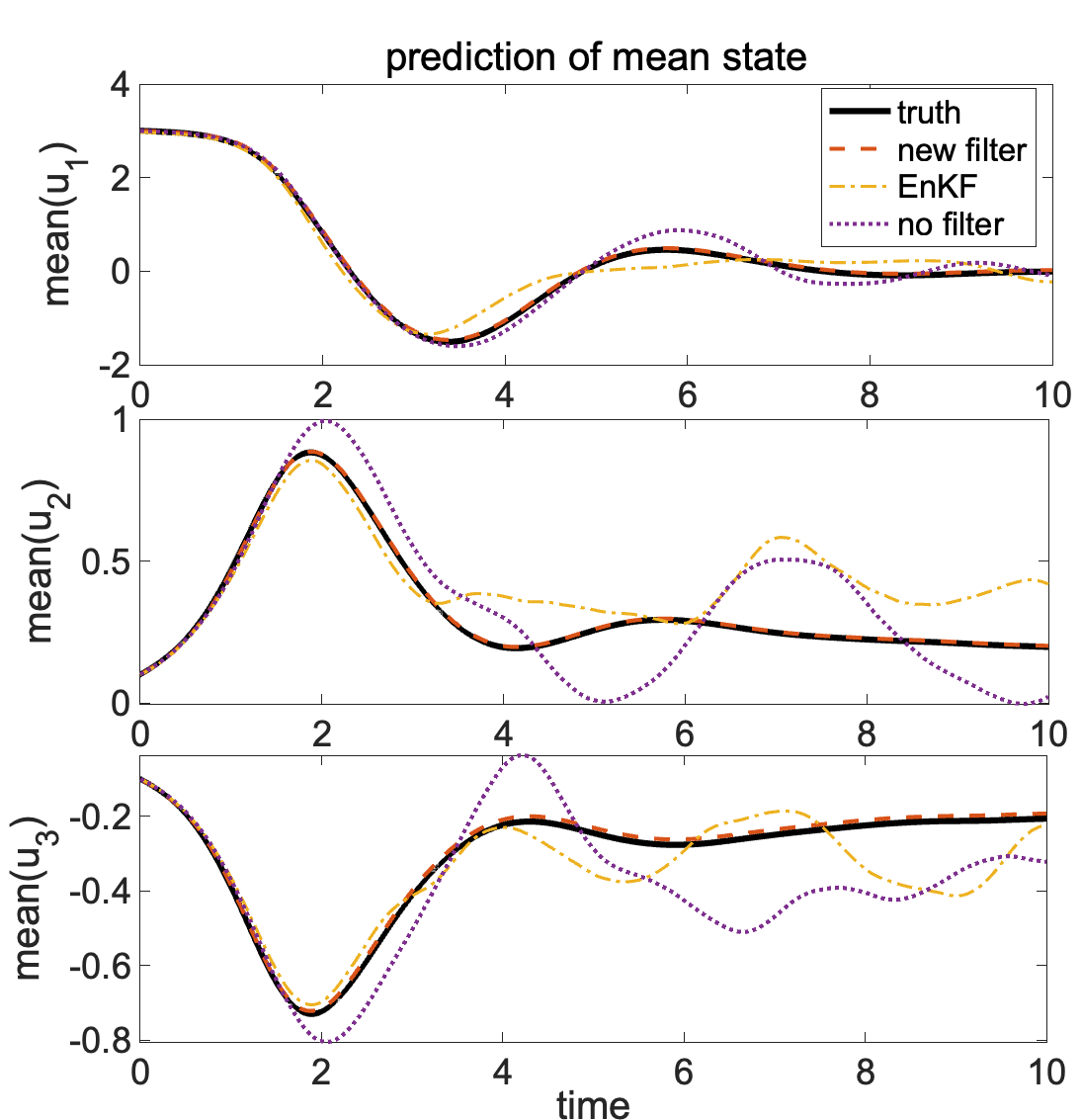}\includegraphics[scale=0.3]{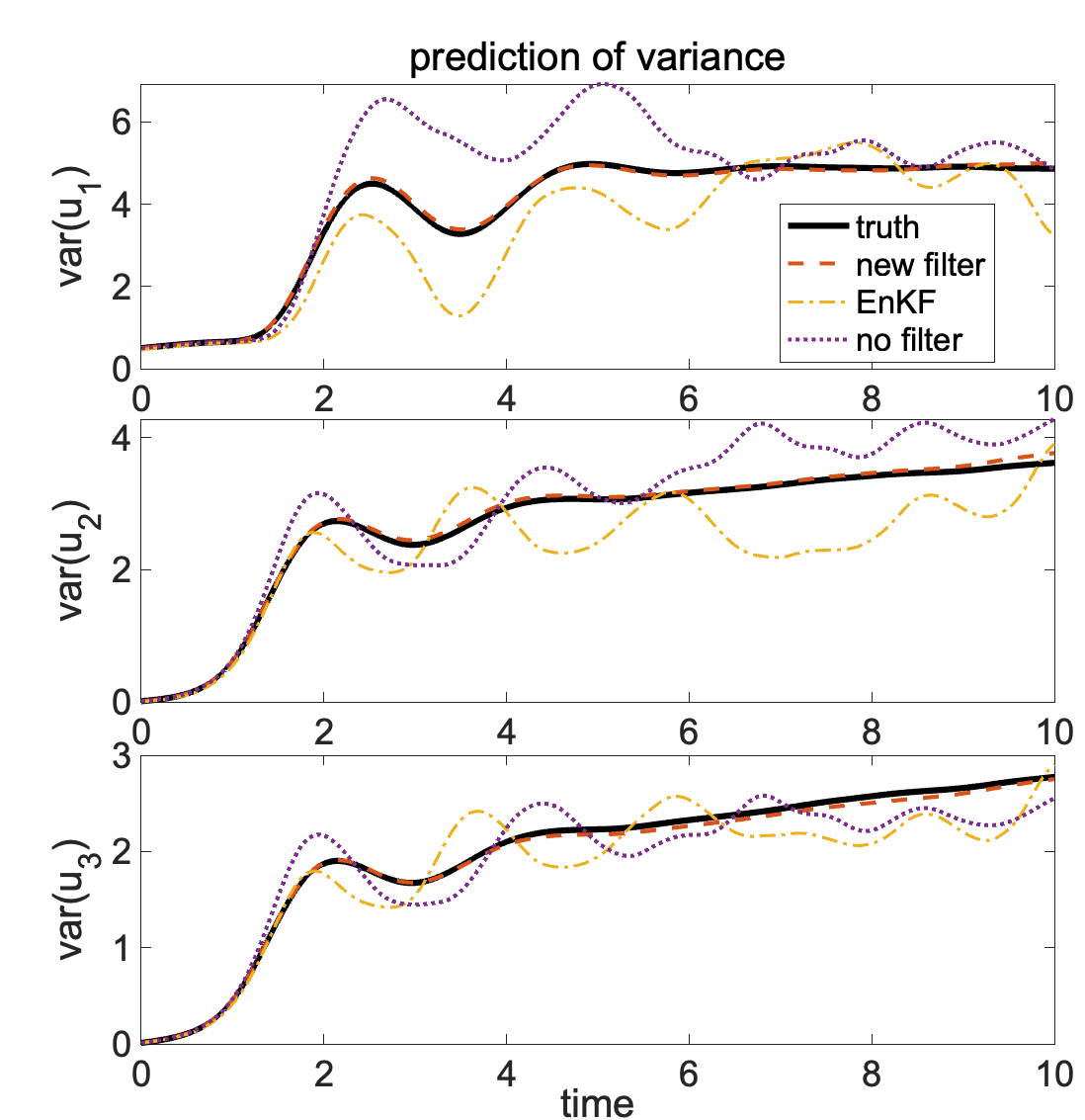}\includegraphics[scale=0.3]{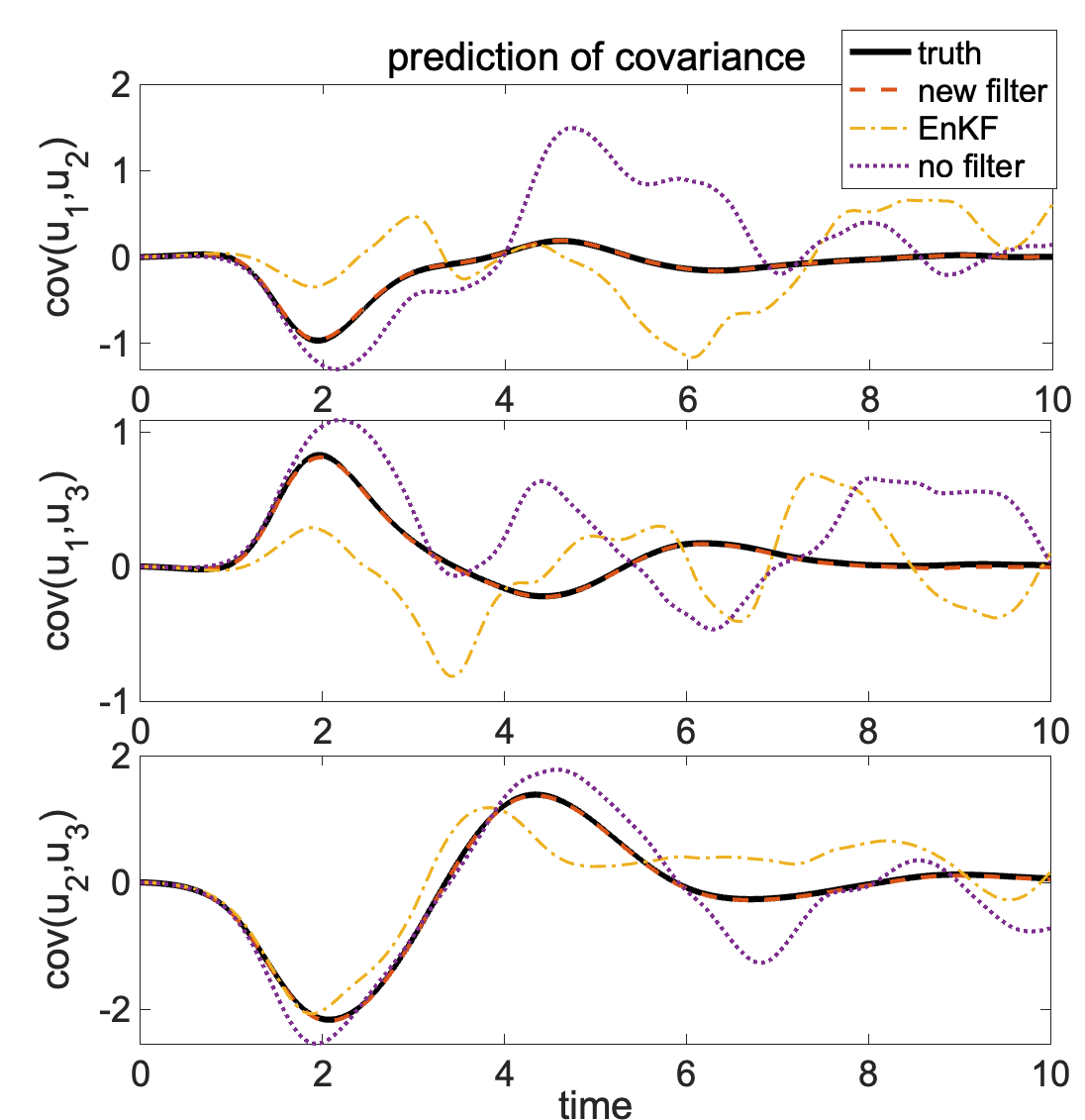}

\caption{Statistical prediction of the mean, variance, and covariance in regime
II of the triad system, with the same setup as in Figure~\ref{fig:Statistical-prediction2}.\protect\label{fig:Statistical-prediction2}}
\end{figure}

\begin{figure}
\includegraphics[scale=0.3]{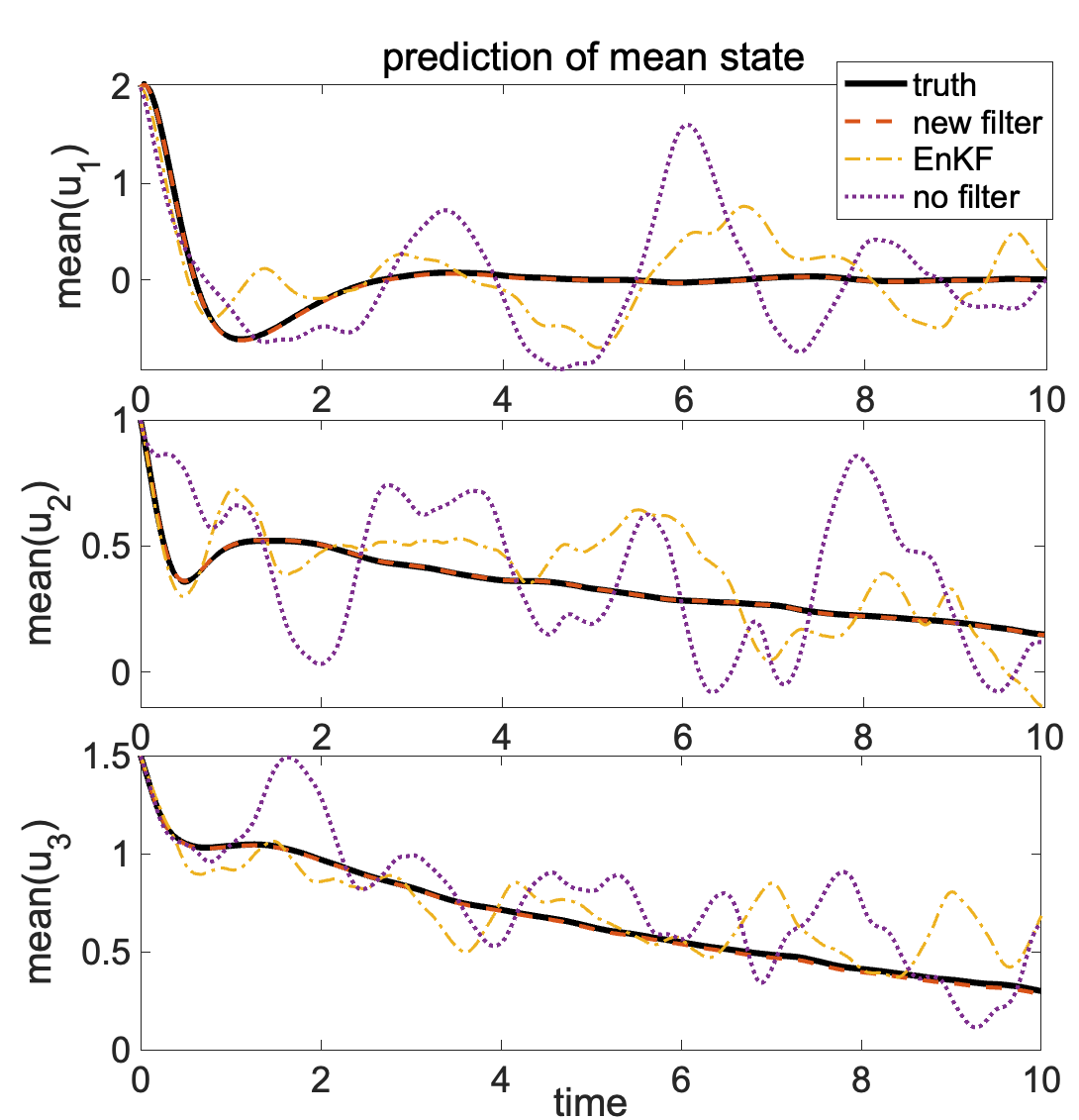}\includegraphics[scale=0.3]{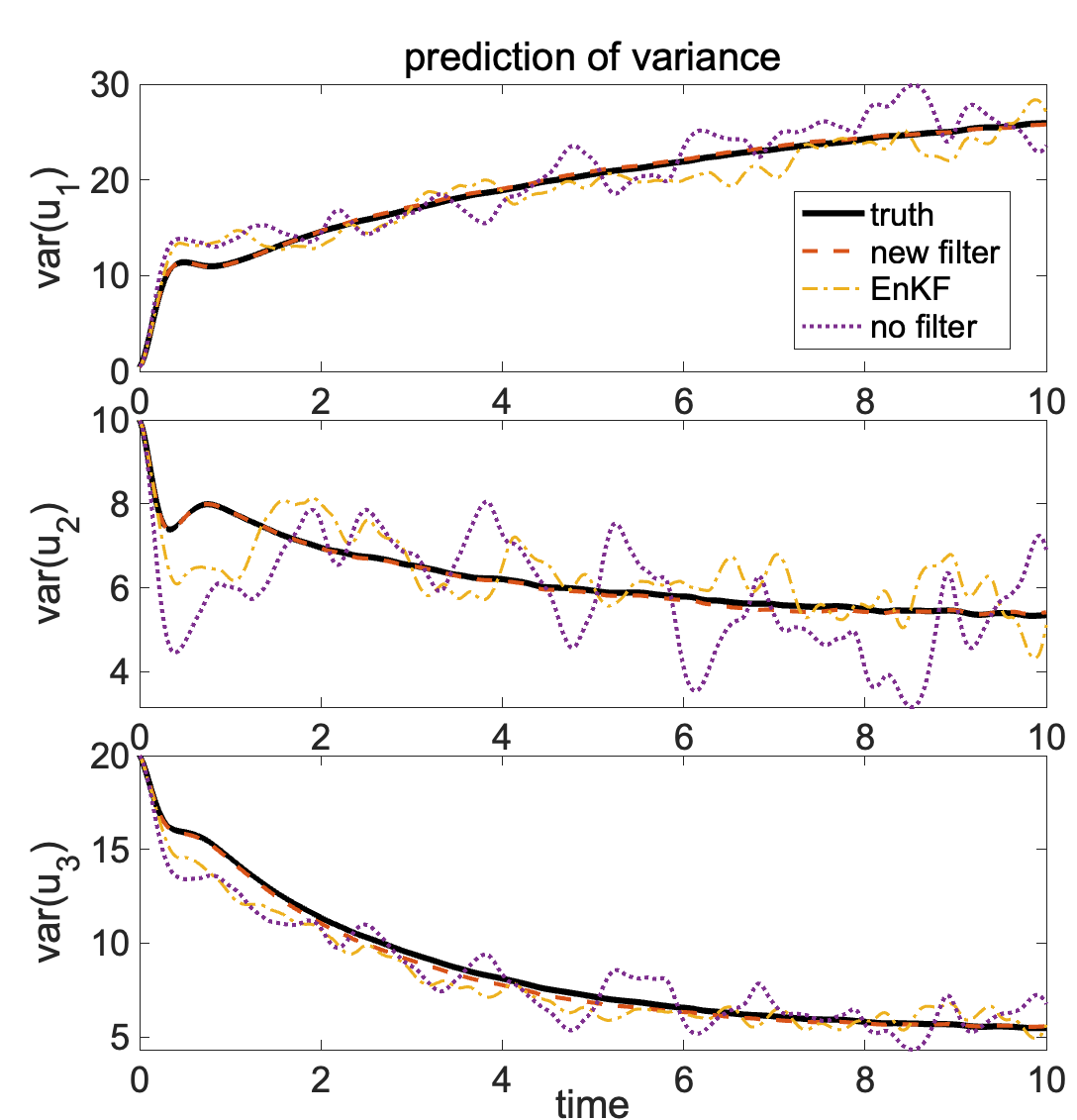}\includegraphics[scale=0.3]{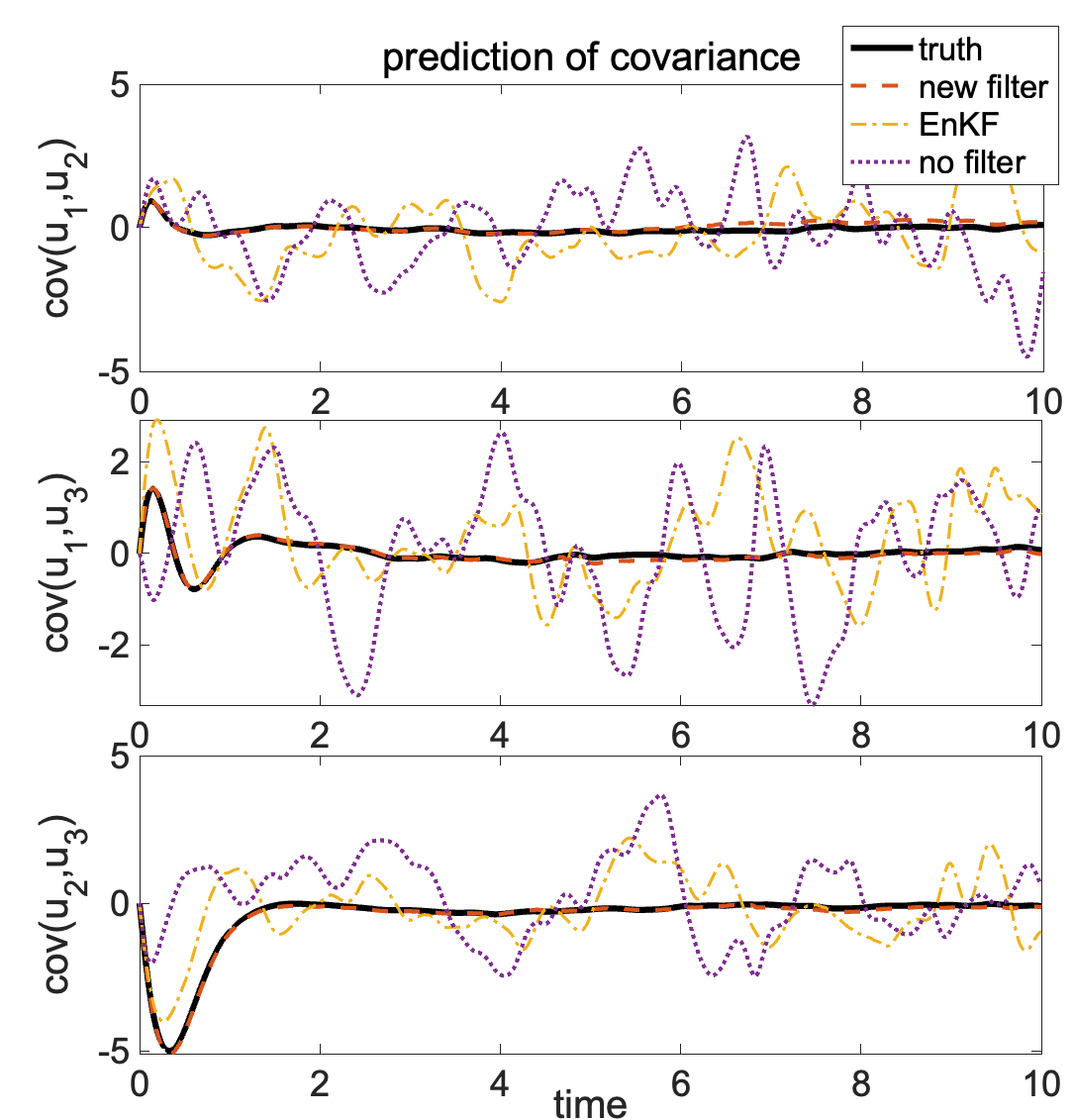}

\caption{Statistical prediction of the mean, variance, and covariance in regime
III of the triad system, with the same setup as in Figure~\ref{fig:Statistical-prediction2}.\protect\label{fig:Statistical-prediction3}}
\end{figure}

In Figure~\ref{fig:Statistical-prediction1}-\ref{fig:Statistical-prediction3},
we plot the model predictions of the mean, variance, and cross-covariance
between the three modes $u_{1},u_{2},u_{3}$ in the three typical
test regimes respectively. The true statistics are compared with the
forecast model without filter \eqref{eq:num_model} and two data assimilation
models. The first model is the standard EnKF \eqref{eq:EnKF} using
only the low-order information and a constant Kalman gain in the filter
update, while higher-order moments are considered according to the
nonlinear observation operators in the new high-order filter model
\eqref{eq:filter_scheme}. The truth is captured by running the original
triad system \eqref{eq:triad} using a very large ensemble size $\mathrm{MC}=1\times10^{5}$.
Only a small ensemble size $N=100$ is used in the model forecasts
for all the tests. Frequent observation data is generated with the
time integration step $\Delta t=1\times10^{-3}$. First in regime
I, the model state will converge to the final near-Gaussian equilibrium
probability distribution. However, this regime demonstrates strong
interactions between the linear operator $L$ and the quadratic nonlinear
operator $B$. This can be illustrated by the persistent positive
growth rate in Figure~\ref{fig:Statistical-forecasts}. The competing
effects lead to strong oscillatory motions between the three modes
indicating frequent exchange of energy between the scales. Non-Gaussian
distributions will also be generated during the transient evolution
of the states. As a result, even starting from accurate initial value
large errors will gradually develop in the direct forecast model without
filter in both the mean and covariance. The low-order EnKF model can
correct the errors a bit from the forecast but still largely deviates
from the true statistical values. In contrast, the high-order filter
maintains the high accuracy in the predictions during the entire evolution
time. In regime II, we focus on the nonlinear effect in the model
driving strong cascade from mode $u_{1}$ to $u_{2},u_{3}$. In this
case, the system is dominated by the nonlinear coupling, and an accurate
characterization of the high-order feedbacks in the statistical equations
will play a central role in achieving good prediction result. As illustrated
in Figure ~\ref{fig:Statistical-prediction2}, the predictions from
the direct forecast model and low-order EnKF quickly diverge from
the truth due to their insufficient sampling of the target probability
distributions. This indicates that the stochastic samples in these
models failed to correctly recover the higher-order moments in the
nonlinear feedback terms in the statistical equations. Again, the
high-order data assimilation scheme keeps stable and accurate predictions
in both the mean and covariance up to the final prediction time. Finally,
regime III sets a most challenging test case containing inherent internal
instability from the linear operator. The nonlinear term then is needed
to introduce the stabilizing effect that needs to be accurately quantified.
Similarly to the other two test cases, we observe that the direct
forecast model and EnKF fail to track the target trajectories of the
mean and covariance with quick divergence to the truth due to the
strong instability quickly amplifying the errors, while the high-order
data assimilation model maintains its high skill against the persistently
instability using only a very small sample size.

\subsection{Prediction of probability distributions and non-Gaussian features}

\begin{figure}
\subfloat[Regime I]{\includegraphics[scale=0.27]{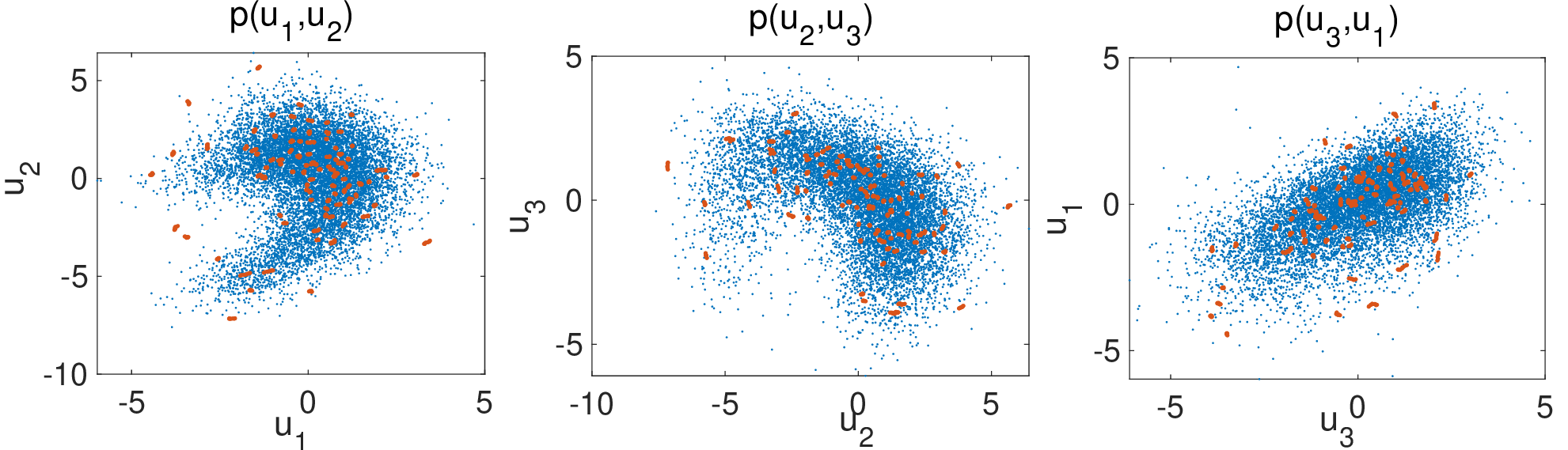}\includegraphics[scale=0.27]{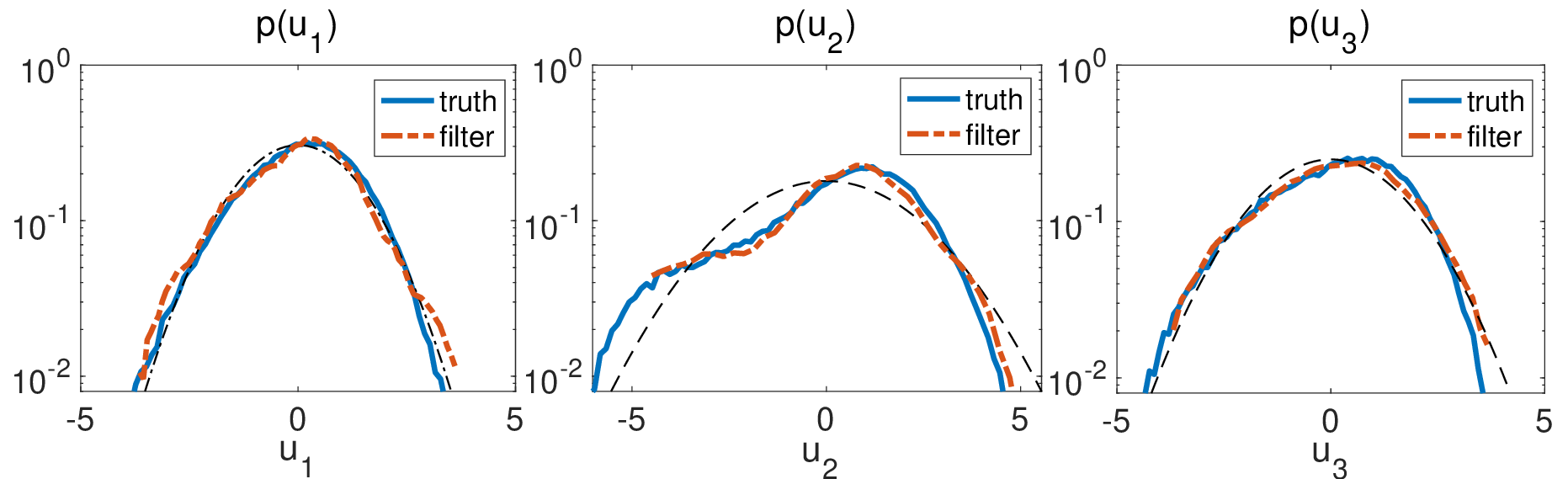}

}

\subfloat[Regime II]{\includegraphics[scale=0.27]{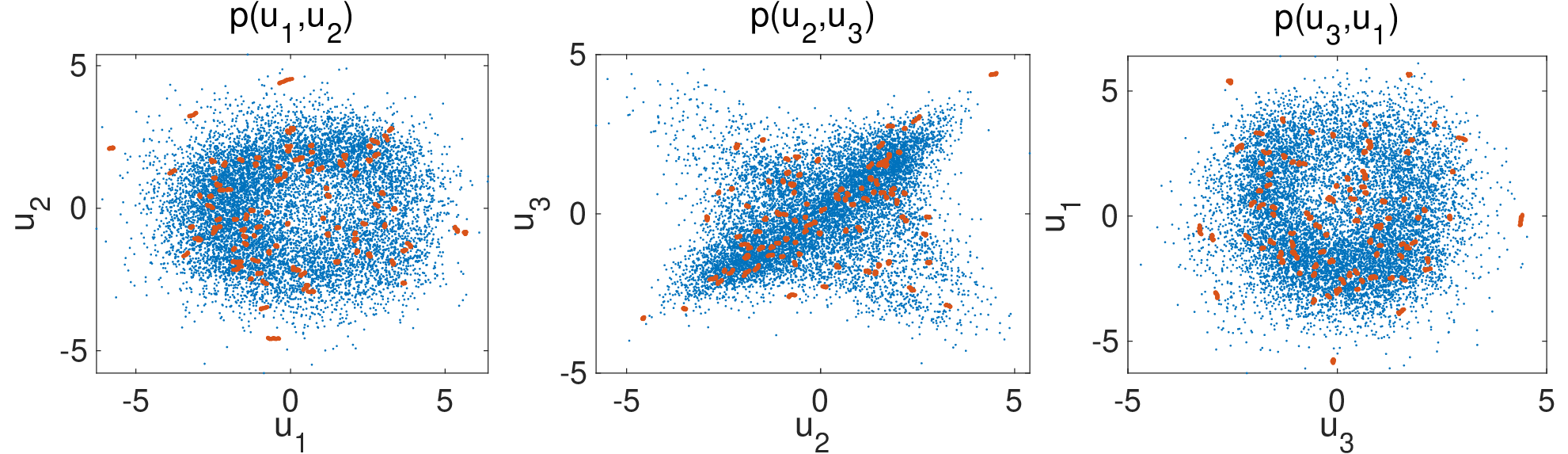}\includegraphics[scale=0.27]{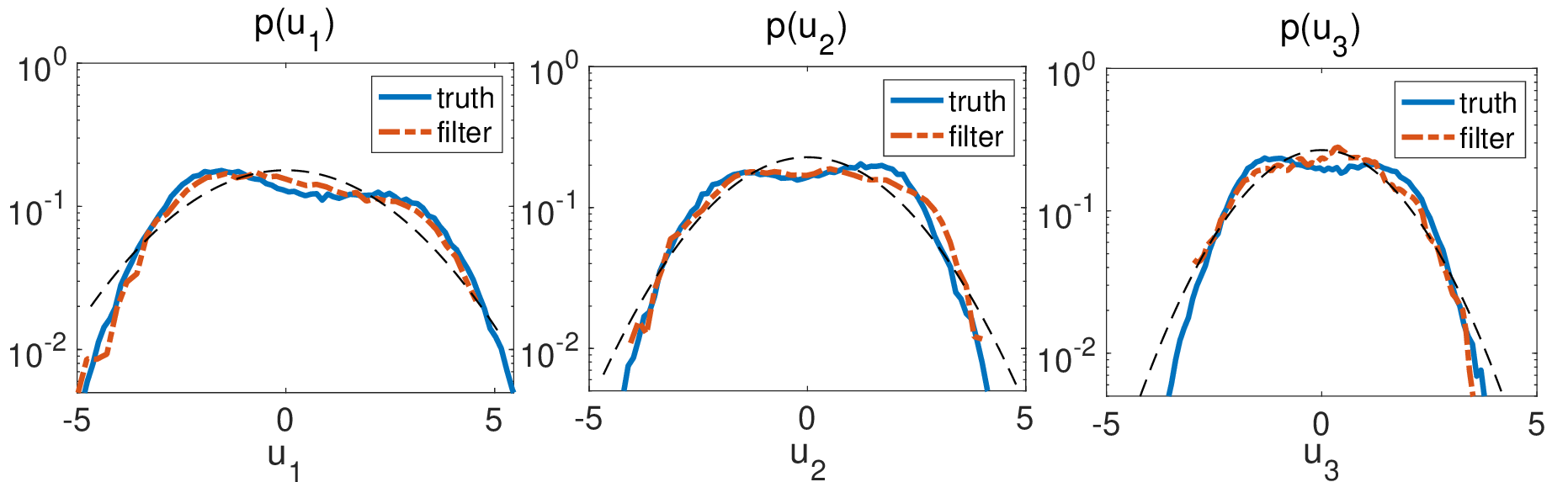}

}

\subfloat[Regime III]{\includegraphics[scale=0.27]{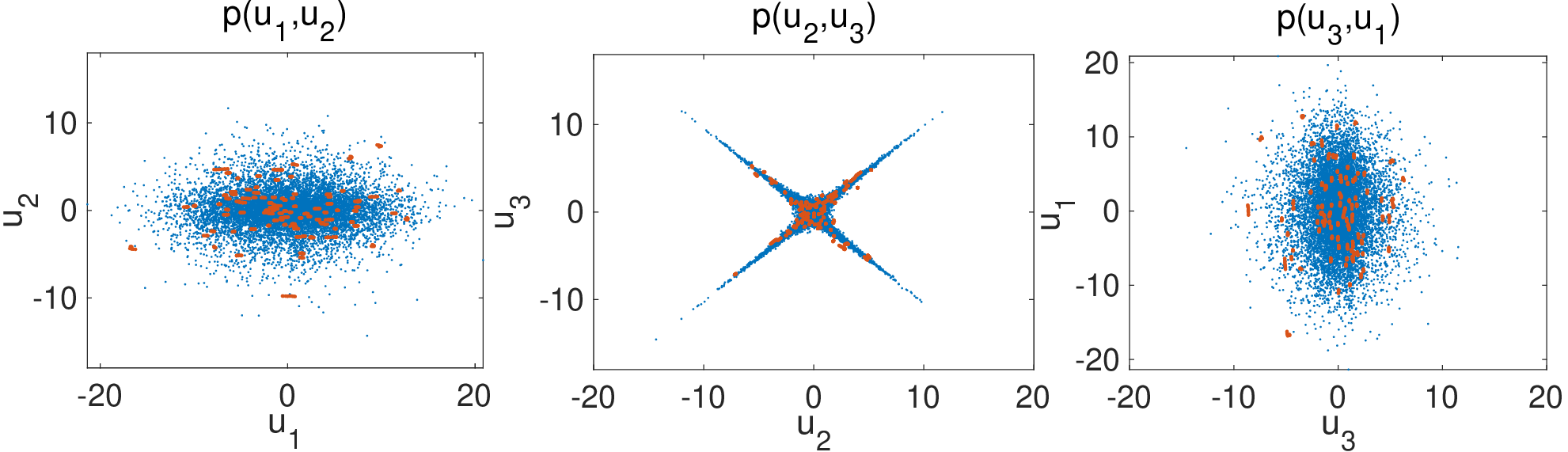}\includegraphics[scale=0.27]{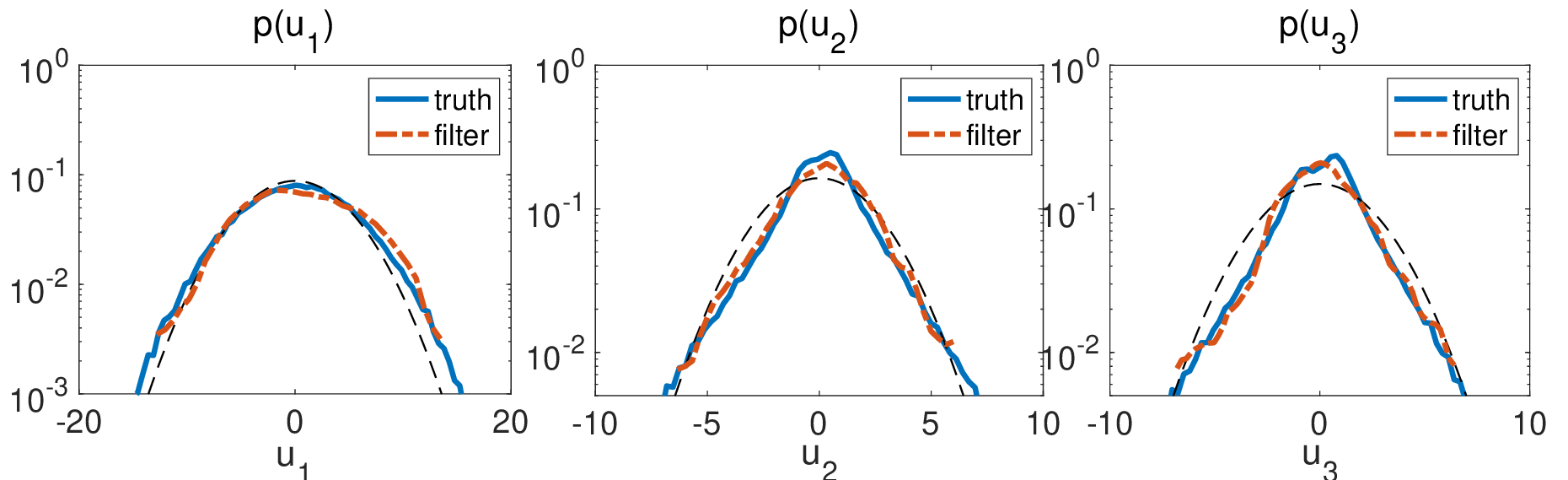}

}

\caption{Probability distributions of the model states at time $t=5$ in the
high-order data assimilation model using $N=100$ samples. The 2D
scatter plots of the truth (blue) are compared with the ensemble filter
solution (red) as well as the 1D marginal distributions. Gaussian
density functions with the same variance are shown in dashed black
lines. \protect\label{fig:Prediction-pdf}}
\end{figure}

The successful prediction of the mean and covariance confirms that
the high-order data assimilation model is able to generate accurate
samples covering the entire spread of the probability distribution
containing essential non-Gaussian statistical features. As a more
detailed illustration of the model forecast skill, we show the sampled
probability distributions in Figure~\ref{fig:Prediction-pdf}. The
projected distributions of the joint states are plotted in scatter
plots together with the marginal PDFs of the three states $u_{1},u_{2},u_{3}$.
First, it can be observed that the typical non-Gaussian features are
generated in all the three test regimes demonstrating highly skewed
or fat-tailed PDFs. These features make important contribution in
the high-order feedback terms in the statistical equations, thus failing
to accurately characterizing their effects in the finite ensemble
approximation will lead to quick divergence that is shown in the statistical
prediction. This explains the large errors and unstable performance
observed in the direct model forecast and EnKF shown in Figure \ref{fig:Statistical-prediction1}-\ref{fig:Statistical-prediction3}
due to the insufficient sampling of these key probability distribution
structures. On the other hand, it shows that the high-order data assimilation
model drives enough samples to the suitable extreme locations so that
the entire non-Gaussian PDFs are will represented among all the test
regimes. This guarantees the high skill of the data assimilation model
to recover the key model statistics including high-order moment information
without requiring a large ensemble size. The uniformly high accuracy
and stability of the new high-order filtering model among all the
three test regimes with distinctive statistical features demonstrate
the universal skill and robustness of the proposed filtering model.

\begin{figure}
\begin{centering}
\includegraphics[scale=0.33]{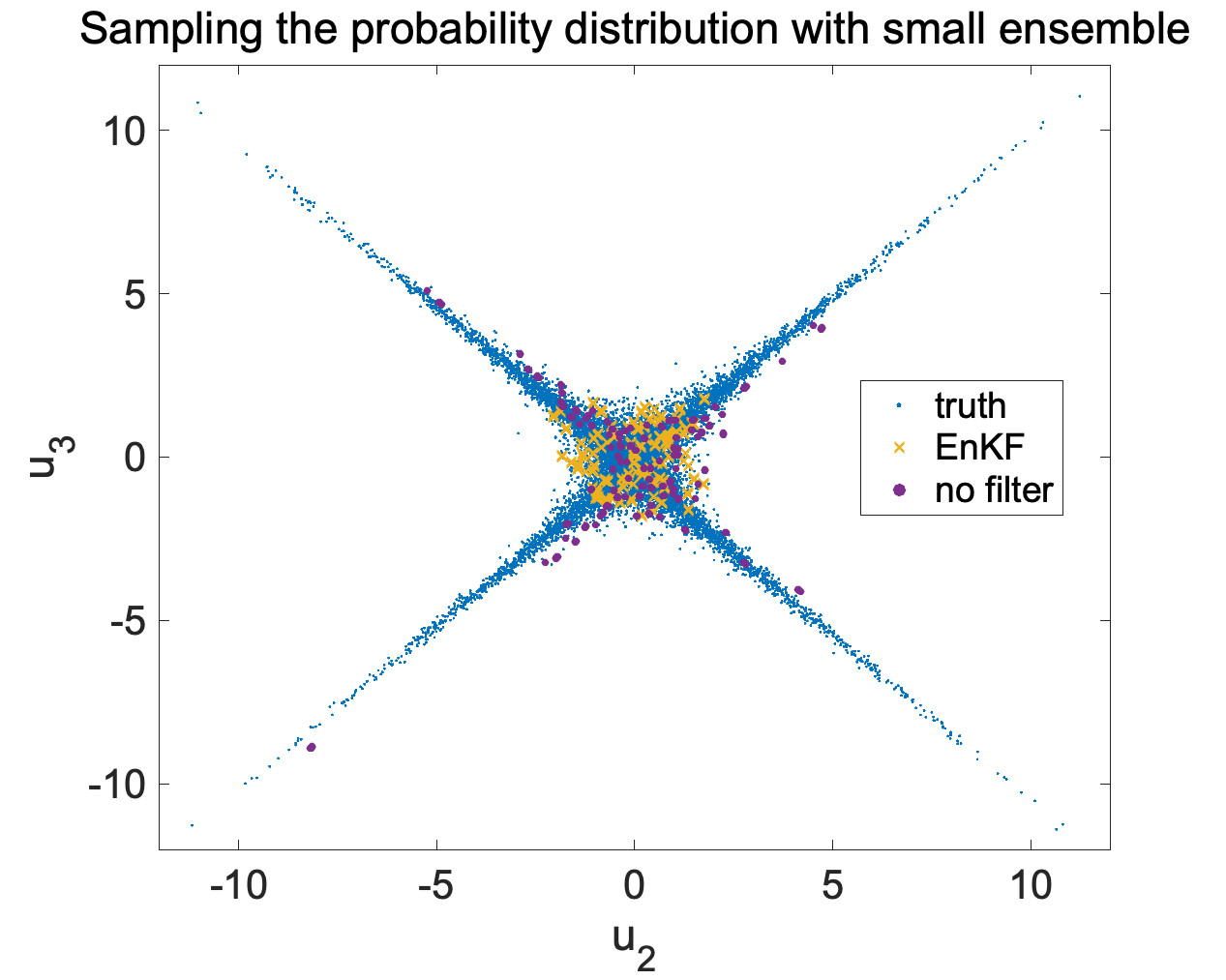}
\par\end{centering}
\caption{Sampling of the target model distribution from a small ensemble forecast
using the direct forecast model without filter \eqref{eq:num_model}
and the standard EnKF \eqref{eq:EnKF}.\protect\label{fig:Sampling}}

\end{figure}

To demonstrate more clearly how the direct forecast model and EnKF
approach fail to reach the accurate time-series predictions in Figure~\ref{fig:Statistical-prediction1}-\ref{fig:Statistical-prediction3},
we show one snapshot of the finite ensemble estimate of the probability
distributions. In particular, Figure~\ref{fig:Sampling} plots the
scatter plots of samples representing the joint PDF of $u_{2}$ and
$u_{3}$ in the most non-Gaussian regime III. It is shown that the
extended four branches of the PDF tail structures are largely missed
in the two models. In the direct forecast model, the small number
of samples cannot sufficiently cover the regions containing extreme
events, and only a few samples can reach the extended wings of the
distribution. The corrections from the EnKF however draw the samples
even closer to a Gaussian distribution rather than reaching the non-Gaussian
features. In contrast, as shown in Figure~\ref{fig:Prediction-pdf},
the new high-order data assimilation model achieves a much better
characterization of the key structures in the probability distribution,
thus guarantees accurate prediction of the statistics. This example
demonstrates the crucial role of accurately sampling the non-Gaussian
PDFs in achieving accurate statistical prediction involving the nonlinear
dynamics.

\begin{figure}
\includegraphics[scale=0.3]{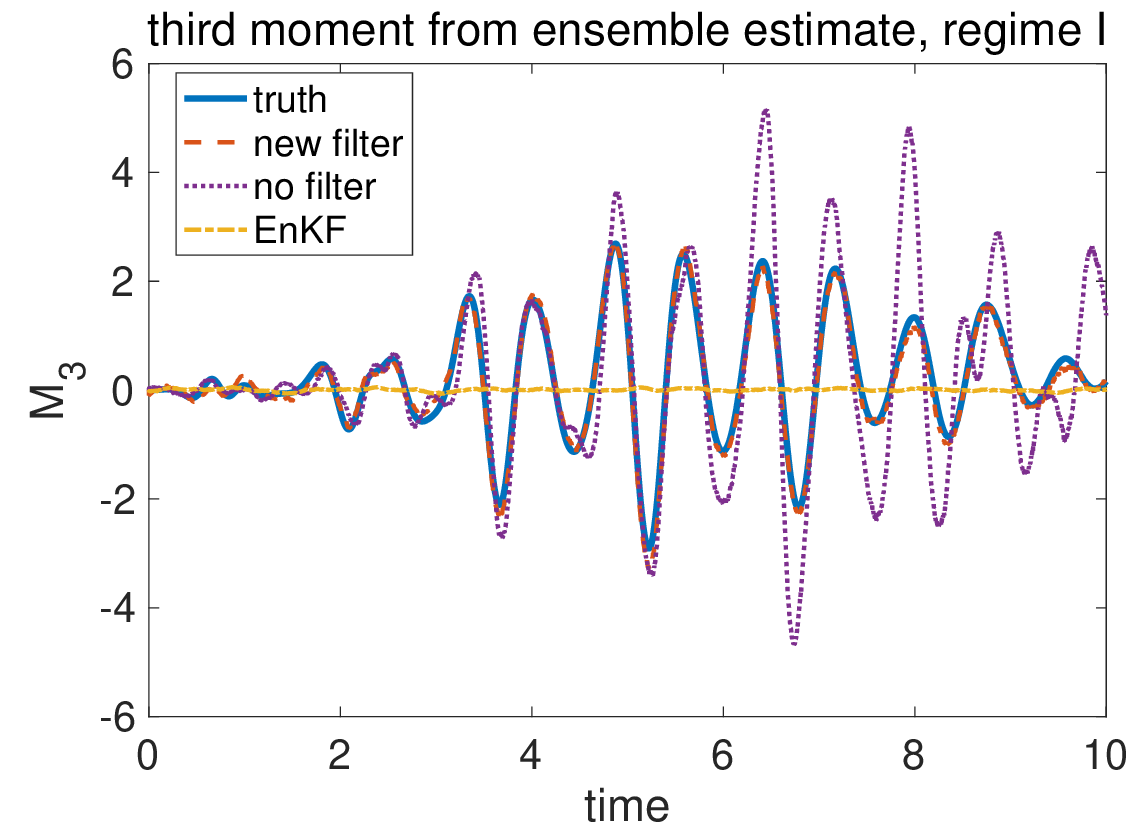}\includegraphics[scale=0.3]{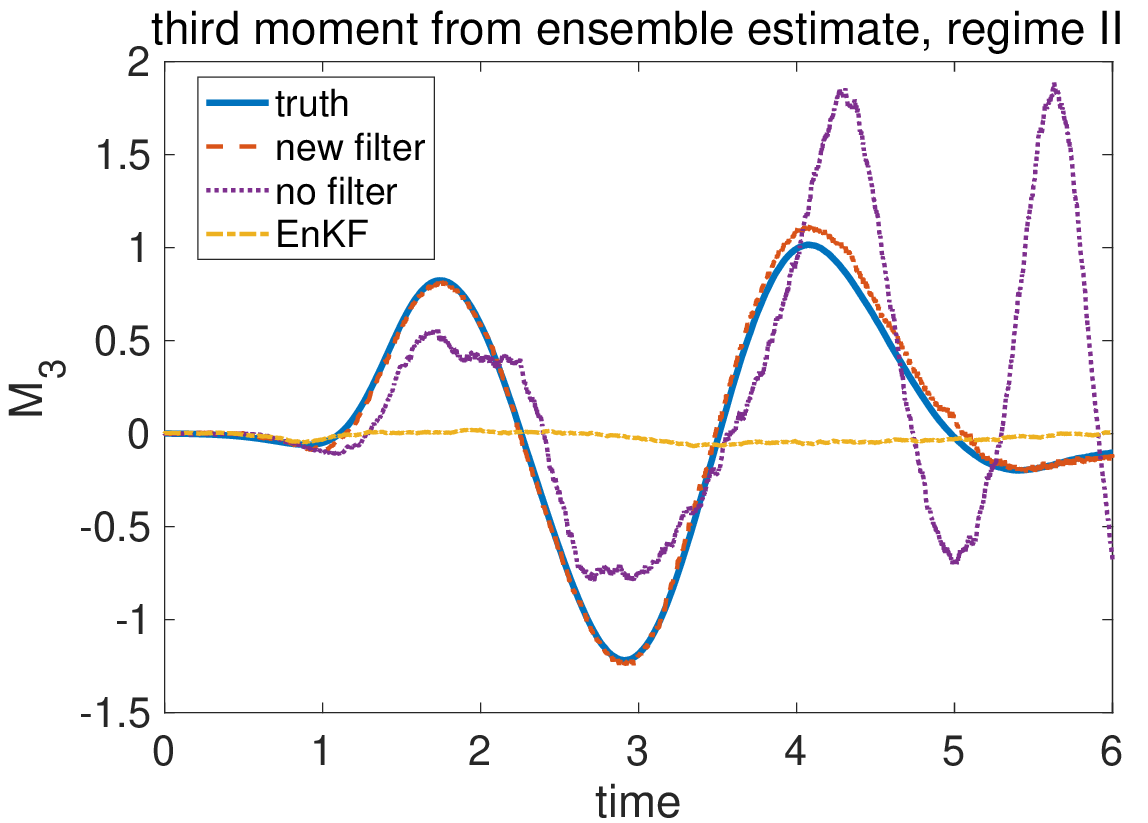}\includegraphics[scale=0.3]{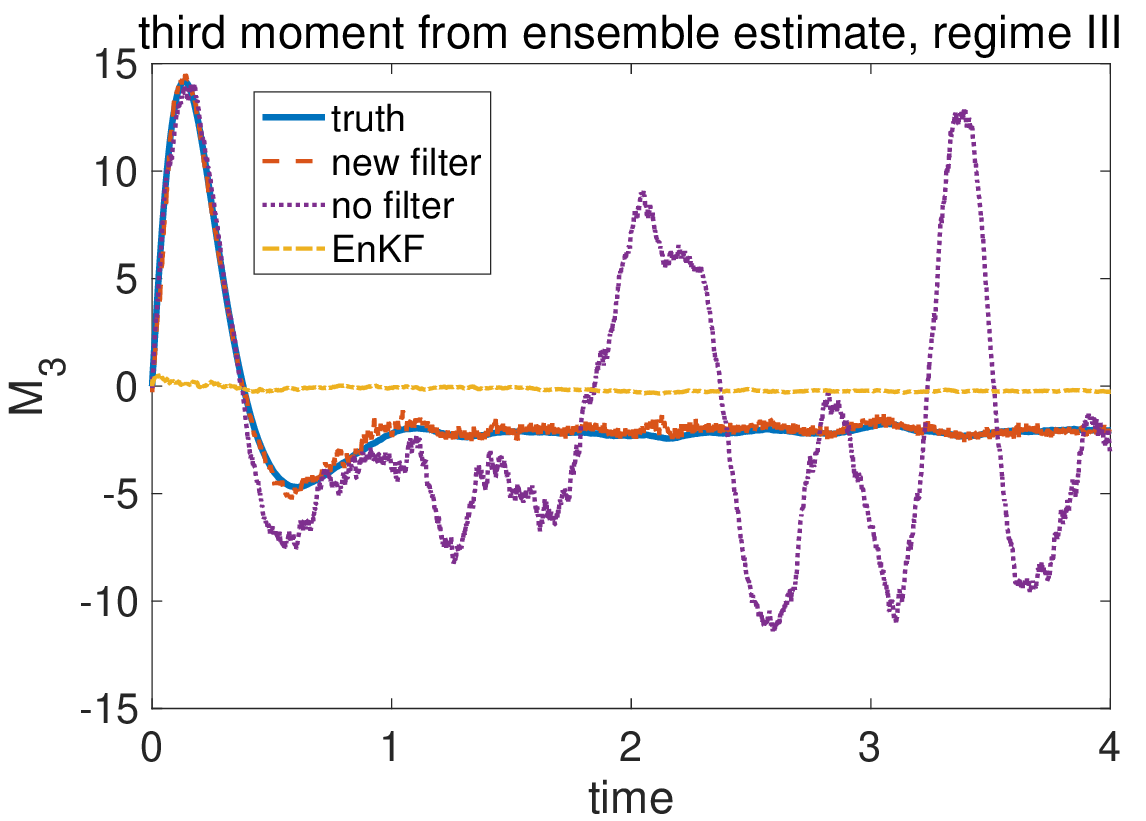}

\caption{Prediction of the third moment $M_{3}=\mathbb{E}\left(u_{1}u_{2}u_{3}\right)$
through the ensemble approximation using the different models in the
three test regimes.\protect\label{fig:Prediction-3mom}}

\end{figure}

As a further illustration of the model prediction of higher-order
moments, Figure~\ref{fig:Prediction-3mom} plots the ensemble recovery
of the third-order moments $M_{3}=\mathbb{E}\left(u_{1}u_{2}u_{3}\right)$
in the three test regimes of the triad model. $M_{3}$ appears in
the dynamical equations for the variances and plays a central role
of balancing the instability from the linear couple terms (see the
explicit statistical equations in \eqref{eq:triad_cor}). First notice
that non-zeros values in $M_{3}$ emerge in all three regimes, showing
the non-negligible role of this high-order feedback term. However,
in the direct forecast model without filter, it can be observed that
the sample estimates of $M_{3}$ are largely missed especially in
the bursts of extreme values. This leads to the final large errors
in the statistical predictions in the mean and covariance shown in
Figure~\ref{fig:Statistical-prediction1}-\ref{fig:Statistical-prediction3}
as well as confirming the biased estimate of the PDF in Figure~\ref{fig:Sampling}.
From the EnKF prediction, on the other hand, near zero values are
assigned to $M_{3}$ from the samples due to the Gaussian nature of
this filter. This is also consistent with the PDF shown in Figure~\ref{fig:Sampling}
and explains the lack of skill in EnKF prediction of the key statistics.
In contrast, the new high-order filter accurately tracks the true
values of $M_{3}$ in the time-series, thus guarantees the successful
predictions of the key statistics.

Finally, we also check the filter performance using different observation
time frequencies $\Delta t_{\mathrm{obs}}$ and using different ensemble
sizes $N$. The total root mean square errors (RMSE) for the predicted
mean and variance are listed in Table~\ref{tab:Prediction-errors-dt}
and \ref{tab:Prediction-errors-N}. For a better quantification of
the sample approximation of the PDFs, we also compute the relative
entropy between the truth and the ensemble estimate. As expected,
using shorter observation time $\Delta t_{\mathrm{obs}}$ and a larger
number of sample $N$ will increase the prediction accuracy, while
the good performance is maintained even with less frequent observations
and an even smaller ensemble size. It further confirms the robust
performance of the high-order data assimilation model to successfully
recover the leading statistics and generate samples that better represent
the key non-Gaussian features in the probability distributions.

\begin{table}
{\footnotesize{}%
\begin{tabular}{ccccccccccccc}
\toprule 
 & \multicolumn{3}{c}{{\footnotesize regime I}} &  & \multicolumn{3}{c}{{\footnotesize regime II}} &  & \multicolumn{3}{c}{{\footnotesize regime III}} & \tabularnewline
\midrule 
{\footnotesize$\Delta t_{\mathrm{obs}}$} & {\footnotesize$0.001$} & {\footnotesize$0.01$} & {\footnotesize$0.05$} & {\footnotesize no filter} & {\footnotesize$0.001$} & {\footnotesize$0.01$} & {\footnotesize$0.05$} & {\footnotesize no filter} & {\footnotesize$0.001$} & {\footnotesize$0.01$} & {\footnotesize$0.05$} & {\footnotesize no filter}\tabularnewline
\midrule
\midrule 
{\footnotesize RMSE in mean} & {\footnotesize 0.0336} & {\footnotesize 0.0584} & {\footnotesize 0.0843} & {\footnotesize 0.3089} & {\footnotesize 0.0216} & {\footnotesize 0.1298} & {\footnotesize 0.1422} & {\footnotesize 0.4799} & {\footnotesize 0.0209} & {\footnotesize 0.1420} & {\footnotesize 0.2403} & {\footnotesize 0.6318}\tabularnewline
\midrule 
{\footnotesize RMSE in variance} & {\footnotesize 0.2631} & {\footnotesize 0.4365} & {\footnotesize 0.5720} & {\footnotesize 1.1285} & {\footnotesize 0.2336} & {\footnotesize 0.6307} & {\footnotesize 0.7680} & {\footnotesize 1.0987} & {\footnotesize 0.5599} & {\footnotesize 1.8973} & {\footnotesize 2.9110} & {\footnotesize 3.9701}\tabularnewline
\midrule 
{\footnotesize Relative entropy} & {\footnotesize 0.2541} & {\footnotesize 0.2105} & {\footnotesize 0.1766} & {\footnotesize 0.4309} & {\footnotesize 0.2408} & {\footnotesize 0.2226} & {\footnotesize 0.2153} & {\footnotesize 0.5475} & {\footnotesize 0.1593} & {\footnotesize 0.2123} & {\footnotesize 0.1409} & {\footnotesize 0.4168}\tabularnewline
\bottomrule
\end{tabular}}{\footnotesize\par}

\caption{Prediction errors with different observation times $\Delta t_{\mathrm{obs}}$
in three test regimes.\protect\label{tab:Prediction-errors-dt}}

\end{table}

\begin{table}
{\footnotesize{}%
\begin{tabular}{ccccccccccccc}
\toprule 
 & \multicolumn{4}{c}{{\footnotesize regime I}} & \multicolumn{3}{c}{{\footnotesize regime II}} &  & \multicolumn{3}{c}{{\footnotesize regime III}} & \tabularnewline
\midrule 
{\footnotesize$N$} & {\footnotesize$50$} & {\footnotesize$100$} & {\footnotesize$200$} & {\footnotesize$500$} & {\footnotesize$50$} & {\footnotesize$100$} & {\footnotesize$200$} & {\footnotesize$500$} & {\footnotesize$50$} & {\footnotesize$100$} & {\footnotesize$200$} & {\footnotesize$500$}\tabularnewline
\midrule
\midrule 
{\footnotesize RMSE in mean} & {\footnotesize 0.0946} & {\footnotesize 0.0336} & {\footnotesize 0.0152} & {\footnotesize 0.0126} & {\footnotesize 0.0679} & {\footnotesize 0.0249} & {\footnotesize 0.0207 } & {\footnotesize 0.0195} & {\footnotesize 0.0368} & {\footnotesize 0.0209} & {\footnotesize 0.0054} & {\footnotesize 0.0024}\tabularnewline
\midrule 
{\footnotesize RMSE in variance} & {\footnotesize 0.7008} & {\footnotesize 0.2631} & {\footnotesize 0.1716} & {\footnotesize 0.1261} & {\footnotesize 0.3729} & {\footnotesize 0.2651} & {\footnotesize 0.1045} & {\footnotesize 0.0773} & {\footnotesize 1.3474} & {\footnotesize 0.5599} & {\footnotesize 0.4406} & {\footnotesize 0.2379}\tabularnewline
\midrule 
{\footnotesize Relative entropy} & {\footnotesize 0.4868} & {\footnotesize 0.2541} & {\footnotesize 0.1906} & {\footnotesize 0.1536} & {\footnotesize 0.4494} & {\footnotesize 0.3463} & {\footnotesize 0.1546} & {\footnotesize 0.0498} & {\footnotesize 0.4375} & {\footnotesize 0.1593} & {\footnotesize 0.0757} & {\footnotesize 0.0461}\tabularnewline
\bottomrule
\end{tabular}}{\footnotesize\par}

\caption{Prediction errors with different ensemble sizes $N$ using the data
assimilation model in three test regimes.\protect\label{tab:Prediction-errors-N}}
\end{table}

\section{Summary\protect\label{sec:Summary}}

In this paper, we developed an explicit high-order data assimilation
framework for effective ensemble prediction of probability distributions
exhibiting highly non-Gaussian statistics. By leveraging observation
data from lower-order statistical moments, the stability and accuracy
of statistical predictions are significantly enhanced using a computational
affordable finite ensemble approach. Specifically, detailed filtering
operators are derived utilizing the explicit quadratic and cubic structures
of the nonlinear coupling terms, resulting in a straightforward numerical
implementation without high computational cost. We performed comprehensive
numerical experiments using an illustrative triad system, which generates
representative turbulent phenomena across different statistical regimes,
to systematically evaluate the skill of the numerical scheme. Inherent
computational barriers for accurate statistical prediction with the
finite ensemble approaches are demonstrated under this simple test
model. Direct numerical comparisons demonstrate that accurately capturing
non-Gaussian distributions is essential for precise statistical prediction
in highly nonlinear dynamics under restricted sampling constraint.
The filtering updates within the proposed data assimilation model
consistently show robust performance in capturing the various types
of non-Gaussian features across multiple tested statistical regimes
requiring only on a small sample size and observation data of leading-moments.
In contrast, traditional ensemble Kalman filter approaches with near-Gaussian
assumptions typically fail to capture such crucial high-order statistical
information.

Future work of this research will include further enhance the proposed
data assimilation strategy by addressing computational challenges
in more complex turbulent systems. In dealing with the high computational
cost of a high-dimensional state, we plan to combine our approach
with high-order moment closure methods incorporating the random batch
approximations \cite{qi2023high,qi2023random}. The additional model
reduction strategies will create practical and computational efficient
algorithms suitable for high-dimensional problems. Immediate applications
of these developments include statistical forecasting in geophysical
turbulence \cite{qi2018rigorous} and modeling viscoelastic fluids
\cite{bao2025deterministic}. Further research may also explore the
extension and validation of the modeling and computational framework
in broader classes of multiscale dynamical systems, potentially enabling
more accurate predictions of extreme events and other complex phenomena.

\section*{Acknowledge}

The research of J.-G. L. is partially supported by the NSF Grant DMS-2106988.
D. Q. is partially supported by ONR Grant N00014-24-1-2192, and NSF
Grant DMS-2407361 and OAC-2232872.

\appendix
\renewcommand\theequation{A\arabic{equation}}
\setcounter{equation}{0}

\section{Proofs of theorems\protect\label{sec:Proofs-of-theorems}}
\begin{proof}
[Proof of Lemma \ref{lem:observation-functions}]By taking partial
derivatives using the explicit expressions in \eqref{eq:obs_func},
we get
\begin{align*}
\partial_{j}H_{k}^{m} & =\gamma_{kjq}z_{q}+\gamma_{kpj}z_{p},\\
\partial_{j}H_{kl}^{v} & =\left(\gamma_{kjq}z_{q}z_{l}+\gamma_{kpj}z_{p}z_{l}+\gamma_{kpq}z_{p}z_{q}\delta_{lj}\right)\\
 & \;+\left(\gamma_{ljq}z_{q}z_{k}+\gamma_{lpj}z_{p}z_{k}+\gamma_{lpq}z_{p}z_{q}\delta_{kj}\right).
\end{align*}
Above, for convenience double appearance of the subindex implies the
summation about the index. Next, multiplying $z_{j}$ and taking the
summation about $j$ yield
\begin{align*}
z^{T}\nabla H_{k}^{m}=z_{j}\partial_{j}H_{k}^{m} & =\gamma_{kjq}z_{j}z_{q}+\gamma_{kpj}z_{p}z_{j}=2H_{k}^{m},\\
z^{T}\nabla H_{kl}^{v}=z_{j}\partial_{j}H_{kl}^{v} & =\left(\gamma_{kjq}z_{l}+\gamma_{ljq}z_{k}\right)z_{j}z_{q}+\left(\gamma_{kpj}z_{l}+\gamma_{lpj}z_{k}\right)z_{p}z_{j}\\
 & \;+\left(\gamma_{kpq}z_{l}+\gamma_{lpq}z_{k}\right)z_{p}z_{q}=3H_{kl}^{v}.
\end{align*}
\end{proof}
\
\begin{proof}
[Proof of Proposition \ref{prop:Kalman_gain}]We can check the solution
\eqref{eq:kalman_algm} by directly putting the expressions back into
the equation \eqref{eq:kalman_gain}. Therefore, for the mean observation
we have
\begin{align*}
\tilde{\mathbb{E}}\left[\left(K^{m}\right)^{T}\nabla H^{m}\right] & =\frac{1}{2}\Gamma_{m}^{-2}\tilde{\mathbb{E}}\left[\left(H^{m}-\bar{H}^{m}\right)\left(Z^{T}\nabla H^{m}\right)\right]\\
 & =\Gamma_{m}^{-2}\tilde{\mathbb{E}}\left[\left(H^{m}-\bar{H}^{m}\right)\left(H^{m}\right)^{T}\right]\\
 & =\Gamma_{m}^{-2}\tilde{\mathbb{E}}\left[\left(H^{m}-\bar{H}^{m}\right)\left(H^{m}-\bar{H}^{m}\right)^{T}\right]=\Gamma_{m}^{-2}C^{H^{m}}.
\end{align*}
Above, the second equality uses the first identity in \eqref{eq:obs_symm}
and the third uses $\tilde{\mathbb{E}}\left[H^{m}-\bar{H}^{m}\right]=0$.
In the same way, we can check the case for variance observation only
with a difference in the coefficient
\begin{align*}
\tilde{\mathbb{E}}\left[\left(K^{v}\right)^{T}\nabla H^{v}\right] & =\frac{1}{3}\Gamma_{v}^{-2}\tilde{\mathbb{E}}\left[\left(H^{v}-\bar{H}^{v}\right)\left(Z^{T}\nabla H^{v}\right)\right]\\
 & =\Gamma_{v}^{-2}\tilde{\mathbb{E}}\left[\left(H^{v}-\bar{H}^{v}\right)\left(H^{v}\right)^{T}\right]\\
 & =\Gamma_{v}^{-2}\tilde{\mathbb{E}}\left[\left(H^{v}-\bar{H}^{v}\right)\left(H^{v}-\bar{H}^{v}\right)^{T}\right]=\Gamma_{v}^{-2}C^{H^{v}}.
\end{align*}
\end{proof}
\
\begin{proof}
[Proof of Proposition \ref{prop:drift}]Using the explicit formula
in \eqref{eq:kalman_algm}, we can compute for the mean observation
case with $\left(\Gamma_{m}^{-2}\right)_{pq}=\gamma_{pq}^{-2}$
\[
\left(K^{m}\Gamma_{m}^{2}K^{mT}\right)_{ij}=\frac{1}{4}Z_{i}Z_{j}\sum_{p,q}\gamma_{pq}^{-2}H_{p}^{m\prime}\left(Z\right)H_{q}^{m\prime}\left(Z\right),
\]
where we denote $H^{\prime}\left(Z\right)=H\left(Z\right)-\bar{H}$.
By taking the divergence on the above identity, we can compute
\begin{align}
\nabla\cdot\left(K^{m}\Gamma_{m}^{2}K^{mT}\right)_{i}= & \sum_{j}\partial_{j}\left(K^{m}\Gamma_{m}^{2}K^{mT}\right)_{ij}=\frac{1}{4}\sum_{j,p,q}\gamma_{pq}^{-2}\partial_{j}\left[Z_{i}Z_{j}H_{p}^{m\prime}\left(Z\right)H_{q}^{m\prime}\left(Z\right)\right]\nonumber \\
= & \frac{1}{4}\sum_{j,p,q}\gamma_{pq}^{-2}\left[\delta_{ij}Z_{j}H_{p}^{m\prime}\left(Z\right)H_{q}^{m\prime}\left(Z\right)+Z_{i}H_{p}^{m\prime}\left(Z\right)H_{q}^{m\prime}\left(Z\right)\right.\nonumber \\
 & \qquad\left.Z_{i}H_{q}^{m\prime}\left(Z\right)Z_{j}\partial_{j}H_{p}^{m}\left(Z\right)+Z_{i}H_{p}^{m\prime}\left(Z\right)Z_{j}\partial_{j}H_{q}^{m}\left(Z\right)\right]\nonumber \\
= & \frac{1}{4}\sum_{p,q}\gamma_{pq}^{-2}Z_{i}\left[H_{p}^{m\prime}\left(Z\right)H_{q}^{m\prime}\left(Z\right)+dH_{p}^{m\prime}\left(Z\right)H_{q}^{m\prime}\left(Z\right)\right]\nonumber \\
 & +\frac{1}{4}\sum_{p,q}\gamma_{pq}^{-2}Z_{i}\left[2H_{p}^{m}\left(Z\right)H_{q}^{m\prime}\left(Z\right)+2H_{q}^{m}\left(Z\right)H_{p}^{m\prime}\left(Z\right)\right]\nonumber \\
= & \frac{5+d}{4}Z_{i}\sum_{p,q}\gamma_{pq}^{-2}H_{p}^{m\prime}\left(Z\right)H_{q}^{m\prime}\left(Z\right)+\frac{1}{2}Z_{i}\sum_{p,q}\gamma_{pq}^{-2}\left[\bar{H}_{p}^{m}H_{q}^{m\prime}\left(Z\right)+\bar{H}_{q}^{m}H_{p}^{m\prime}\left(Z\right)\right].\label{eq:am-1}
\end{align}
The second from last equality above again uses the identity \eqref{eq:obs_symm},
$\sum_{j}Z_{j}\partial_{j}H^{m}=2H^{m}$. Notice that additional term
in the last line above due to the mean term $\bar{H}^{m}$. In the
same way, we can use \eqref{eq:obs_symm} again and find 
\begin{align*}
\left(\nabla\cdot\left(K^{m}\right)^{T}\right)_{q} & =\sum_{j}\partial_{j}\left(\frac{1}{2}ZH^{m\prime}\left(Z\right)^{T}\Gamma_{m}^{-2}\right)_{jq}=\frac{1}{2}\sum_{j,p}\gamma_{pq}^{-2}\partial_{j}\left[Z_{j}H_{p}^{m\prime}\left(Z\right)\right]\\
 & =\frac{1}{2}\sum_{j,p}\gamma_{pq}^{-2}\left[H_{p}^{m\prime}\left(Z\right)+Z_{j}\partial_{j}H_{p}^{m}\left(Z\right)\right]\\
 & =\frac{1}{2}\sum_{p}\gamma_{pq}^{-2}\left[dH_{p}^{m\prime}\left(Z\right)+2H_{p}^{m}\left(Z\right)\right]\\
 & =\frac{d+2}{2}\sum_{p}\gamma_{pq}^{-2}H_{p}^{m\prime}\left(Z\right)+\sum_{p}\gamma_{pq}^{-2}\bar{H}_{p}^{m}.
\end{align*}
Thus the second term in \eqref{eq:drift} gives with the above identity
\begin{align}
\left(K^{m}\Gamma_{m}^{2}\nabla\cdot\left(K^{m}\right)^{T}\right)_{i} & =\frac{1}{2}Z_{i}\left(H^{m\prime}\right)^{T}\nabla\cdot\left(K^{m}\right)^{T}\nonumber \\
 & =\frac{d+2}{4}Z_{i}\sum_{p,q}\gamma_{pq}^{-2}H_{p}^{m\prime}H_{q}^{m\prime}+\frac{1}{2}Z_{i}\sum_{p,q}\gamma_{pq}^{-2}\bar{H}_{p}^{m}H_{q}^{m\prime}.\label{eq:am-2}
\end{align}
Combining the final results in \eqref{eq:am-1} and \eqref{eq:am-2},
we find
\begin{align*}
a_{i}^{m} & =\nabla\cdot\left(K^{m}\Gamma_{m}^{2}K^{mT}\right)_{i}-K_{m}\Gamma_{m}^{2}\nabla\cdot\left(K^{mT}\right)_{i}\\
 & =\frac{3}{4}Z_{i}\sum_{p}\gamma_{pq}^{-2}H_{p}^{m\prime}H_{q}^{m\prime}+\frac{1}{2}Z_{i}\sum_{p,q}\gamma_{pq}^{-2}\bar{H}_{q}^{m}H_{p}^{m\prime}\left(Z\right)\\
 & =\frac{1}{4}Z_{i}\sum_{p}\gamma_{pq}^{-2}H_{p}^{m\prime}\left(3H_{q}^{m\prime}+2\bar{H}_{q}^{m}\right)\\
 & =\frac{1}{4}Z_{i}\sum_{p}\gamma_{pq}^{-2}\left(H_{p}^{m}-\bar{H}_{p}^{m}\right)\left(3H_{q}^{m}-\bar{H}_{q}^{m}\right).
\end{align*}
This gives the expression for the drift term of the mean in \eqref{eq:drift_algm}.
Repeating the same procedure for the variance, we can arrive at the
expression for $a^{v}$ in a similar fashion.
\end{proof}
\
\begin{proof}
[Proof of Proposition \ref{prop:filtering_update}]Using the explicit
expressions derived in \eqref{eq:kalman_algm} and \eqref{eq:drift_algm}
as well as the discrete integration of \eqref{eq:update_stat}, we
can compute the filtering update terms from the observation of the
mean as
\begin{align*}
a^{m}\mathrm{\Delta}t+K_{}^{m}\Delta I_{}^{m}= & \frac{1}{4}Z\left[H^{m}\left(Z\right)-\bar{H}^{m}\right]^{T}\Gamma_{m}^{-2}\left[3H^{m}\left(Z\right)-\bar{H}^{m}\right]\Delta t\\
 & +\frac{1}{2}Z\left[H^{m}\left(Z\right)-\bar{H}^{m}\right]^{T}\Gamma_{m}^{-2}\left[\Delta\bar{u}_{}-\left[H^{m}\left(Z\right)+h_{m}\left(\bar{u}^{N}\right)\right]\Delta t\right]\\
= & \frac{3}{4}ZH^{m\prime}\left(Z\right)^{T}\Gamma_{m}^{-2}H^{m\prime}\left(Z\right)\Delta t+\frac{1}{2}ZH^{m\prime}\left(Z\right)^{T}\Gamma_{m}^{-2}\bar{H}^{m}\Delta t\\
 & +\frac{1}{2}ZH^{m\prime}\left(Z\right)^{T}\Gamma_{m}^{-2}\left[\Delta\bar{u}-\Delta\bar{u}^{N}-H^{m\prime}\left(Z\right)\Delta t\right]\\
= & \frac{1}{4}ZH^{m\prime}\left(Z\right)^{T}\Gamma_{m}^{-2}H^{m\prime}\left(Z\right)\Delta t+\frac{1}{2}ZH^{m\prime}\left(Z\right)^{T}\Gamma_{m}^{-2}\bar{H}^{m}\Delta t+\frac{1}{2}ZH^{m\prime}\left(Z\right)^{T}\Gamma_{m}^{-2}\left(\Delta\bar{u}-\Delta\bar{u}^{N}\right).
\end{align*}
Above, we define the mean and fluctuation terms, $\bar{H}^{m}=\mathbb{E}^{N}\left[H^{m}\left(\tilde{\mathbf{Z}}\right)\right]$
and $H^{m\prime}=H^{m}\left(\tilde{Z}\right)-\bar{H}^{m}$, w.r.t.
the empirical ensemble averages. Similarly, we can compute
\begin{align*}
a^{v}\mathrm{\Delta}t+K^{v}\Delta I^{v}= & \frac{1}{9}Z\left[H^{v}\left(Z\right)-\bar{H}^{v}\right]^{T}\Gamma_{v}^{-2}\left[4H^{v}\left(Z\right)-\bar{H}^{v}\right]\\
 & +\frac{1}{3}Z\left[H^{v}\left(Z\right)-\bar{H}^{v}\right]^{T}\Gamma_{v}^{-2}\left[\Delta R-\left[H^{v}\left(Z\right)+h_{v}\left(\bar{u}^{N},R^{N}\right)\right]\Delta t\right]\\
= & \frac{4}{9}ZH^{v\prime}\left(Z\right)^{T}\Gamma_{v}^{-2}H^{v\prime}\left(Z\right)\Delta t+\frac{1}{3}ZH^{v\prime}\left(Z\right)^{T}\Gamma_{v}^{-2}\bar{H}^{v}\Delta t\\
 & +\frac{1}{3}ZH^{v\prime}\left(Z\right)^{T}\Gamma_{v}^{-2}\left[\Delta R-\Delta R^{N}-H^{v\prime}\left(Z\right)\Delta t\right]\\
= & \frac{1}{9}ZH^{v\prime}\left(Z\right)^{T}\Gamma_{v}^{-2}H^{v\prime}\left(Z\right)\Delta t+\frac{1}{3}ZH^{v\prime}\left(Z\right)^{T}\Gamma_{v}^{-2}\bar{H}^{v}\Delta t+\frac{1}{3}ZH^{v\prime}\left(Z\right)^{T}\Gamma_{v}^{-2}\left(\Delta R-\Delta R^{N}\right).
\end{align*}
\end{proof}
\renewcommand\theequation{B\arabic{equation}}
\setcounter{equation}{0}

\section{Details on the triad system\protect\label{sec:Details-on-triad}}

Here we provide more details on the dynamical and statistical properties
on the triad model \eqref{eq:triad}.

\subsection{A direct link to geophysical turbulent fluid\protect\label{subsec:A-direct-link}}

We can consider the quasi-geostrophic (QG) potential vorticity equation
with forcing and dissipation defined on a two-dimensional periodic
domain $\mathbf{x}\in\left[-\pi,\pi\right]\times\left[-\pi,\pi\right]$
\begin{equation}
\frac{\partial q}{\partial t}+\nabla^{\bot}\psi\cdot\nabla q=\nu\Delta q,\quad\Delta\psi=q.\label{eq:qg}
\end{equation}
Under projection to the Fourier spectral modes $\mathbf{k}=\left(k_{x},k_{y}\right)$
inside a set of finite wavenumber truncation $\mathcal{K}$, the streamfunction
$\psi$ and potential vorticity $q$ can be expressed as
\[
\psi=\sum_{\mathbf{k}\in\mathcal{K}}\hat{\psi}_{\mathbf{k}}e^{i\mathbf{k\cdot x}},\quad q=\sum_{\mathbf{k}\in\mathcal{K}}\left(-\left|\mathbf{k}\right|\right)\hat{\psi}_{\mathbf{k}}e^{i\mathbf{k\cdot x}}.
\]
The QG system \eqref{eq:qg} then can be expressed for each spectral
mode $\hat{\psi}_{\mathbf{k}}$ under the above decomposition as
\[
\frac{\mathrm{d}\hat{\psi}_{\mathbf{k}}}{\mathrm{d}t}+\sum_{\mathbf{k=-m-n}}\frac{\left|\mathbf{n}\right|^{2}}{\left|\mathbf{k}\right|^{2}}\mathbf{m}^{\bot}\cdot\mathbf{n}\hat{\psi}_{\mathbf{m}}\hat{\psi}_{\mathbf{n}}=-\nu\left|\mathbf{k}\right|^{2}\hat{\psi}_{\mathbf{k}}.
\]
Therefore, we have the barotropic triads of three wavenumber components,
$\hat{\psi}_{\mathbf{k}},\hat{\psi}_{\mathbf{m}},\hat{\psi}_{\mathbf{n}}$,
obeying the selecting rule $\mathbf{k+m+n=0}$. Consider an initial
condition in which only these three components of a particular triad
are excited, then these three modes will only interact with each other
while no other modes will get excited due to the particular triad
relations as the system evolves in time. By projecting the above equation
to the active triad modes, we get the dynamical equations for the
selected modes neglecting the forcing and dissipation terms on the
right hand side
\begin{equation}
\frac{\mathrm{d}\hat{\psi}_{\mathbf{k}}}{\mathrm{d}t}+A_{\mathbf{kmn}}\hat{\psi}_{\mathbf{m}}\hat{\psi}_{\mathbf{n}}=0,\quad\mathbf{k+m+n=0},\label{eq:baro_triad}
\end{equation}
where $A_{\mathbf{kmn}}=\frac{\left|\mathbf{n}\right|^{2}}{\left|\mathbf{k}\right|^{2}}\mathbf{m^{\bot}\cdot n}$
is the triad interaction coefficient with the detailed symmetry $A_{kmn}+A_{mnk}+A_{nkm}=0$,
showing the conservation of kinetic energy, 
\[
\frac{\mathrm{d}}{\mathrm{d}t}\left(\left|\mathbf{k}\right|^{2}\left|\hat{\psi}_{\mathbf{k}}\right|^{2}+\left|\mathbf{m}\right|^{2}\left|\hat{\psi}_{\mathbf{m}}\right|^{2}+\left|\mathbf{n}\right|^{2}\left|\hat{\psi}_{\mathbf{n}}\right|^{2}\right)=0.
\]
The typical forward and backward cascades of energy and enstrophy
in turbulent flow are characterized by the triad interactions between
the three models. Hence from the above discussion, in the two-dimensional
QG turbulence, the nonlinear energy transfer is exactly governed by
the barotropic triads the same as \eqref{eq:triad} in the nonlinear
interaction part. More detailed characterization of these coupling
effects with link to geophysical turbulence can be found in \cite{majda2018strategies}.

\subsection{Statistical and dynamical properties of the triad system}

The triad system \eqref{eq:triad} is subject to stochasticity from
the initial state and external forcing.The probability density function
$p\left(\mathbf{u},t\right)$ associated with the triad equations
satisfies the following Fokker-Planck equation
\begin{equation}
\partial_{t}p=-\left(B\left(\mathbf{u},\mathbf{u}\right)+\Lambda\mathbf{u}\right)\cdot\nabla_{\mathbf{u}}p+\sum_{i=1}^{3}\left(d_{i}p+\frac{1}{2}\sigma_{i}^{2}\partial_{u_{i}}^{2}p\right),\label{eq:FP}
\end{equation}
with initial state $p\left(\mathbf{u},t\right)\mid_{t=0}=p_{0}\left(\mathbf{u}\right)$.
While the original triad system \eqref{eq:triad} is nonlinear, the
statistical dynamics \eqref{eq:FP} becomes linear equation for the
smooth PDF $p$. However, in general the explicit solution of the
Fokker-Planck equation is still difficult to achieve due to the nonlinear
interaction terms in the triad system.

\subsubsection{Equilibrium invariant measure with equipartition of energy}

Under a special arrangement about the damping and noise coefficients,
one special solution of a Gaussian invariant measure, $p_{\mathrm{eq}}$,
can be reached at the equilibrium. Assume that the damping operator
$d_{i}$ and random noise forcing $\sigma_{i}$ satisfy the following
relation in each component
\begin{equation}
\sigma_{\mathrm{eq}}^{2}=\frac{\sigma_{1}^{2}}{2d_{1}}=\frac{\sigma_{2}^{2}}{2d_{2}}=\frac{\sigma_{3}^{2}}{2d_{3}}.\label{eq:equi_measure}
\end{equation}
Therefore, a Gaussian invariant measure can be found with equipartition
of energy in each component, that is,
\begin{equation}
p_{\mathrm{eq}}\left(\mathbf{u}\right)=C_{\mathrm{eq}}^{-1}\exp\left(-\frac{1}{2}\sigma_{\mathrm{eq}}^{-2}\left|\mathbf{u}\right|^{2}\right).\label{eq:invar_m}
\end{equation}
Above $\sigma_{\mathrm{eq}}^{2}$ is the equilibrium variance in the
Gaussian invariant distribution $p_{\mathrm{eq}}$ that controls the
variability in each mode. To see this, we can substitute the invariant
measure (\ref{eq:invar_m}) back into the Fokker-Planck equation (\ref{eq:FP}).
It is a special case from the Theorem in \cite{majda2016introduction}.
In the general case with additional external forcing and inhomogeneous
structure, energy is injected into the modes and transferred to each
other due to the nonlinear quadratic interaction through more complicated
mechanism, thus strong nonlinear non-Gaussian statistics with energy
cascade and internal instabilities can be generated. Detailed energy
mechanism and stability for the triad system can be found in \cite{majda2018strategies,majda2019linear}.

\subsubsection{Typical dynamical regimes in the triad system}

Though simple in appearance, the triad system (\ref{eq:triad}) has
representative statistical features including energy cascade between
modes and internal instabilities that can be created in this simple
set-up. A fundamental factor in the triad system is the internal instabilities
that make the mean unstable over various directions in phase space
as is typical for anisotropic fully turbulent systems. Elementary
intuition about energy transfer in such models can be gained by looking
at the special situation with only the nonlinear interactions in (\ref{eq:triad}).
We examine the linear stability of the fixed point, $\bar{\mathbf{u}}=\left(\bar{u}_{1},0,0\right)^{T}$.
Elementary calculations show that the perturbation $\delta u_{1}$
satisfies $\frac{d\delta u_{1}}{dt}=0$ while the perturbations $\delta u_{2},\delta u_{3}$
satisfy the second-order equations
\[
\frac{d^{2}}{dt^{2}}\left(\delta u_{2}\right)=\left(B_{2}B_{3}\bar{u}_{1}^{2}\right)\delta u_{2},\quad\frac{d^{2}}{dt^{2}}\left(\delta u_{3}\right)=\left(B_{2}B_{3}\bar{u}_{1}^{2}\right)\delta u_{3},
\]
so that we find that there is instability in the states $u_{2},u_{3}$
from a non-zero $\bar{u}_{1}$ if $B_{2}B_{3}>0$. Combined with the
energy conservation principle $B_{1}+B_{2}+B_{3}=0,$ we find that
from the initial state $\left(\bar{u}_{1},0,0\right)$
\begin{eqnarray}
 & \mathrm{the\:energy\:in\:}\delta u_{2},\delta u_{3}\mathrm{\:grows\:provided\:that}\label{eq:triad_stability}\\
 & B_{1}\mathrm{\:has}\;\mathrm{the\:opposite\:sign\:with\:}B_{2}\mathrm{\:and\:}B_{3}.\nonumber 
\end{eqnarray}
The elementary analysis in (\ref{eq:triad_stability}) suggests that
we can expect a flow or cascade of energy from $u_{1}$ to $u_{2}$
and $u_{3}$ where it is dissipated provided the interaction coefficient
$B_{1}$ has the opposite sign from $B_{2}$ and $B_{3}$. Then energy
cascades can be induced from the strongly forced unstable energetic
mode to the stable less energetic modes with stronger damping effects. 

\subsection{Moment equations for the triad system}

Here, we provide the detailed moment equations for the mean mean and
covariances of the triad state $\mathbf{u}$. First, the mean state
$\bar{\mathbf{u}}=\left(u_{1},u_{2},u_{3}\right)^{T}$ of the triad
model can be written as
\begin{equation}
\begin{aligned}\mathrm{d}\bar{u}_{1}= & \left[\left(-d_{1}\bar{u}_{1}-\lambda_{3}\bar{u}_{2}+\lambda_{2}\bar{u}_{3}\right)+B_{1}\bar{u}_{2}\bar{u}_{3}+B_{1}\left\langle u_{2}^{\prime}u_{3}^{\prime}\right\rangle \right]\mathrm{d}t,\\
\mathrm{d}\bar{u}_{2}= & \left[\left(\lambda_{3}\bar{u}_{1}-d_{2}\bar{u}_{2}-\lambda_{1}\bar{u}_{3}\right)+B_{2}\bar{u}_{1}\bar{u}_{3}+B_{2}\left\langle u_{1}^{\prime}u_{3}^{\prime}\right\rangle \right]\mathrm{d}t,\\
\mathrm{d}\bar{u}_{3}= & \left[\left(-\lambda_{2}\bar{u}_{1}+\lambda_{1}\bar{u}_{2}-d_{3}\bar{u}_{3}\right)+B_{3}\bar{u}_{1}\bar{u}_{2}+B_{3}\left\langle u_{1}^{\prime}u_{2}^{\prime}\right\rangle \right]\mathrm{d}t,
\end{aligned}
\label{eq:triad_mean}
\end{equation}
where we use $\left\langle \cdot\right\rangle $ to represent the
expectation. Correspondingly, the stochastic fluctuation $\mathbf{u}^{\prime}=\left(u_{1}^{\prime},u_{2}^{\prime},u_{3}^{\prime}\right)^{T}$
of the triad model satisfies the following set of SDEs
\begin{equation}
\begin{aligned}\mathrm{d}u_{1}^{\prime}= & \left[\left(-d_{1}u_{1}^{\prime}-\lambda_{3}u_{2}^{\prime}+\lambda_{2}u_{3}^{\prime}\right)+B_{1}\left(\bar{u}_{2}u_{3}^{\prime}+\bar{u}_{3}u_{2}^{\prime}\right)+B_{1}\left(u_{2}^{\prime}u_{3}^{\prime}-c_{1}\right)\right]\mathrm{d}t+\sigma_{1}\mathrm{d}W_{1},\\
\mathrm{d}u_{2}^{\prime}= & \left[\left(\lambda_{3}u_{1}^{\prime}-d_{2}u_{2}^{\prime}-\lambda_{1}u_{3}^{\prime}\right)+B_{2}\left(\bar{u}_{1}u_{3}^{\prime}+\bar{u}_{3}u_{1}^{\prime}\right)+B_{2}\left(u_{1}^{\prime}u_{3}^{\prime}-c_{2}\right)\right]\mathrm{d}t+\sigma_{2}\mathrm{d}W_{2},\\
\mathrm{d}u_{3}^{\prime}= & \left[\left(-\lambda_{2}u_{1}^{\prime}+\lambda_{1}u_{2}^{\prime}-d_{3}u_{3}^{\prime}\right)+B_{3}\left(\bar{u}_{1}u_{2}^{\prime}+\bar{u}_{2}u_{1}^{\prime}\right)+B_{3}\left(u_{1}^{\prime}u_{2}^{\prime}-c_{3}\right)\right]\mathrm{d}t+\sigma_{3}\mathrm{d}W_{3}.
\end{aligned}
\label{eq:triad_fluc}
\end{equation}
Above, the stochastic equations are coupled with the cross-covariances
$\mathbf{c}=\left(c_{1},c_{2},c_{3}\right)^{T}$ that satisfy the
following statistical equations
\begin{equation}
\begin{aligned}\mathrm{d}c_{1}= & \left[-\left(d_{2}+d_{3}\right)c_{1}+\lambda_{3}c_{2}-\lambda_{2}c_{3}+\lambda_{1}\left(r_{2}-r_{3}\right)\right.\\
 & \left.+\left(B_{2}\bar{u}_{1}r_{3}+B_{3}\bar{u}_{1}r_{2}\right)+\left(B_{2}\bar{u}_{3}c_{2}+B_{3}\bar{u}_{2}c_{3}\right)+\left(B_{2}\left\langle u_{1}^{\prime}u_{3}^{\prime2}\right\rangle +B_{3}\left\langle u_{1}^{\prime}u_{2}^{\prime2}\right\rangle \right)\right]\mathrm{d}t,\\
\mathrm{d}c_{2}= & \left[-\lambda_{3}c_{1}-\left(d_{1}+d_{3}\right)c_{2}+\lambda_{1}c_{3}+\lambda_{2}\left(r_{3}-r_{1}\right)\right.\\
 & \left.+\left(B_{1}\bar{u}_{2}r_{3}+B_{3}\bar{u}_{2}r_{1}\right)+\left(B_{1}\bar{u}_{3}c_{1}+B_{3}\bar{u}_{1}c_{3}\right)+\left(B_{1}\left\langle u_{2}^{\prime}u_{3}^{\prime2}\right\rangle +B_{3}\left\langle u_{1}^{\prime2}u_{2}^{\prime}\right\rangle \right)\right]\mathrm{d}t,\\
\mathrm{d}c_{3}= & \left[\lambda_{2}c_{1}-\lambda_{1}c_{2}-\left(d_{1}+d_{2}\right)c_{3}+\lambda_{3}\left(r_{1}-r_{2}\right)\right.\\
 & \left.\left(B_{1}\bar{u}_{3}r_{2}+B_{2}\bar{u}_{3}r_{1}\right)+\left(B_{1}\bar{u}_{2}c_{1}+B_{2}\bar{u}_{1}c_{2}\right)+\left(B_{1}\left\langle u_{2}^{\prime2}u_{3}^{\prime}\right\rangle +B_{2}\left\langle u_{1}^{\prime2}u_{3}^{\prime}\right\rangle \right)\right]\mathrm{d}t.
\end{aligned}
\label{eq:triad_xcor}
\end{equation}
And correspondingly, the statistical equations for the variances $\mathbf{r}=\left(r_{1},r_{2},r_{3}\right)^{T}$
satisfy the following equations
\begin{equation}
\begin{aligned}\mathrm{d}r_{1}= & 2\left[\left(-d_{1}r_{1}+\lambda_{2}c_{2}-\lambda_{3}c_{3}\right)+B_{1}\left(\bar{u}_{2}c_{2}+\bar{u}_{3}c_{3}\right)+B_{1}\left\langle u_{1}^{\prime}u_{2}^{\prime}u_{3}^{\prime}\right\rangle +\sigma_{1}^{2}\right]\mathrm{d}t,\\
\mathrm{d}r_{2}= & 2\left[\left(-\lambda_{1}c_{1}-d_{2}r_{2}+\lambda_{3}c_{3}\right)+B_{2}\left(\bar{u}_{1}c_{1}+\bar{u}_{3}c_{3}\right)+B_{2}\left\langle u_{1}^{\prime}u_{2}^{\prime}u_{3}^{\prime}\right\rangle +\sigma_{2}^{2}\right]\mathrm{d}t,\\
\mathrm{d}r_{3}= & 2\left[\left(\lambda_{1}c_{1}-\lambda_{2}c_{2}-d_{3}r_{3}\right)+B_{3}\left(\bar{u}_{1}c_{1}+\bar{u}_{2}c_{2}\right)+B_{3}\left\langle u_{1}^{\prime}u_{2}^{\prime}u_{3}^{\prime}\right\rangle +\sigma_{3}^{2}\right]\mathrm{d}t.
\end{aligned}
\label{eq:triad_cor}
\end{equation}

\bibliographystyle{plain}
\bibliography{refs}

\end{document}